\documentclass[a4paper,11pt,reqno,noindent]{amsart}
\usepackage[centertags]{amsmath}
\usepackage{amsfonts,amssymb,amsthm} 
\usepackage[dvips]{graphicx} 
\usepackage{psfrag}
\usepackage[english]{babel}
\usepackage{newlfont}
\usepackage{color}
\usepackage{dsfont}
\usepackage[body={15cm,21.5cm},centering]{geometry} 
\usepackage{fancyhdr}
\usepackage{esint}

\usepackage[active]{srcltx}

\newtheorem{theo}{Theorem}
\newtheorem{lemma}{Lemma}
\newtheorem{prop}{Proposition}
\newtheorem{corollary}{Corollary}

\theoremstyle{definition}
\newtheorem{rem}{Remark}
\newtheorem{defi}{Definition}

\numberwithin{equation}{section}
\numberwithin{lemma}{section} 
\numberwithin{defi}{section} 
\numberwithin{rem}{section} 

\newcommand \dps{\displaystyle }

\newcommand{\loc}{\mathrm{loc}}

\newcommand{\ho}{\mathrm{hom}}

\newcommand{\R}{\mathbb{R}}
\newcommand{\Z}{\mathbb{Z}}
\newcommand{\N}{\mathbb{N}}

\newcommand{\Id}{\text{Id}}
\newcommand{\e}{\varepsilon}

\newcommand{\calH}{\mathcal{H}}
\newcommand{\calF}{\mathcal{F}}

\newcommand{\calP}{\mathcal{P}}
\newcommand{\calE}{\mathcal{E}_{L,T}}
\newcommand{\DD}{\mathrm{D}}
\newcommand{\calL}{\mathcal{L}}

\mathchardef\emptyset="001F

\newcommand{\ext}[1]{\overline{#1}}

\newcommand{\var}[1]{\mathrm{var}\left[#1\right]}
\newcommand{\varO}[1]{\mathrm{var}_0\left[#1\right]}
\newcommand{\cov}[2]{\mathrm{cov}\left[#1;#2\right]}
\newcommand{\Expec}[1]{\left\langle #1 \right\rangle}
\newcommand{\expec}[1]{\left\langle #1 \right\rangle}
\newcommand{\step}[1]{\noindent \textit{Step} #1.}
\newcommand{\supp}[1]{\mathrm{supp} \left(#1 \right)}
\newcommand{\moy}[1]{\left\langle\!\left\langle #1 \right\rangle\!\right\rangle_L}
\newcommand{\osc}[2]{\underset{\dps #1}{\mathrm{osc}} \,#2\,}

\title[Corrector equation in stochastic homogenization]
{Quantitative results on the corrector equation in stochastic homogenization}
\author[A. Gloria \& F. Otto]{Antoine Gloria \& Felix Otto}
\date{\today}
\address[Antoine Gloria]{Universit\'e Libre de Bruxelles (ULB) \\ Brussels, Belgium \\ and Team MEPHYSTO \\  Inria Lille - Nord Europe \\ Villeneuve d'Ascq, France}
\email{agloria@ulb.ac.be}
\address[Felix Otto]{Max-Planck-Institut f\"ur Mathematik in den Naturwissenschaften \\ Leipzig, Germany}
\email{otto@mis.mpg.de}
\begin{document}
\maketitle

\begin{center}
\begin{minipage}{13cm}
\small{
\noindent {\bf Abstract.} 
We derive optimal estimates in stochastic homogenization of linear elliptic equations in divergence form
in dimensions $d\ge 2$. In previous works we studied the model problem of a \emph{discrete} elliptic equation on $\Z^d$. 
Under
the assumption that a spectral gap estimate holds in probability, we proved that there exists a stationary corrector field in dimensions $d>2$ and that the energy density of that corrector behaves as if it had finite range of correlation in terms of the variance of spatial averages --- the latter decays at the rate of the central limit theorem.
In this article we extend these results, and several other estimates, to
the case of a \emph{continuum} linear elliptic equation whose (not necessarily symmetric) coefficient field satisfies a \textit{continuum}
 version of the 
spectral gap estimate. In particular, our results cover the example of Poisson random inclusions.
\vspace{10pt}


\vspace{6pt}
}

\end{minipage}
\end{center}

\bigskip

\section{Introduction}

\noindent We establish quantitative results on the corrector equation for the stochastic homogenization of linear elliptic equations in divergence form, 
when 
the diffusion coefficients satisfy a spectral gap estimate in probability.
Let $\Omega$ be the set of admissible coefficients $A:\R^d\to \R^{d\times d}$ which are measurable and take values into the set of uniformly bounded and elliptic matrices (see Section~\ref{sec:framework} for details).
Consider a probability measure on $\Omega$ (which we call an ensemble) whose expectation is denoted by $\expec{\cdot}$.
Let $D$ be a bounded domain.
Since the seminal contributions of Papanicolaou and Varadhan in \cite{Papanicolaou-Varadhan-79} and of Kozlov 
in \cite{Kozlov-79}, it is known that if the ensemble is stationary and ergodic, then for all $f\in H^{-1}(D)$ and almost every realization of $A$, the weak solution $u_\e \in H^1_0(D)$ of the elliptic equation
\begin{equation*}
-\nabla \cdot A(\frac{\cdot}{\e})\nabla u_\e\,=\,f 
\end{equation*}
weakly converges in $H^1(D)$, as $\e$ vanishes, to the unique weak solution $u_\ho\in H^1_0(D)$ of the deterministic elliptic equation 
\begin{equation*}
-\nabla \cdot A_\ho \nabla u_\ho\,=\,f.
\end{equation*}
The matrix $A_\ho$ is a deterministic and constant elliptic matrix.
As a by-product of the analysis, it is shown that $A_\ho$ is characterized by the formula
\begin{equation}\label{eq:intro-0}
A_\ho \xi\,=\,\expec{A(0)(\xi+\nabla \ext \phi(0))},
\end{equation}
for all $\xi\in \R^d$, where 
$\ext\phi$ is the so-called corrector in direction $\xi$.
It is the unique random field
taking values in $H^1_\loc(\R^d)$ whose realization solves almost surely the \emph{corrector equation}
\begin{equation}\label{eq:intro-1}
-\nabla\cdot A(\xi+\nabla \ext\phi)\,=\,0 
\end{equation}
in the sense of distributions on $\R^d$, such that $\ext \phi(0)=0$ almost surely (at every point $x\in \R^d$ the quantity $\ext\phi(x)$ is almost surely well-defined), and $\nabla \ext\phi$ is stationary and has bounded second moment.
In order to prove the homogenization result, and the existence 
of the corrector field $\ext\phi$, both Papanicolaou \& Varadhan and Kozlov 
rewrite equations \eqref{eq:intro-1} in the probability space $L^2(\Omega)$ (see Section~\ref{sec:framework} for details), where it naturally lives.
In the periodic case --- which can be recast in this setting --- this space is simply $L^2(\mathbb{T})/\R$,  with $\mathbb{T}$ the $d$-dimensional torus. In this case, \eqref{eq:intro-1} reduces
to an elliptic equation on the torus, for which we have the Poincar\'e inequality at our disposal.
In the general ergodic case this nice picture breaks down, and the absence of Poincar\'e's inequality in the infinite-dimensional space
$\Omega$ makes the analysis 
of the corrector equation more subtle. 
To circumvent the lack of coercivity of the elliptic operator in probability, 
these authors add a zero-order term of magnitude $T^{-1}>0$ to the equation, and consider
the unique stationary field with bounded second moment and vanishing expectation $\ext\phi_T$ that solves the modified corrector equation
\begin{equation}\label{eq:intro-2}
T^{-1}\ext\phi_T-\nabla\cdot A(\xi+\nabla \ext\phi_T)\,=\,0 
\end{equation}
in the sense of distributions on $\R^d$ almost surely.
The existence and uniqueness of $\ext\phi_T$ are a direct consequence of the Lax-Milgram theorem. 
In addition,  the a priori estimate
\begin{equation*}
T^{-1}\expec{\ext\phi_T^2(0)}+\expec{|\nabla \ext\phi_T(0)|^2}\,\lesssim\,1 
\end{equation*}
is enough to pass to the limit as $T\uparrow +\infty$ in the equation, and allows one to define
 $\nabla \ext\phi$
as the weak limit of  $\nabla \ext\phi_T$ --- which is a stationary gradient field. Yet 
one loses control of $\expec{\ext\phi_T^2(0)}$, and it is not known whether there exist a stationary random field
$\ext\psi$ such that  $\nabla \ext\psi =\nabla \ext\phi$.

\medskip
\noindent As far as rates are concerned there are only few contributions in the literature.
A first general comment is that ergodicity alone is not enough to obtain convergence rates, so that 
mixing properties have to be assumed on the coefficients $A$.
Besides the optimal estimates in the one-dimensional case by Bourgeat and Piatnitskii \cite{Bourgeat-99},
the first and still unsurpassed contribution in the linear case is due to Yurinski{\u\i} who proved in \cite[(0.10)]{Yurinskii-86} that for $d>2$ and for mixing coefficients with an algebraic decay (not necessarily integrable), there exists $\gamma>0$
such that 
\begin{equation}\label{eq:intro-3}
\expec{|u_\e-u_\ho|^2}\,\lesssim \,\e^{\gamma}.
\end{equation}

\medskip
\noindent The focus of the present paper is not on the homogenization error
$\expec{|u_\e-u_\ho|^2}$, but rather on the corrector field and its decorrelation properties.
As shown in the case of discrete elliptic equations in \cite{Gloria-Neukamm-Otto-11}, this is indeed a first step towards the quantification of the homogenization error.

\medskip
\noindent 
The key ingredient of our analysis is a proxy for Poincar\'e's inequality in probability, in the form of a spectral
gap estimate, which generalizes to the \emph{continuum} setting the estimate
\begin{equation}\label{eq:intro-5}
\var{X} \,\leq \, \sum_{e} \expec{\sup_{\dps a(e)}\left|\frac{\partial X}{\partial a(e)}\right|^2}\var{a}
\end{equation}
we used in the case of a \emph{discrete} elliptic equation, see \cite{Gloria-Otto-09,Gloria-Otto-09b}.
Although this estimate may seem to crucially rely on the fact that there are only countably many random
variables $\{a(e)\}_e$, this is not the case.
In the continuum setting \eqref{eq:intro-5} can indeed be replaced by
\begin{equation}\label{eq:intro-6}
\var{X} \,\lesssim \, \int_{\R^d} \expec{\bigg(\osc{A|_{B(z)}}{X}\bigg)^2}dz,
\end{equation}
where $\osc{A|_{B(z)}}{X}$ denotes the oscillation of $X$ with respect to the restriction of $A$ onto the ball $B(z)=\{z'\,|\,|z-z'|<1\}$ centered at $z$ and of radius 1.
Whereas \eqref{eq:intro-5} holds
for independent and identically distributed coefficients, \eqref{eq:intro-6}  holds for instance
for the Poisson inclusions process.

\medskip
\noindent
With this single ingredient \eqref{eq:intro-6} of probability theory, and in line with the discrete case \cite{Gloria-Otto-09,Gloria-Otto-09b}, we shall prove using linear elliptic PDE theory 
 that all the moments $\expec{|\ext\phi_T(0)|^q}$ ($q>0$) of the modified corrector are bounded for $d>2$ independently of $T$. This implies in particular the existence of a \emph{stationary} corrector, see  Proposition~\ref{prop:main-1} and Corollary~\ref{cor:main-1} below.
Let $\phi'$ denote the adjoint corrector in direction $\xi'$, that is, the corrector associated with the transpose coefficients $A^*$ of $A$.
In terms of quantitative estimates we shall prove for $d>2$ that the variance of smooth averages of the energy density $(\xi'+\nabla \ext\phi')\cdot A(\xi+\nabla \ext\phi)$ of the corrector on balls of radius $L$ decays at the rate $L^{-d}$ of the central limit theorem (as if the energy density had finite range of correlation, which it has not), see Theorem~\ref{th:main-2}.
Last we shall give optimal estimates of the convergence of the gradient $\nabla \ext\phi_T$ of the modified corrector towards the gradient $\nabla \ext\phi$ of the corrector, and of the approximation $\expec{(\xi'+\nabla \ext\phi_T'(0))\cdot A(0)(\xi+\nabla \ext\phi_T(0))}$ of the homogenized coefficients towards the homogenized coefficients $\xi'\cdot A_\ho \xi$, see  Theorem~\ref{th:main-2} and Proposition~\ref{cor:main-2}.

\medskip
\noindent It is worth noticing that our results hold for random diffusion coefficients which are merely measurable.
In particular, what matters for the estimates is only the correlation length of the random coefficient field, not the potentially smaller length scale given by the spatial variations of the coefficients.

\medskip
\noindent Before we conclude this introduction, let us mention the recent contribution by Armstrong and Smart. In \cite{Armstrong-Smart-14}, they develop a quantitative stochastic homogenization theory for (nonlinear) convex integral functionals based on a quantification of the subadditive ergodic theorem for fields with finite range of dependence, and 
get suboptimal algebraic rates of convergence for the Dirichlet problem. 

\bigskip
\noindent Throughout the paper, we make use of the following notation:
\begin{itemize}
\item $d\geq 2$ is the dimension;
\item $\N_0$ denotes the set of non-negative integers, and $\N$ the set of positive integers;
\item $\expec{\cdot}$ is the expectation;
\item $\var{\cdot}$ is the variance associated with the ensemble average;
\item $\cov{\cdot}{\cdot}$ is the covariance associated with the ensemble average;
\item $\lesssim$ and $\gtrsim$ stand for $\leq$ and $\geq$ up to a multiplicative constant which only depends on the dimension $d$, the ellipticity constant $\lambda$ (see \eqref{def:alpha-beta} below), the spectral gap constants $\rho$ and
$\ell$ (see Definition~\ref{lem:var-estim}) if not otherwise stated;
\item when both $\lesssim$ and $\gtrsim$ hold, we simply write $\sim$;
\item we use $\gg $ instead of $\gtrsim$ to specify that the multiplicative constant is large w.~r.~t. 1 (although finite);
\item for all $R>0$, and $z\in \R^d$, $B_R(z):=\{z'\,|\, |z-z'|<R\}$, and $B_R:=B_R(0)$.
\end{itemize}

\section{Main results}\label{sec:main-results}

\subsection{General framework}\label{sec:framework}

In this paragraph we recall some standard results due to Papanicolaou and Varadhan \cite{Papanicolaou-Varadhan-79}.
We start with the definition of the random coefficient field.

\medskip

\noindent 
We let  $\lambda\in(0,1]$ denote an ellipticity constant
which is fixed  throughout the paper, and set
\begin{align}
  \Omega_0:=\Big\{\,A_0\in\R^{d\times d}\,:\,&A_0\text{ is bounded, i.~e. }|A_0\xi| \leq |\xi|\text{ for all $\xi\in\R^d$,}\nonumber \\
  &A_0\text{ is  elliptic, i.~e. } \lambda |\xi|^2\leq \xi\cdot A_0\xi \text{ for all $\xi\in\R^d$}\,\Big\}.\label{def:alpha-beta} 
\end{align}
We equip $\Omega_0$ with the usual
topology of $\R^{d\times d}$. A \textit{coefficient field},
denoted by $A$, is a Lebesgue-measurable function on $\R^d$ taking values in $\Omega_0$. 
We then define 
\begin{equation*}
\Omega:=\{\text{measurable maps }A:\R^d\to \Omega_0\},
\end{equation*}
which we equip with the $\sigma$-algebra $\calF$ 
that makes the evaluations $A\mapsto \int_{\R^d} A_{ij}(x)\chi(x)dx$ measurable for all $i,j\in \{1,\dots,d\}$ and 
all smooth functions $\chi$ with compact support.  This makes  $\calF$ countably generated.

\medskip

\noindent Following the convention in
statistical mechanics, we describe a \textit{random coefficient field} by equipping $(\Omega,\calF)$
with an ensemble $\expec{\cdot}$ (the expected value). Following \cite{Papanicolaou-Varadhan-79}, we shall assume that $\expec{\cdot}$ is stochastically continuous: For all $\delta>0$ and $x\in \R^d$,
$$
\lim_{|h|\downarrow 0}  \expec{\mathds{1}_{\{A\,:\,|A(x+h)-A(x)|>\delta\}}}\,=\,0
$$
We shall always assume that $\expec{\cdot}$ is
\textit{stationary}, i.~e. for all translations $z\in\R^d$ the
coefficient fields $\{\R^d\ni x\mapsto  A(x)\}$ and $\{\R^d\ni
x\mapsto  A(x+z)\}$ have the same joint distribution under $\expec{\cdot}$. 
Let
$\tau_z:\Omega\to\Omega,\; A(\cdot)\mapsto A(\cdot+z)$ denote the shift
by $z$, then $\expec{\cdot}$ is stationary if and only if $\tau_z$ is
$\expec{\cdot}$-preserving for all shifts $z\in\R^d$.
The stochastic continuity assumption ensures that the map $ \R^d \times \Omega \to \Omega, (x,A)\mapsto \tau_x A$
is measurable (where $\R^d$ is equipped with the $\sigma$-algebra of Lebesgue measurable sets).

\medskip

\noindent A random variable is a measurable function on $(\Omega,\calF)$.
We denote by $\calH=L^2(\Omega,\calF,\expec{\cdot})$ the Banach space of square integrable random variables, that is, those random variables $\zeta$ such that $\expec{\zeta^2} <\infty$. This is a Hilbert space for the scalar product $(\zeta,\chi) \mapsto \expec{\zeta \chi}$.
By definition of $\calF$, any random variable of $\calH$ can be approximated by a random variable of $\calH$ that only depends on the value of $A\in \Omega$ on bounded domains.
A \textit{random field} $\tilde\zeta$ is a measurable function on $\R^d\times \Omega$.
To any random variable $\zeta:\Omega\to\R$ we associate
a \textit{$\expec{\cdot}$-stationary extension} $\ext
\zeta:\R^d\times \Omega\to\R$ via $\ext\zeta(x,A):=\zeta(A(\cdot+x))$.
Conversely, we say that a random field is
\textit{$\expec{\cdot}$-stationary} if it can be represented in that
form. 
If $\tilde \zeta$ is a stationary field, then $\tilde \zeta(x,A)=\zeta(\tau_xA)$ for some random variable $\zeta$, so
that for all $x\in \R^d$, $\tilde \zeta(x,\cdot)$ is measurable on $(\Omega,\calF)$ by the measurability of the map $(x,A) \mapsto \tau_xA$ on $\R^d\times \Omega$.
If $\expec{\cdot}$ is stationary, then the ensemble average of
a stationary random field $\ext\zeta$ is independent of $x\in\R^d$;
therefore we simply write $\expec{\ext\zeta}$ instead of $\expec{\ext\zeta(x)}$. 

\medskip

\noindent Stationarity allows one to define a differential calculus on $\calH$.
As shown in \cite[Section~2]{Papanicolaou-Varadhan-79}, since $\expec{\cdot}$ is stochastically continuous, one may define a differential operator $\DD$ on $\calH$ by its components $\DD_i$ in direction $e_i$ for all $i\in \{1,\dots,d\}$ as follows:
\begin{equation*}
\DD_i\zeta (A)\,:=\, \lim_{h\to 0} \frac{\zeta(\tau_{he_i} A)-\zeta(A)}{h} 
\,=\,\lim_{h\to 0} \frac{\ext \zeta (h e_i,A)- \ext \zeta(0,A)}{h}\,=\, \nabla \ext \zeta(0,A).
\end{equation*}
The domain $\calH^1$ of $D$ is closed and dense in $\calH$. It is a Hilbert space for the inner product $(\zeta, \chi) \mapsto \expec{\zeta \chi }+\expec{D\zeta \cdot D\chi}$.

\medskip

\noindent We say that a stationary ensemble is ergodic if the only elements of $\calF$ that are invariant by the shift group $(\tau_z)_{z\in \R^d}$ have probability 0 or 1.
\begin{lemma}[corrector] \cite[Theorem~2]{Papanicolaou-Varadhan-79}\label{lem:corr}
Let $\expec{\cdot}$ be an ergodic stationary ensemble.
Then for all directions $\xi\in\R^d$, $|\xi|=1$,
there exists a unique random field $\ext \phi$ in $H^1_\loc(\R^d,\calH)$ which solves the corrector equation
\begin{equation}\label{eq:corr}
-\nabla\cdot A \left(\xi+ \nabla \ext \phi \right)\,=\,0 
\end{equation}
in the sense of distributions on $\R^d$ and satisfies $\ext \phi(0)=0$, both almost surely,
and such that $\nabla \ext \phi$ is the stationary extension of the field $\nabla \ext \phi(0,\cdot) \in \calH$, with $\expec{\nabla \ext \phi(0,\cdot)}=0$.
In particular, $\expec{|\nabla  \ext\phi(0,\cdot)|^2}\lesssim 1$.
\qed
\end{lemma}
\noindent We also recall the standard definition of the modified corrector:
\begin{lemma}[modified corrector]\cite[Proof of Theorem~2]{Papanicolaou-Varadhan-79}\label{lem:app-corr}
Let $\expec{\cdot}$ be a stationary ensemble.
Then for all $T>0$ and all directions $\xi\in\R^d$, $|\xi|=1$,
there exists a unique random field $\phi_{T}\in \calH^1$ with vanishing expectation, whose stationary extension  $\ext \phi_T$ solves the modified corrector equation
\begin{equation}\label{eq:app-corr}
T^{-1}\ext \phi_{T}-\nabla\cdot A\left(\xi+ \nabla \ext \phi_{T}\right)\,=\,0 
\end{equation}
distributionally on $\R^d$ almost surely, and such that 
 $T^{-1} \expec{\phi_{T}^2}+\expec{{|\DD \phi_{T}|}^2}\lesssim 1$.
\qed
\end{lemma}
\noindent Note that $\ext \phi_T$ is stationary, whereas $\ext \phi$ is not.
\begin{rem}\label{rem:app-corr-eq-proba}
The field $\phi_T$ can be defined as the unique solution in $\calH^1$ of: For all $\zeta \in \calH^1$,
\begin{equation}\label{eq:app-corr-eq-proba}
\expec{T^{-1}\zeta\phi_T+\DD \zeta \cdot A(0)\DD \phi_T}\,=\,-\expec{\DD \zeta\cdot A(0)\xi}.
\end{equation}
\qed
\end{rem}
\begin{rem}\label{rem:adjoint-corr}
If $A$ is replaced by its pointwise transpose $A^*$ in Lemmas~\ref{lem:corr} and~\ref{lem:app-corr}, the associated correctors are called adjoint correctors. For all $\xi'\in \R^d$ and $T>0$, the adjoint corrector $\phi'$ and modified adjoint corrector $\phi'_T$ are suitable solutions
of 
\begin{equation*}
\begin{array}{rcl}
-\nabla \cdot A^* \left(\xi'+\nabla \ext \phi'\right)&=&0,\\
T^{-1}\ext \phi_{T}'-\nabla\cdot A^*\left(\xi'+ \nabla \ext \phi_{T}'\right)&=&0 .
\end{array}
\end{equation*}
\qed
\end{rem}
\begin{defi}[homogenized coefficients]
Let $\expec{\cdot}$ be an ergodic stationary ensemble, let $\xi,\xi'\in \R^d$, and $\ext\phi$ and $\ext\phi'$ be the corrector and adjoint corrector of  Lemma~\ref{lem:corr} and Remark~\ref{rem:adjoint-corr}.
We define the 
homogenized $d\times d$-matrix $A_\ho$ in directions $\xi'$ and $\xi$ by 
\begin{equation}\label{eq:fo-homog}
\xi'\cdot A_\ho\xi\,=\,\expec{(\xi'+\nabla \ext\phi'(0))\cdot A(0)(\xi+\nabla\ext\phi(0))}\,=\,\xi'\cdot \expec{A(0)(\xi+\nabla\ext\phi(0))}.
\end{equation}
\qed
\end{defi}

\subsection{Statement of the main results}\label{sec:stat-results}

To obtain quantitative results, we assume in addition to stationarity and ergodicity that $\expec{\cdot}$ has a spectral gap in the following sense.
\begin{defi}[spectral gap (SG)]\label{lem:var-estim} 
We say that an ensemble $\expec{\cdot}$ satisfies (SG) if there exist
some $\rho>0$ and $\ell<\infty$ such that 
for all measurable functions $X$ on $(\Omega,\calF)$ we have
\begin{equation}\label{eq:var-estim}
\var{X}\;\leq \;\frac{1}{\rho} \left\langle\int_{\mathbb{R}^d}
\left(\osc{A|_{B_\ell(z)}}{X}\right)^2dz\right\rangle,
\end{equation}
where $\osc{A|_{B_\ell(z)}}{X}$ denotes the oscillation
of $X$ with respect to $A$ restricted onto the ball
$B_\ell(z)$ of radius $\ell$ and center at $z\in\mathbb{R}^d$:
\begin{eqnarray}
\left(\osc{A|_{U}}{X}\right)(A)
&=&\left(\sup_{\dps A|_{U}}X\right)(A)
-\left(\inf_{\dps A|_{U}}X\right)(A)\nonumber\\
&=&\sup\left\{X(\tilde A)|\tilde A\in \Omega,\;
\tilde A|_{\mathbb{R}^d\setminus U}=A|_{\mathbb{R}^d\setminus U}\right\}\nonumber\\
&&-
\inf\left\{X(\tilde A)|\tilde A\in\Omega,\;
\tilde A|_{\mathbb{R}^d\setminus U}=A|_{\mathbb{R}^d\setminus U}\right\}.\label{Lc.4}
\end{eqnarray}
Note that for $U\subset\mathbb{R}^d$,
$\osc{A|_{U}}{X}\in [0,+\infty]$ itself is a random variable, which is not necessarily measurable so that
the expectation of the RHS of \eqref{eq:var-estim} is understood as an outer expectation.
\qed
\end{defi}
\noindent By scaling, the choice of the radius 1 is no loss of generality in Definition~\ref{lem:var-estim}.
As the following lemma shows,  (SG) is stronger than ergodicity.
\begin{lemma}\label{lem:SG-ergo}
Let $\expec{\cdot}$ be a stationary ensemble that satisfies (SG) for some $\rho>0$ and $\ell<\infty$.
Then $\expec{\cdot}$ is ergodic.
\qed
\end{lemma}
\noindent The first main result of this paper shows that the variance of smooth averages of the energy density of the modified
corrector on a domain of size $L$ decays according to the central limit theorem scaling $L^{-d}$.
\begin{theo}\label{th:main-1}
Let $\expec{\cdot}$ be 
a stationary ensemble that satisfies (SG),
and let $\ext \phi,\ext \phi'$ and $\ext\phi_T,\ext\phi_T'$ denote the corrector and adjoint corrector, and modified corrector and modified adjoint corrector
for direction $\xi,\xi'\in \R^d$,  $|\xi|=|\xi'|= 1$, and $T>0$, cf. Lemmas~\ref{lem:corr} and~\ref{lem:app-corr}, and Remark~\ref{rem:adjoint-corr}.
We define for all $L>0$ the random matrix $A_{T,L}$ characterized by
\begin{eqnarray*}
\xi'\cdot A_{T,L}\xi:=\int_{\R^d} \big(T^{-1}\ext \phi'_T(x)\ext \phi_T(x)+(\xi'+\nabla\ext \phi'_T(x))\cdot A(x)(\xi+\nabla\ext \phi_T(x))\big)\eta_L(x)\,dx, 
\end{eqnarray*}
where $x\mapsto \eta_L(x)$ is a smooth averaging function on $B_L$ such that $\int_{\R^d}\eta_L(x)dx=1$ and $\sup |\nabla \eta_L| \lesssim L^{-d-1}$. 
Then, for all $T \gg 1$,
\begin{equation}\label{eq:estim-var-hom}
\var{\xi'\cdot A_{T,L}\xi} \,\lesssim\, 
\left\{
\begin{array}{ll}
d=2: & L^{-2} \ln (2+\frac{\sqrt{T}}{L}),\\
d>2: & L^{-d}.
\end{array}
\right.
\end{equation}
In particular, by letting $T\uparrow +\infty$ in  \eqref{eq:estim-var-hom}, the variance estimate
holds for the energy density of the correctors $\ext\phi'$ and $\ext \phi$ themselves for $d>2$.
\qed
\end{theo}
\noindent The main ingredient to the proof of Theorem~\ref{th:main-1} is
of independent interest. It states that all finite stochastic moments
of the modified corrector $\phi_T$ are bounded independently
of $T$ for $d>2$ and grow at most logarithmically in $T$ for
$d=2$. 
\begin{prop}\label{prop:main-1}
Let $\expec{\cdot}$ be 
a stationary ensemble that satisfies (SG),
and let $\phi_T$ denote the modified corrector 
for direction $\xi\in\R^d$, $|\xi|=1$.
Then for all $q\geq 1$ and for all $T\gg 1$ 
\begin{equation}\label{eq:estim-prop}
\expec{|\phi_T|^{q}}^\frac{1}{q} \,\lesssim  \, \left\{
\begin{array}{ll}
d=2: & (\ln T)^{\frac{1}{2}}, \\
d>2: &1,
\end{array}
\right.
\end{equation}
where the multiplicative constant depends on $q$, next to $\lambda$, $\rho$, $\ell$, and $d$.
In addition, for all $q\ge 1$ and for all $R \geq 1$, 
\begin{equation}\label{eq:estim-prop2}
\expec{\Big(\int_{B_R}|\nabla \ext\phi_T(y)|^{2}dy\Big)^{\frac{q}{2}}}^\frac{1}{q} \,\lesssim  \, 1,
\end{equation}
where the multiplicative constant depends on $q$ and $R$, next to $\lambda$,  $\rho$, $\ell$,  and $d$.
\qed
\end{prop}
\begin{rem}
Since (SG) is invariant by transposition of $A$, all the estimates obtained on the modified corrector and on the corrector hold as well for the modified adjoint corrector and the adjoint corrector under the same assumptions on $A$. \qed
\end{rem}

\noindent For $d>2$ we also proved in \cite{Gloria-Otto-14} the corresponding versions of Theorem~\ref{th:main-1} and Proposition~\ref{prop:main-1} for the approximation of the corrector using periodic boundary conditions on cubes of side length $L$. As opposed to the present proof, the proof in \cite{Gloria-Otto-14} does not make use of Green's functions and relies on the De Giorgi-Nash-Moser regularity theory.

\medskip

\noindent
As a direct corollary of Proposition~\ref{prop:main-1} and of Lemma~\ref{lem:app-corr}, we obtain the following existence and
uniqueness result for stationary solutions of the corrector equation \eqref{eq:corr} for $d>2$, which settles a long-standing open question.
\begin{corollary}\label{cor:main-1}
Let $\expec{\cdot}$ be 
a stationary ensemble that satisfies (SG).
Then, for $d>2$ and for all directions $\xi\in\R^d$, $|\xi|=1$, there exists a unique random field $\phi\in \calH^1$
with vanishing expectation whose stationary extension $\ext\phi$ solves the corrector equation 
\begin{equation*}
-\nabla\cdot A(\xi+\nabla \ext \phi)\,=\,0
\end{equation*}
distributionally on $\R^d$ almost surely.
In particular, $\expec{\phi^2+|\DD \phi|^2}\,\lesssim\,1$.
\qed
\end{corollary}

\noindent The proof of this result as a corollary of Proposition~\ref{prop:main-1} is elementary and left to the reader.
Our second main result quantifies the difference between $A_\ho$ and an approximation of $A_\ho$ obtained
using $\nabla \phi_T$ instead of $\nabla \phi$, that we call the systematic error.
\begin{theo}\label{th:main-2}
Let $\expec{\cdot}$ be 
a stationary ensemble that satisfies (SG),
and let $\phi_T,\phi_T'$ denote the modified corrector
and modified adjoint corrector for directions $\xi,\xi'\in\R^d$, respectively, $|\xi|=|\xi'|=1$, and $T>0$.
The approximation $A_T$ of the homogenized matrix $A_\ho$ defined by 
\begin{equation*}
\xi'\cdot A_{T}\xi \,:=\,\expec{(\xi'+\DD \phi_T')\cdot A(0)(\xi+\DD \phi_{T})}\label{eq:def-AT}
\end{equation*}
satisfies for  $T\gg 1$
\begin{equation}\label{eq:estim-err-hom}
|A_\ho-A_T| \,\lesssim\,\left\{
\begin{array}{ll}
d=2: & T^{-1}, \\
d=3: & T^{-\frac{3}{2}},\\
d=4: & T^{-2}\ln T,\\
d>4: & T^{-2}.
\end{array}
\right.
\end{equation}
\qed
\end{theo}
\noindent Note that estimate \eqref{eq:estim-err-hom} saturates at $d=4$. Higher order approximations of $A_\ho$ using the modified correctors $\phi_T$ and extrapolation techniques have been introduced
by Mourrat and the first author in  \cite[Proposition~2]{Gloria-Mourrat-10}. We proved in \cite{Gloria-Neukamm-Otto-14} in the discrete setting that the optimal scaling of the systematic error is $T^{-\frac{d}{2}}$ even beyond $d=4$ and that it can be reached in any dimension for approximations 
of sufficiently high order. We believe that the corresponding continuum version of these estimates also holds true.

\medskip
\noindent
Theorem~\ref{th:main-2} is a direct consequence of the following proposition, which quantifies the convergence
of the gradient of the modified corrector to its weak limit.
\begin{prop}\label{cor:main-2}
Let $\expec{\cdot}$ be
a stationary ensemble that satisfies (SG),
and let $\phi_T$ denote the modified corrector for direction $\xi\in \R^d$,  $|\xi|= 1$, $T>0$, and let 
$\nabla \ext\phi(0)$ denote the weak limit
of $\DD\phi_T$ in $\calH$.
Then for all $T\gg 1$,
\begin{equation}\label{eq:prop1}
\expec{|\DD \phi_T-\nabla \ext\phi(0)|^2} \,\lesssim \, \left\{
\begin{array}{rcl}
d=2&:&T^{-1},\\
d=3&:&T^{-\frac{3}{2}}, \\
d=4&:&T^{-2}\ln T, \\
d>4&:&T^{-2}.
\end{array}
\right. 
\end{equation}
\qed
\end{prop}
\begin{rem}\label{rem:cvg-rate-grad}
For $d>2$, if we denote by $\phi$ the \emph{stationary} corrector of Corollary~\ref{cor:main-1}, 
we also have
\begin{equation}\label{eq:prop1bis}
\expec{(\phi_T-\phi)^2} \,\lesssim \, \left\{
\begin{array}{rcl}
d=3&:&T^{-\frac{1}{2}}, \\
d=4&:&T^{-1}\ln T, \\
d>4&:&T^{-1}.
\end{array}
\right. 
\end{equation}
\qed
\end{rem}

\medskip

\noindent In the case when the coefficients $A$ are symmetric, the operator $\calL=-\DD \cdot A(0) \DD$ defines a quadratic form on $\calH^1$. We denote by $\calL$ its Friedrichs extension on $\calH$ as well.
Since $\calL$ is a self-adjoint non-negative operator, by the spectral theorem, it admits the spectral resolution
\begin{equation}\label{eq:spectral-reso}
\calL\,=\,\int_0^\infty \lambda G(d\lambda).
\end{equation}
We obtain as a by-product of the proof of Proposition~\ref{cor:main-2} the following bounds on the bottom of the spectrum of  $\calL$ projected on $\mathfrak{d}=-\DD \cdot A(0)\xi \in (\calH^1)'$:
\begin{corollary}\label{coro:spectr-bottom}
Let $\expec{\cdot}$ be
a stationary ensemble taking values in the set of symmetric matrices and that satisfies (SG), let $\xi\in \R^d$ with $|\xi|= 1$, and $\mathfrak{d}=-\DD \cdot A(0)\xi$.
Then the spectral resolution $G$ of $\calL=-\DD\cdot A(0)\DD$ satisfies for all $\nu>0$:
\begin{equation}
\expec{\mathfrak{d}G(d\lambda)\mathfrak{d}}([0,\nu])\,\lesssim\,
\left\{
\begin{array}{rcl}
2\leq d<6&:&\nu^{\frac{d}{2}+1},\\
d=6&:&\nu^{4}|\log \nu|,\\
d>6&:&\nu^4.
\end{array}
\right.
\end{equation}
\qed
\end{corollary}
\noindent In the discrete setting we proved  in \cite{Gloria-Neukamm-Otto-14}, using a semi-group approach, that 
$$\expec{\mathfrak{d}G(d\lambda)\mathfrak{d}}([0,\nu])\,\lesssim\,\nu^{\frac{d}{2}+1}$$
holds for all $d\geq 2$.
The method we use here could be pushed forward to prove similar estimates for all $d\geq 6$.

\medskip

\noindent Before we turn to the structure of the proofs, let us comment on the interest of these results.
As in \cite{Gloria-Otto-09,Gloria-Otto-09b}, our main concern here is the approximation of the homogenized coefficients $A_\ho$.
As discussed in \cite{Gloria-Otto-09,Gloria-Otto-09b,Gloria-10} in the discrete setting, the modified correctors $\ext\phi_T$ and $\ext\phi_T'$ can be replaced on some ball $B_L$ by approximations $\ext\phi_{T,R}$ and $\ext\phi_{T,R}'$ computed on a larger ball $B_R$ with homogeneous Dirichlet boundary conditions up to an error of infinite order measured in units of
$\frac{R-L}{\sqrt{T}}$. This holds as well in the continuum setting and we shall consider that we have access to the modified correctors
$\ext\phi_T$ and $\ext\phi_T'$ on $B_L$ in practice.
A natural approximation of $A_\ho$ is then given by
\begin{equation*}
\xi'\cdot \tilde A_{T,L}\xi \,:=\,\int_{B_L}(\xi'+\nabla \ext \phi_{T}'(x))\cdot A(x)(\xi+\nabla \ext\phi_{T}(x))\eta_L(x)dx,
\end{equation*}
where $\eta_L$ is as in Theorem~\ref{th:main-1}.
By stationarity, the error between $\tilde A_{T,L}$ (which is a random variable) and $A_\ho$ satisfies
\begin{equation*}
\expec{(\xi'\cdot \tilde A_{T,L}\xi-\xi'\cdot A_{\ho}\xi)^2}\,=\,\var{\xi'\cdot \tilde A_{T,L}\xi}+\big(\xi'\cdot (A_T-A_\ho)\xi\big)^2.
\end{equation*}
The square root of the first term is called the random error, and the square root of the second term, the systematic error.
The systematic error is estimated in Theorem~\ref{th:main-2}, whereas the random error is estimated
in  Theorem~\ref{th:main-1} as the following remark  shows.
\begin{rem}\label{Remark1}
While it is natural to include the zero-order term 
$T^{-1}\langle \phi_T'\phi_T\rangle$ into the definition of
the energy density, it is not essential for our result. 
Here comes the reason: By a simplified version
of the string of arguments which lead to Theorem~\ref{th:main-1}
we can show that the variance of the zero-order
term is estimated by
$$
\var{\int_{\R^d}\ext\phi'_T(x)\ext\phi_T(x)\eta_L(x)dx}
\;\lesssim\;
\left\{
\begin{array}{ll}
d=2: & \ln T, \\
d>2: & L^{2-d}.
\end{array}
\right. 
$$
This is of higher order than \eqref{eq:estim-var-hom} for $L\lesssim T$.
When approximating $\ext\phi_T$ and $\ext\phi_T'$ on $B_L$ by some $\phi_{T,R}$ and $\phi_{T,R}'$ on a bounded domain $B_R$, one needs
$R-L\gg \sqrt{T}$ for the error due to the artificial boundary conditions to be small.
Taking $R\sim L$, this yields $L\gg \sqrt{T}$, which is compatible with the regime $L\lesssim {T}$.
\qed
\end{rem}

\subsection{The example of the Poisson inclusions process}

\begin{figure}
\centering
\begin{minipage}[b]{0.49\linewidth}
\centering
\includegraphics[scale=.3]{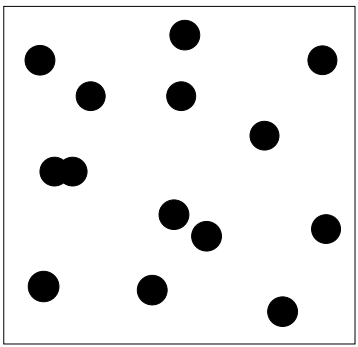}
\caption{Poisson random inclusions}
\label{fig:cont2}
\end{minipage}
\end{figure}

\begin{defi}\label{def:A-poisson}
By the ``Poisson ensemble'' we understand the following probability measure on $\Omega$:
Let the configuration of points $\calP:=\{x_n\}_{n\in \N}$ on $\R^d$
be distributed according to the Poisson
point process with density one. This means the following
\begin{itemize}
\item 
For any two disjoint (Lebesgue measurable) subsets $D$ and $D'$ of $\R^d$ we have that
the configuration of points in $D$ and the configuration of points in $D'$ are independent.
In other words, if $X$ is a function of $\calP$ that depends on $\calP$ only through $\calP|_{D}$
and $X'$ is a function of $\calP$ that depends on $\calP$ only through $\calP|_{D'}$ we have
\begin{equation}\label{L4.23}
\Expec{XX'}_0=\Expec{X}_0\Expec{X'}_0,
\end{equation}
where $\Expec{\cdot}_0$ denotes the expectation w.~r.~t. the Poisson point process.
\item For any (Lebesgue measurable) bounded subset $D$ of $\R^d$, 
the number of points in $D$ is Poisson distributed;
the expected number is given by the Lebesgue measure of $D$.
\end{itemize}
With any realization $\calP=\{x_n\}_{n\in \N}$ of the Poisson point process, 
we associate the coefficient field $A\in\Omega$ (see Figure~\ref{fig:cont2} for a typical realization) via
\begin{equation}\label{eq:A-poisson}
A(x)=\left\{\begin{array}{ccc}
\lambda&\mbox{if}&x\in\bigcup_{n=1}^{\infty}B(x_n)\\
1&\mbox{else}\end{array}\right\}\Id.
\end{equation}
This defines a probability measure $\expec{\cdot}$ on $\Omega$ by ``push-forward'' of $\Expec{\cdot}_0$.
\qed
\end{defi}
\noindent We then have:
\begin{lemma}\label{lem:SG-poisson}
The Poisson ensemble is
stationary and
satisfies (SG) with constants
$\rho=\ell=1$. 
\qed
\end{lemma}
\noindent For a direct proof of Lemma~\ref{lem:SG-poisson} (with suboptimal constants $\rho$ and $\ell$) relying on a martingale decomposition approach, we refer to \cite{Gloria-Otto-14}.
The present version (with optimal constants $\rho=\ell=1$) follows from the well-known Poincar\'e inequality for the Poisson point process: For all measurable functions $X$ of the Poisson point process, we have
\begin{equation}\label{eq:sg-Penrose}
\varO{X}\,\leq\, \int_{\R^d} \expec{(X(\cdot\cup\{x\})-X)^2}_0dx,
\end{equation}
see for instance \cite{Wu-00,Last-Penrose-11}.
For all measurable functions of $A$, we then have
\begin{multline*}
\var{X}=\varO{X\circ A}\,\stackrel{\eqref{eq:sg-Penrose}}{\leq}  \, \int_{\R^d} \expec{(X\circ A(\cdot\cup\{x\})-X\circ A)^2}_0dx
\\
\leq \, \int_{\R^d} \expec{\left(\osc{A|_{B(x)}}{X}\right)^2}_0dx\,=\, \int_{\R^d} \expec{\left(\osc{A|_{B(x)}}{X}\right)^2}dx,
\end{multline*}
where the last two expectations are outer expectations.

\medskip

\noindent General constructions of ensembles $\expec{\cdot}$ from the Poisson point process ensemble $\expec{\cdot}_0$, as well as weighted nonlocal versions of (SG), are discussed in \cite{Duerinckx-Gloria-14b}.


\numberwithin{corollary}{section} 

\subsection{Structure of the proofs and statement of the auxiliary results}\label{sec:structure-proofs}

The proof of Proposition~\ref{prop:main-1} is new and gives optimal scalings in any dimension (contrary to the approach of \cite{Gloria-Otto-09}).
Proposition~\ref{prop:main-1} is a direct consequence of the following two lemmas (and Jensen's inequality in probability).
The first lemma shows that the estimate  \eqref{eq:estim-prop} is a consequence of \eqref{eq:estim-prop2}
for all $q$ large enough.
\begin{lemma}\label{lem:gradphi_to_phi}
Let $\expec{\cdot}$ be 
a stationary ensemble that satisfies (SG),
and let $\phi_T$ denote the modified corrector 
for direction $\xi\in\R^d$, $|\xi|=1$.
Then there exists $\bar q\geq 1$ such that for all $q\geq \bar q$ and for all $T\gg 1$ and $R\gtrsim 1$, 
\begin{equation}\label{eq:estim-lem-gradphi_to_phi}
\expec{|\phi_T|^{2q}}^{\frac{1}{q}} \,\lesssim  \, \expec{\Big(\fint_{B_R}|\nabla \ext\phi_T(y)|^{2}dy\Big)^{q}}^{\frac{1}{q}}\left\{
\begin{array}{ll}
d=2: & \ln T , \\
d>2: &1,
\end{array}
\right.
\end{equation}
where the multiplicative constant depends on $q$, next to $\lambda$, $\rho$, $\ell$, and $d$.
\qed
\end{lemma}
\noindent The second lemma yields  \eqref{eq:estim-prop2}.
\begin{lemma}\label{lem:gradphi_estim}
Let $\expec{\cdot}$ be 
a stationary ensemble that satisfies (SG),
and let $\phi_T$ denote the modified corrector 
for direction $\xi\in\R^d$, $|\xi|=1$.
Then for all $q\geq 1$ and for all $T\gg 1$ and for all $R\gtrsim 1$, 
\begin{equation}\label{eq:estim-lem-gradphi_estim}
\expec{\Big(\fint_{B_R}|\nabla \ext\phi_T(y)|^{2}dy\Big)^{\frac{q}{2}}}^\frac{1}{q}\,\lesssim\,1
\end{equation}
where the multiplicative constants depend on $q$, next to $\lambda$, $\rho$, $\ell$, and $d$.
\qed
\end{lemma}
\begin{rem}
For $d>2$, by Young's inequality, Lemma~\ref{lem:gradphi_estim} is a consequence of Lemma~\ref{lem:gradphi_to_phi} itself
and of the following Caccioppoli inequality in probability for the modified corrector:
For all $q\in \N$, 
\begin{equation}\label{eq:phi-grad-phi}
\expec{\phi_T^{2q}|\DD \phi_T|^2} \lesssim \expec{\phi_T^{2q}},
\end{equation}
as we used  in \cite{Gloria-Otto-09}.
For $d=2$, however, this argument does not provide the optimal power of the logarithm in  \eqref{eq:estim-prop}  nor the optimal scaling in \eqref{eq:estim-prop2} for $d=2$, whence the more subtle approach developed here.
\qed
\end{rem}

\medskip
\noindent In order to prove Lemma~\ref{lem:gradphi_to_phi} we shall apply~(SG)
 to powers of the modified corrector $\phi_T$.
Compared to the discrete setting, we display a significantly simplified proof which avoids the involved induction argument we used in \cite{Gloria-Otto-09}. To this aim, we first derive a ``$q$-version'' of the spectral gap estimate, a continuum analogue of the spectral gap estimate of \cite{Gloria-Neukamm-Otto-14}.
\begin{corollary}[$q$-(SG)]\label{coro:var-estim} 
If $\expec{\cdot}$ satisfies (SG) with constants $\rho>0$ and $\ell<\infty$, then we have for all $q\ge 1$
and all random variables $X$
\begin{equation}\label{eq:cor-var-estim}
\expec{(X-\expec{X})^{2q}}^{\frac{1}{q}}\;\lesssim \;\left\langle\Big(\int_{\mathbb{R}^d}
\Big(\osc{A|_{B_{\tilde \ell}(z)}}{X}\Big)^2dz\Big)^q\right\rangle^{\frac{1}{q}},
\end{equation}
with $\tilde \ell=2\ell$, where the multiplicative constant depends on $q$ and $\rho$.
\qed
\end{corollary}
\noindent 
We note that we lose a factor of two on the radius when passing to the $q$-version 
of (SG) in Corollary~\ref{coro:var-estim}. It is obvious that the original (SG) also holds with radius 2.
From now on, we will use both with radius 2.

\medskip
\noindent
In order to obtain explicit formulas for the oscillation of $\phi_T$,  
we consider an alternative definition for $\ext \phi_T$ that extends the definition
of modified correctors for any $A\in \Omega$ (and not only for almost every $A$). It is as follows.
\begin{lemma}\label{lem:def-phiT-conv}
For all $A\in\Omega$, $T>0$, and $\xi \in \R^d$ with $|\xi|=1$, there exists a unique distributional solution on $\R^d$ of the equation
\begin{equation}\label{eq:app-corr-distribution}
T^{-1}\ext \phi_T-\nabla\cdot A(\xi+\nabla \ext \phi_T)\,=\,0 
\end{equation}
in the class of functions $\chi$ in $H^1_\loc(\R^d)$ such that 
$\limsup_{t\uparrow \infty}\fint_{B_t} (\chi^2+|\nabla \chi|^2)dx \,<\,\infty$.
In addition, this solution satisfies
\begin{equation}\label{eq:uniqueness-cond-refined}
\sup_{z\in \R^d} \int_{B_{\sqrt{T}}(z)} (T^{-1}\ext \phi_T^2+|\nabla \ext \phi_T|^2)dx \,\lesssim \, \sqrt{T}^d.
\end{equation}
\qed
\end{lemma}
\noindent By definition of the $\sigma$-algebra, it is clear that square local averages of $\ext \phi_T$ 
and of $\nabla \ext \phi_T$ are measurable on $(\Omega,\calF)$.
For almost all $A\in \Omega$, the Birkoff ergodic theorem shows that  $\ext \phi_T$ defined in Lemma~\ref{lem:app-corr} satisfies 
$$
\lim_{t\uparrow \infty} \fint_{B_t} (\ext \phi_T^2+|\nabla \ext\phi_T|^2)dx\,=\,\expec{\phi_T^2+|\DD \phi_T|^2}<\infty,
$$
and satisfies \eqref{eq:app-corr-distribution} in the sense of distributions. Hence $\ext \phi_T(\cdot;A)$ coincides 
with the solution of Lemma~\ref{lem:def-phiT-conv} for almost all $A\in \Omega$.

\medskip

\noindent When applying Lemma~\ref{lem:var-estim} to powers of $\ext \phi_T(0;\cdot)$, the sensitivity of $\ext \phi_T(0;A)$ with respect to the coefficients $A$ appears
and needs to be controlled. Our estimates involve Green's functions, whose well-known properties are recalled in the following definition.
\begin{defi}[Green's function]\label{def:Green}
For all $A\in \Omega$ and every $0<T<\infty$, 
there exists a unique function $G_T(x,y;A)\ge 0$ with the following properties
\begin{itemize}
\item Qualitative continuity off the diagonal, that is,
\begin{equation}\label{gu.14}
\{(x,y)\in\mathbb{R}^d\times\mathbb{R}^d|x\not=y\}\ni(x,y)\mapsto G_T(x,y;A)\quad
\mbox{is continuous}.
\end{equation}
\item Upper pointwise bounds on $G_T$:
\begin{equation}\label{eq:ptwise-decay-estim}
G_T(x,y;A)\,\lesssim\,g_T(x-y)\,:=\,\exp(-c\frac{|x-y|}{\sqrt{T}})
\left\{\begin{array}{ccc}
\ln(2+\frac{\sqrt{T}}{|x-y|})&\mbox{for}&d=2\\
|x-y|^{2-d}&\mbox{for}&d>2
\end{array}\right\},
\end{equation}
where here and in the sequel the rate constant $c>0$ in the exponential
is generic and may change from term to term, but only depends on $d$ and $\lambda$.
\item Averaged bounds on $\nabla_x G_T$ and $\nabla_y G_T$:
\begin{eqnarray}
\left(R^{-d}\int_{R<|x-y|<2R}|\nabla_xG_T(x,y;A)|^2dx\right)^\frac{1}{2}&\lesssim& \exp(-c\frac{R}{\sqrt{T}}) R^{1-d},\label{gu.11}\\
\left(R^{-d}\int_{R<|y-x|<2R}|\nabla_yG_T(x,y;A)|^2dy\right)^\frac{1}{2}&\lesssim& \exp(-c\frac{R}{\sqrt{T}}) R^{1-d}.\label{gu.12}
\end{eqnarray}
\item Differential equation:
We note that \eqref{eq:ptwise-decay-estim} and \eqref{gu.11} \& \eqref{gu.12} imply that 
$\mathbb{R}^d\ni x\mapsto(G_T(x,y;A),\nabla_xG_T(x,y;A))$ and 
$\mathbb{R}^d\ni y\mapsto(G_T(x,y;A),\nabla_yG_T(x,y;A))$ are (locally) integrable. Hence even for
discontinuous $A$, we may formulate the requirement
\begin{eqnarray}
T^{-1}G_T-\nabla_x\cdot A(x)\nabla_xG_T=\delta(x-y)&&\mbox{distributionally in}\;\mathbb{R}_x^d,\label{eq:Green}\\
T^{-1}G_T-\nabla_y\cdot A^*(y)\nabla_yG_T=\delta(y-x)&&\mbox{distributionally in}\;\mathbb{R}_y^d,\label{eq:Green-t}
\end{eqnarray}
where $A^*$ denotes the transpose of $A$.
\end{itemize}
We note that the uniqueness statement implies $G_T(x,y;A^*)=G_T(y,x;A)$ so that $G_T$ is symmetric 
when $A$ is symmetric.
\qed
\end{defi}
\noindent
Although these results are well-known, we did not find suitable references dealing with the massive term. We display in the appendix a self-contained proof using as only ingredient the De Giorgi-Nash-Moser theory, and inspired by \cite{Grueter-Widman-82}.

\medskip
\noindent
The following result is of independent interest. 
It quantifies the sensitivity of solutions of linear elliptic PDEs with respect to the coefficient field.
\begin{lemma}\label{lem:diff-phi}
Let $A\in \Omega$, and let $G_T$ and $\phi_T$ be the associated Green function and modified corrector for $T>0$ and $\xi\in\R^d$, $|\xi|=1$. 
Then, for all $x,z\in \R^d$, $R\sim 1$, and $T>0$, we have
\begin{equation}\label{eq:diff-phi-1}
 \osc{A|_{B_R(z)}}{\ext \phi_T(x)} \,\lesssim \, {\ext h}_T(z,x) \left(\int_{B_{3R}(z)}|\nabla \ext\phi_T(y)|^2dy+1 \right)^{\frac{1}{2}},
\end{equation}
where ${\ext h}_T$ is given by
\begin{equation}\label{eq:def-frak-h}
{\ext h}_T(z,x)\,:=\,\left\{
\begin{array}{lcl}
\dps \left(\int_{B_R(z)}|\nabla_y G_{T}(y,x)|^2dy \right)^{\frac{1}{2}} & \text{for} & |z-x| \geq  2R\\
1 &\text{for} & |z-x| < 2R
\end{array}
\right\} \,\lesssim \,1.
\end{equation}
In addition, we also have
\begin{equation}\label{eq:sup-phi-1}
\sup_{\dps A|_{B_R(z)}}\int_{B_R(x)}|\nabla \ext\phi_T(y)|^2dy \,\lesssim \, \int_{B_R(x)}|\nabla \ext\phi_T(y)|^2dy+\int_{B_R(z)}|\nabla \ext\phi_T(y)|^2dy+1.
\end{equation}
\qed
\end{lemma}
\noindent Although this lemma holds for measurable coefficients, we first prove it under an additional smoothness assumption on $A$. This assumption is then removed by an approximation
argument: The pointwise convergence of $\ext\phi_T$ and $G_T$ under the convergence of $A$ follows from the De Giorgi-Nash-Moser theory (in the form of a uniform H\"older estimate).
This is a difference with the discrete setting for which (discrete) gradients of a function $X$ are controlled by the
function $X$ itself and Green functions are not singular --- so that smoothness is not an issue.

\medskip

\noindent As can be already seen on Lemma~\ref{lem:diff-phi}, not only the Green function itself but also
its gradient appears in the estimates.
On the one hand we shall need local estimates which are uniform w.~r.~t. the conductivity function:
\begin{lemma}\label{lem:diff-Green}
Let $A\in\Omega$, and for all $\tilde A\in \Omega$ let $G_T(\cdot,\cdot;\tilde A)$ be the Green function
associated with $\tilde A$, $T>0$.
Then, for all $R\sim 1$, and for all $x,z\in\R^d$ with $|x-z|>R$, we have 
\begin{equation}\label{eq:bd-G(x,e)}
\sup_{\dps \begin{array}{c} \tilde A\in \Omega,\\ \tilde A|_{\R^d\setminus B_R(z)}=A|_{\R^d\setminus B_R(z)}\end{array}} \int_{B_R(z)}|\nabla_yG_T(y,x;\tilde A)|^2dy \,\lesssim \,\int_{B_R(z)}|\nabla_yG_T(y,x;A)|^2dy.
\end{equation}
\qed
\end{lemma}
\noindent On the other hand we shall make use of both integrated and pointwise estimates
on the gradient of the Green function: optimal quenched but integrated or annealed but pointwise control with an exponent $2p$ slightly larger than $2$ --- Meyers' type estimates --- and a suboptimal but quenched and pointwise control.
\begin{lemma}[optimal quenched integrated estimates of gradients]\label{lem:int-grad}
Let $A\in\Omega$ and $G_T$ be its associated Green function, $T>0$. 
Then, there exists $\bar{p}>1$ depending only on $\lambda$, and $d$ such that for all $\bar{p}\geq  p\geq 1$ and $R>0$, we have
\begin{eqnarray}
\left(R^{-d}\int_{R< |y|\leq 2R}|\nabla_yG_T(y,0) |^{2p}dy\right)^{\frac{1}{2p}}&\lesssim& R^{1-d} \exp(-c\frac{R}{\sqrt{T}})\label{eq:int-grad},
\\
\left(R^{-2d}\int_{B_R}\int_{8R< |y|\leq 16R}|\nabla\nabla G_T(y,x) |^{2p}dydx\right)^{\frac{1}{2p}} &
\lesssim & R^{-d} \exp(-c\frac{R}{\sqrt{T}}),\label{eq:int-grad2}
\end{eqnarray}
where $\nabla \nabla$ denotes the mixed second gradient.
\qed
\end{lemma}
\noindent For the proof of \eqref{eq:int-grad} in Lemma~\ref{lem:int-grad}, we refer the reader to 
the corresponding results \cite[Lemmas~2.7 and~2.9]{Gloria-Otto-09} in the discrete setting, the proofs of which are first presented in the continuum setting considered here (where algebraic decay can be replaced by the exponential decay
stated here). For \eqref{eq:int-grad2}, which we shall only use to prove the following lemma, the proof is similar and the Meyers' argument is used twice: once on each variable.
\begin{lemma}[optimal annealed pointwise estimates of gradients]\label{lem:annealed-estim}
Let $\expec{\cdot}$ be 
a stationary ensemble, and for all $A\in\Omega$ denote by $G_T$ the associated Green function, $T>0$. 
Then, there exists $\bar{p}>1$ depending only on $\lambda$, and $d$ such that for all $\bar{p}\geq  p\geq 1$ 
and all $|y|\gg 1$,
\begin{eqnarray}
\expec{|\nabla_yG_T(y,0)|^{2p}}^{\frac{1}{2p}}&\lesssim& |y|^{1-d}\exp(-c\frac{|y|}{\sqrt{T}}),\label{eq:annealed-estim} \\
\expec{|\nabla\nabla G_T(y,0)|}&\lesssim& |y|^{-d}\exp(-c\frac{|y|}{\sqrt{T}}).
\label{eq:annealed-estim2}
\end{eqnarray}
\qed
\end{lemma}
\noindent For $p=1$, \eqref{eq:annealed-estim} is a consequence of the annealed estimates \cite{Delmotte-Deuschel-05}
by Delmotte and Deuschel on the parabolic Green function for stationary ensembles.
We prove Lemma~\ref{lem:annealed-estim} by combining the Meyers' estimates of Lemma~\ref{lem:int-grad} with the elliptic approach of the Delmotte-Deuschel result developed by Marahrens and the second author
in  \cite{Marahrens-Otto-14}.
Although the estimate~\eqref{eq:annealed-estim2} on the mixed second derivative is not used in this article, it is stated here for future reference. 
\begin{lemma}[suboptimal quenched pointwise estimates of gradients]\label{lem:ptwise-estim}
Let $A\in\Omega$ and $G_T$ be its associated Green function, $T>0$. 
Then, there exists $\alpha>0$ depending only on $\lambda$ such that for all 
$R\sim 1$ and all $|z|>2R$, we have
\begin{equation}\label{eq:Green-ptwise-estim}
\Big(\int_{B_R(z)}|\nabla_{z'} G_T(z',0)|^2dz'\Big)^{\frac{1}{2}}
\,
\lesssim 
\,
\left\{
\begin{array}{lll}
|z|^{-\alpha} \exp(-c\frac{|z|}{\sqrt{T}}) &\text{ for }&d=2,
\\
|z|^{2-d} \exp(-c\frac{|z|}{\sqrt{T}}) &\text{ for }&d>2.
\end{array}
\right\}
\end{equation}
\qed
\end{lemma}
\noindent This lemma (which is suboptimal for $d>2$ but sufficient for our purpose) 
follows from Caccioppoli's inequality and the following precised energy estimate 
that we shall use in the proof of Lemma~\ref{lem:gradphi_estim}.
\begin{lemma}[Precised energy estimates]\label{lem:modif-ener}
There exists an exponent $\alpha(d,\lambda)>0$ such
that for all $A\in \Omega$,
\begin{itemize}
\item for all $R\ge 1$, $T>0$ and any function $v\in H^1(B_R)$ satisfying 
\begin{equation}\label{eq:modif-ener-3.1}
T^{-1}v-\nabla\cdot A\nabla v=0
\end{equation}
we have
\begin{equation}\label{eq:modif-ener-1.1}
\left(\int_{B_1}(T^{-1}v^2+|\nabla v|^2)dx\right)^\frac{1}{2}\lesssim 
R^{-\alpha}\left(\int_{B_R}(T^{-1}v^2+|\nabla v|^2)dx\right)^\frac{1}{2}.
\end{equation}
\item for all $T>0$ and functions $v\in H^1(\R^d)$ and vector fields $g\in L^2(\R^d,\R^d)$ related by
\begin{equation}\label{eq:modif-ener-3.2}
T^{-1}v-\nabla\cdot A\nabla v=\nabla\cdot g,
\end{equation}
and all radii $R$, we have
\begin{equation}\label{eq:modif-ener-1.2}
\left(\int_{B_R}(T^{-1}v^2+|\nabla v|^2)dx\right)^\frac{1}{2}\lesssim 
\left(\int_{\R^d}(\frac{|x|}{R}+1)^{-2\alpha}|g(x)|^2dx\right)^\frac{1}{2}.
\end{equation}
\item for all $R\ge 1$, $T>0$, the modified corrector $\ext\phi_T$ satisfies
\begin{multline}\label{eq:modif-ener-1.3}
\left(\int_{B_1}(T^{-1}(\xi\cdot x+\ext\phi_T(x))^2+|\xi+\nabla \ext\phi_T(x)|^2)dx\right)^\frac{1}{2}
\\ \lesssim \,
R^{-\alpha} \left(\int_{B_R}(T^{-1}(\xi\cdot x+\ext\phi_T(x))^2+|\xi+\nabla \ext\phi_T(x)|^2+T^{-1}R^2)dx\right)^\frac{1}{2}.
\end{multline}
\end{itemize}
\qed
\end{lemma}
\medskip
\noindent 
As in \cite{Gloria-Otto-09}, the proof of Theorem~\ref{th:main-1} relies on 
(SG)
and on Proposition~\ref{prop:main-1}.
As opposed to our proof in the discrete setting we shall replace the
use of convolution estimates of Green functions (cf. \cite[Lemma~2.10]{Gloria-Otto-09} and \cite[Estimate A.8]{Gloria-Mourrat-10}) by a suitable use of the pointwise estimates of Lemmas~\ref{lem:annealed-estim} and~\ref{lem:ptwise-estim}.

\medskip
\noindent  Our proof of Proposition~\ref{cor:main-2} is new, and significantly differs from the corresponding proofs for the discrete setting in \cite{Gloria-Otto-09b,Gloria-Mourrat-10}.
Since the function $(0,\infty)\to \calH^1, T\mapsto \phi_T$ is smooth, we may define $\psi_T:=T^2\frac{\partial\phi_{T}}{\partial T}\in \calH^1$.
As for the corresponding proof in the discrete setting,  
we have to estimate the quantity $\expec{\phi_T\psi_T}=\cov{\phi_T}{\psi_T}$. In \cite{Gloria-Otto-09b} we used the covariance estimate of \cite[Lemma~3]{Gloria-Otto-09b} as a starting point. 
In the case of the Poisson point process, a corresponding covariance estimate holds as well and is known
as the Harris--FKG inequality, see \cite{Wu-00,Last-Penrose-11}. 
It is however not clear whether (SG) implies a covariance inequality in general.
A first possibility to avoid the use of a covariance estimate is to appeal to spectral theory (in the case of symmetric coefficients) 
to bound this covariance using the variance of $\psi_T$, in the spirit of \cite{Mourrat-10,Gloria-Mourrat-10} in the discrete setting (see also Corollary~\ref{coro:spectr-bottom}).
In our proof of Proposition~\ref{cor:main-2} however, we use neither a covariance estimate nor spectral theory.
To apply Lemma~\ref{lem:var-estim} to $\psi_T$, one needs to control the susceptibility of $\psi_T$ (in the spirit of Lemma~\ref{lem:diff-phi}).
\begin{lemma}\label{lem:diff-psi}
Let $A\in \Omega$, and let $G_T$ and $\ext\phi_T$ be the associated Green function and modified corrector for $\xi\in\R^d$, $|\xi|=1$, and $T>0$.
We set
\begin{equation}\label{eq:diff-psi}
\ext\psi_T\,=\,T^2 \frac{\partial \ext \phi_T}{\partial T},
\end{equation}
and note that  $\ext\psi_T\in H^1_\loc(\R^d)$ is the unique distributional solution in the class of functions $\chi$ in $H^1_\loc(\R^d)$ such that 
$\limsup_{t\uparrow \infty}\fint_{B_t} (\chi^2+|\nabla \chi|^2)dx \,<\,\infty$ of
\begin{equation}\label{eq:equation-psiT}
T^{-1}\ext\psi_T-\nabla \cdot A\nabla \ext\psi_T\,=\,\ext\phi_T.
\end{equation}
For all $R\sim 1$, $T>0$, and for all $x,z\in \R^d$, we have
\begin{multline}\label{eq:diff-psi-1}
\dps\osc{A|_{B_R(z)}}{\ext\psi_T(x)}\,\lesssim \, {\ext h}_{T}(z,x) \left(\int_{B_{3R}(z)}|\nabla \ext\psi_T(z')|^2dz'+ \nu_d(T)\Big(\int_{B_{9R}(z)}|\nabla \ext\phi_T(z')|^2dz' +1 \Big)\right)^{\frac{1}{2}}
\\+\left(\int_{B_{3R}(z)}|\nabla \ext\phi_T(z')|^2dz'+1\right)^{\frac{1}{2}} \int_{\R^d}g_{T}(x-y){\ext h}_{T}(z,y)dy,
\end{multline}
where $\nu_d(T)$ is given by 
\begin{equation}\label{eq:def-nud}
\nu_d(T)\,=\,\left\{
\begin{array}{rcl}
d=2&:& T\ln T,\\
d=3&:& \sqrt{T},\\
d=4&:&\ln T,\\
d>4&:&1,
\end{array}
\right.
\end{equation}
${\ext h}_T$ is as in \eqref{eq:def-frak-h}, 
and $g_T$  as in \eqref{eq:ptwise-decay-estim}.
In addition, we also have
\begin{equation}\label{eq:sup-psi-1}
\sup_{\dps A|_{B_R(z)}}\int_{B_R(z)}|\nabla \ext\psi_T(y)|^2dy \,\lesssim \, \int_{B_R(z)}|\nabla \ext\psi_T(y)|^2dy+\nu_d(T)\Big(\int_{B_{3R}(z)}|\nabla \ext\phi_T(y)|^2dy+1\Big).
\end{equation}
\qed
\end{lemma}
\noindent The oscillation of $\psi_T$ involves integrals of products of the Green function and of its gradient, the expectation of which is controlled using the pointwise estimates \eqref{eq:ptwise-decay-estim} in Definition~\ref{def:Green} on the Green
function and the 
pointwise estimates of Lemmas~\ref{lem:annealed-estim} and~\ref{lem:ptwise-estim} on its gradient.


\section{Proofs of the main results}

\noindent 
In this section we prove Proposition~\ref{prop:main-1} in the form
of Lemmas~\ref{lem:gradphi_to_phi} and~\ref{lem:gradphi_estim}, Theorem~\ref{th:main-1}, Proposition~\ref{cor:main-2}, and Theorem~\ref{th:main-2}. Recall that we assume that (SG) and $q$-(SG) hold with $\tilde \ell=\ell=2$.

\subsection{Proof of Lemma~\ref{lem:gradphi_to_phi}}

In this proof the multiplicative constants in $\lesssim$ may depend on $q\ge 1$.
We split the proof into three steps.

\medskip
\step{1} Application of (SG): Proof of
\begin{equation}\label{eq:prop-01}
\expec{\phi_T^{2q}}^{\frac{1}{q}}\;\lesssim \; \expec{\left(\int_{\R^d} {\ext h}_T^2(z,0)\Big(\int_{B_{6}(z)}|\nabla \ext\phi_T(y)|^2dy+1 \Big) dz\right)^q}^{\frac{1}{q}}.
\end{equation}
By Corollary~\ref{coro:var-estim}, which
we apply to $X=\phi_T$, we have for all $q\ge 1$
\begin{equation*}
\expec{\phi_T^{2q}}^{\frac{1}{q}}\;\lesssim \; \expec{\left(\int_{\R^d} \Big( \osc{A|_{B_2(z)}}{\phi_T}\Big)^2dz\right)^q}^{\frac{1}{q}}.
\end{equation*}
From Lemma~\ref{lem:diff-phi} with $R=2$ we learn that 
\begin{equation*}
\left(\osc{A|_{B_2(z)}}{\ext\phi_T(0)} \right)^2\,\lesssim \, {\ext h}_T^2(z,0)\left(\int_{B_{6}(z)}|\nabla \ext\phi_T(y)|^2dy+1 \right).
\end{equation*}
This yields \eqref{eq:prop-01}.

\step{2} Dyadic decomposition of $\R^d$ and use of stationarity: Proof of
\begin{multline}
\expec{\phi_T^{2q}}^{\frac{1}{q}}\;\lesssim \;  \left( 1+\sum_{i\in \N} (2^i R)^{\frac{d}{q}} \Big( \sup_{A\in \Omega} \int_{2^iR <|z|\leq 2^{i+1}R } {\ext h}_T^{2\frac{q}{q-1}}(z,0)dz\Big)^{\frac{q-1}{q}}  \right) \\
\times \left( \expec{\Big(\int_{B_{6}} |\nabla\ext \phi_T(y)|^{2}dy\Big)^q}+1\right)^{\frac{1}{q}} \label{eq:pr-prop-bded-1}
\end{multline}
for $R\geq 2$ such that Lemma~\ref{lem:int-grad} holds.

\noindent Since we control $\nabla_y G_T(y,0)$ well when integrated over dyadic annuli, we decompose $\R^d$ into the ball $\{|z|\leq 2R\}$
and the annuli $\{2^i R <|z|\leq 2^{i+1}R\}$ for $i\in \N$.
The triangle inequality in $L^q(\Omega)$ on the RHS of \eqref{eq:prop-01} yields
\begin{multline*}
\expec{\phi_T^{2q}}^{\frac{1}{q}}\;\lesssim \; \expec{\left(\int_{|z|\leq 2R} {\ext h}_T^2(z,0) \Big(\int_{B_{6}(z)}|\nabla \ext \phi_T(y)|^2dy+1 \Big) dz\right)^q}^{\frac{1}{q}}
\\
+\sum_{i\in \N} \expec{\left(\int_{2^i R <|z|\leq 2^{i+1}R} {\ext h}_T^2(z,0) \Big(\int_{B_{6}(z)}|\nabla \ext \phi_T(y)|^2dy+1 \Big) dz\right)^q}^{\frac{1}{q}}.
\end{multline*}
Since ${\ext h}_T \,\lesssim\,1$ pointwise by definition, cf. \eqref{eq:def-frak-h}, 
the stationarity of $\nabla \ext\phi_T$ yields
\begin{multline*}
\expec{\left(\int_{|z|\leq 2R}{\ext h}_T^2(z,0) \Big(\int_{B_{6}(z)}|\nabla \ext\phi_T(y)|^2dy+1 \Big) dz\right)^q}^{\frac{1}{q}}
\\\,\lesssim \, \left( \expec{\Big(\int_{B_{6}} |\nabla \ext\phi_T(y)|^{2}dy\Big)^q}+1\right)^{\frac{1}{q}}.
\end{multline*}
For the other terms, we use H\"older's inequality in the $z$-integral with exponents $(\frac{q}{q-1},q)$, bound the integral involving ${\ext h}_T$ by
its supremum over $\Omega$, and then appeal again to the stationarity of $\nabla \ext\phi_T$:
\begin{eqnarray*}
\lefteqn{\expec{\left(\int_{2^i R <|z|\leq 2^{i+1}R} {\ext h}_T^2(z,0) \left(\int_{B_{6}(z)}|\nabla \ext\phi_T(y)|^2dy+1 \right) dz\right)^q}}  \\
&\leq & \left\langle\left(\int_{2^i R <|z|\leq 2^{i+1}R} {\ext h}_T^{2\frac{q}{q-1}}(z,0)dz\right)^{q-1}  \int_{2^i R <|z|\leq 2^{i+1}R} 
\left(\int_{B_{6}(z)}|\nabla \ext\phi_T(y)|^2dy+1 \right)^q dz\right\rangle \\
&\lesssim &\Big( \sup_{A\in \Omega} \int_{2^iR <|z|\leq 2^{i+1}R } {\ext h}_T^{2\frac{q}{q-1}}(z,0)dz\Big)^{q-1} (2^iR)^d \left( \expec{\Big(\int_{B_{6}} |\nabla \ext\phi_T(y)|^{2}dy\Big)^q}+1\right).
\end{eqnarray*}
Estimate~\eqref{eq:pr-prop-bded-1} then follows from summing over $i$.

\medskip
\step{3} Choice of $q$ and estimate of the Green function: Proof of 
\begin{equation}\label{eq:pr-prop-bded-2}
\expec{\phi_T^{2q}}^{\frac{1}{q}} \;\lesssim \;\left( \expec{\Big(\int_{B_{6}} |\nabla \ext\phi_T(y)|^{2}dy\Big)^q}^{\frac{1}{q}}+1\right) \left\{ 
\begin{array}{lll}
d=2&:&\ln T,\\
d>2&:&1,
\end{array}
\right\}
\end{equation}
for all $q$ large enough so that $\frac{q}{q-1}\leq \bar p$, where $\bar p$ is the Meyers exponent of Lemma~\ref{lem:int-grad}.

\noindent By definition \eqref{eq:def-frak-h} of ${\ext h}_T$ and H\"older's inequality, we have for all $i\in \N$
\begin{equation*}
\int_{2^iR <|z|\leq 2^{i+1}R } {\ext h}_T^{2\frac{q}{q-1}}(z,0)dz \,\lesssim \, \int_{(2^i-1)R <|z|\leq (2^{i+1}+1)R } |\nabla G_T(z,0)|^{2\frac{q}{q-1}}dz.
\end{equation*}
Estimate~\eqref{eq:int-grad} yields a bound for all $i\in \N$ which is uniform in $A\in \Omega$:
\begin{equation*}
\int_{2^iR <|z|\leq 2^{i+1}R } {\ext h}_T^{2\frac{q}{q-1}}(z,0)dz \,\lesssim \, (2^iR)^d (2^iR)^{(1-d) 2\frac{q}{q-1}}
\exp(-c\frac{2q}{q-1}\frac{(2^i-1)R}{\sqrt{T}}).
\end{equation*}
Combined with \eqref{eq:pr-prop-bded-1}, and the fact that $R$ is of order $1$, this yields
\begin{eqnarray*}
\expec{\phi_T^{2q}}^{\frac{1}{q}}&\lesssim &  \left( \expec{\Big(\int_{B_{6}} |\nabla \ext\phi_T(y)|^{2}dy\Big)^q}^{\frac{1}{q}}+1\right) \\
&& \qquad \times \sum_{i\in \N_0} \left((2^iR)^{d} \big((2^{i}R)^{d} (2^iR)^{(1-d) 2\frac{q}{q-1}} \big)^{q-1}\right)^{\frac{1}{q}}\exp(-c\frac{2^iR}{\sqrt{T}}). \\
&=& \left( \expec{\Big(\int_{B_{6}} |\nabla \ext\phi_T(y)|^{2}dy\Big)^q}^{\frac{1}{q}}+1\right) \sum_{i\in \N_0} (2^iR)^{2-d} \exp(-c\frac{2^iR}{\sqrt{T}})   .
\end{eqnarray*}
Since for $d>2$, the sum on the RHS is bounded independently of $T$, and for $d=2$ the sum is bounded by $\ln T$, \eqref{eq:pr-prop-bded-2} follows.

\subsection{Proof of Lemma~\ref{lem:gradphi_estim}}

We split the proof into four steps, and combine the approach without Green's functions we developed in \cite{Gloria-Otto-14}
with a compactness argument developed by Bella and the second author for systems \cite{Bella-Otto-14}, which we extend
from bounded domains with periodic boundary conditions to the whole space with the massive term.
In the first step we decompose $\nabla \ext\phi_T$ in Fourier modes, and show it is enough to consider a finite number of Fourier coefficients.
In the second step we estimate the oscillation of the Fourier coefficients, and apply $q$-(SG) and elliptic regularity in the third step to obtain a nonlinear estimate. We conclude in the fourth step.

\medskip

\step{1} Compactness argument: We argue that for any $\delta>0$ and any radius
\begin{equation}\label{eq:co-u.31}
R\le\sqrt{T},
\end{equation}
there exist $N(d,\delta)$ linear functionals $F_0,\cdots, F_{N-1}:H^1(B_{2R})\to \R$ bounded in the sense that
\begin{equation}\label{eq:co-u.17}
|F_n u|\le\left(\int_{B_{2R}}(T^{-1}u^2+|\nabla u|^2)dx\right)^\frac{1}{2},
\end{equation}
and which have the property that
for any pair of functions $u:H^1(B_{2R})$ and $f\in L^2(B_{2R})$ related by
\begin{equation}\label{eq:co-u.15}
T^{-1}u-\nabla\cdot A\nabla u=T^{-1}f
\end{equation}
we have
\begin{equation}\label{eq:co-u.24}
\int_{B_R}({T}^{-1}u^2+|\nabla u|^2)dx\,\lesssim\,
\sum_{n=0}^{N-1}|F_n u|^2+\delta\int_{B_{2R}}|\nabla u|^2dx+\int_{B_{2R}}{T}^{-1}f^2dx,
\end{equation}
where the multiplicative constant is independent of $\delta,T$ and $R$.
We split the estimate \eqref{eq:co-u.24} into an a priori estimate for \eqref{eq:co-u.15}, namely
\begin{equation}\label{eq:co-u.20}
\int_{B_R}({T}^{-1}u^2+|\nabla u|^2)dx\,\lesssim\,
\int_{B_{2R}}({T}^{-1}u^2+R^{-2}(u-\bar u)^2+{T}^{-1}f^2)dx,
\end{equation}
where $\bar u$ denotes the average of $u$ in $B_{2R}$,
and into the construction of the functionals $F_n$ such that for an {\it arbitrary} function $u\in H^1(B_{2R})$
\begin{equation}\label{eq:co-u.16}
\int_{B_{2R}}({T}^{-1}u^2+R^{-2}(u-\bar u)^2)dx\,\lesssim\,
\sum_{n=0}^{N-1}(F_n u)^2+\delta\int_{B_{2R}}|\nabla u|^2dx.
\end{equation}
We start with \eqref{eq:co-u.16}, which thanks to \eqref{eq:co-u.31} we may split into
\begin{equation}\label{eq:co-u.16bis}
R^{-2}\int_{B_{2R}}(u-\bar u)^2dx\,\lesssim\,
\sum_{n=1}^{N-1}(F_n u)^2+\delta\int_{B_{2R}}|\nabla u|^2dx
\end{equation}
and
\begin{equation}\nonumber
{T}^{-1}\int_{B_{2R}}\bar u^2dx \,\le\,(F_0 u)^2,
\end{equation}
where the last estimate is trivially satisfied (as an identity) by defining 
$F_0u=\sqrt{\frac{|B_{2R}|}{T}}\bar u=\frac{\int_{B_{2R}}udx}{\sqrt{T|B_{2R}|}}$,
which by Jensen's inequality satisfies the boundedness condition (\ref{eq:co-u.17}) in the simple
form of $(F_0u)^2\le{T}^{-1}\int_{B_{2R}}u^2dx$.
We thus turn to \eqref{eq:co-u.16bis}; by rescaling length according to $x=R\hat x$,
we may assume that $2R=1$. Let $\{(\lambda_n,u_n)\}_{n=0,1,\cdots}$ denote 
a complete set of increasing eigenvalues and $L^2$-orthonormal eigenfunctions of $-\triangle$ on $B_1$ endowed
with homogeneous Neumann boundary conditions, that is
\begin{equation}\label{eq:co-u.18}
\int_{B_1}\nabla v\cdot\nabla u_ndx=\lambda_n\int_{B_1}v u_ndx\quad
\mbox{for all functions }\;v\in H^1(B_1).
\end{equation}
In particular, we have $\int_{B_1}|\nabla u_n|^2dx=\lambda_n\int_{B_1}u_n^2dx=\lambda_n$.
We also note that $\lambda_1>0$. Hence for all $n\ge 1$ 
\begin{equation}\label{eq:co-u.19}
F_nu=\int_{B_1}\nabla u\cdot\frac{\nabla u_n}{\sqrt{\lambda_n}}dx\quad\mbox{for all functions}\;u\in H^1(B_1)
\end{equation}
defines a linear functional $F_n$ on vector fields that has the boundedness property \eqref{eq:co-u.17}
in form of $(F_nu)^2\le\int_{B_1}|\nabla u|^2dx$.
By completeness of the orthonormal system $\{u_n\}_{n=0,1,\cdots}$, Plancherel
and $u_0=const$, we have
\begin{eqnarray*}
\lefteqn{\int_{B_1}(u-\bar u)^2dx=\sum_{n=1}^\infty\left(\int_{B_1}u u_ndx\right)^2}\\
&\stackrel{(\ref{eq:co-u.18})}{=}&
\sum_{n=1}^\infty\frac{1}{\lambda_n}
\left(\int_{B_1}\nabla u\cdot \frac{\nabla u_n}{\sqrt{\lambda_n}}dx\right)^2\\
&\le&\frac{1}{\lambda_1}\sum_{n=1}^{N-1}\left(\int_{B_1}\nabla u\cdot \frac{\nabla u_n}{\sqrt{\lambda_n}}dx\right)^2
+\frac{1}{\lambda_{N}}\sum_{n=N}^\infty\left(\int_{B_1}\nabla u\cdot \frac{\nabla u_n}{\sqrt{\lambda_n}}dx\right)^2.
\end{eqnarray*}
We note that $\eqref{eq:co-u.18}$ yields that also $\{\frac{\nabla u_n}{\sqrt{\lambda_n}}\}_{n=1,\cdots}$
is orthonormal, so that the above together with definition \eqref{eq:co-u.19} yields
\begin{eqnarray*}
\int_{B_1}(u-\bar u)^2dx
&\le&\frac{1}{\lambda_1}\sum_{n=1}^{N-1}(F_n\nabla u)^2
+\frac{1}{\lambda_{N}}\int_{B_1}|\nabla u|^2dx.
\end{eqnarray*}
Because of $\lim_{N\uparrow\infty}\lambda_{N}=\infty$, this implies \eqref{eq:co-u.16bis} in its
$(R=2)$-version.

\medskip

\noindent We now turn to \eqref{eq:co-u.20}; it is obviously enough to show
\begin{equation}\nonumber
\int_{B_R}|\nabla u|^2dx\lesssim
\int_{B_{2R}}({T}^{-1}(u-f)^2+R^{-2}(u-\bar u)^2)dx.
\end{equation}
By rescaling length according to $x=\sqrt{T}\hat x$, it is enough
to establish the case of $T=1$, that is,
\begin{equation}\label{eq:co-u.20bis}
\int_{B_R}|\nabla u|^2dx\,\lesssim\,
\int_{B_{2R}}((f-u)^2+R^{-2}(u-\bar u)^2)dx.
\end{equation}
We test \eqref{eq:co-u.15} for $T=1$, that is,
\begin{equation}\label{eq:co-u.21}
-\nabla\cdot A\nabla u=f-u
\end{equation}
with $\eta^2(u-\bar u)$, where $\eta$ is a cut-off function for $B_R$ in $B_{2R}$:
\begin{equation}\nonumber
\int_{B_{2R}}\eta^2\nabla u\cdot A\nabla udx=\int_{B_{2R}}\eta(u-\bar u)(-2\nabla\eta\cdot A\nabla u+\eta (f-u))dx,
\end{equation}
which by the properties of $A$ turns into
\begin{equation}\nonumber
\lambda\int_{B_{2R}}\eta^2|\nabla u|^2dx\,\le\,\int_{B_{2R}}\eta|u-\bar u|(2|\nabla\eta||\nabla u|+\eta|f-u|)dx.
\end{equation}
Using Young's inequality this gives
\begin{equation}\nonumber
\int_{B_{2R}}\eta^2|\nabla u|^2dx\,\lesssim\,\int_{B_{2R}}(|\nabla\eta|^2(u-\bar u)^2+\eta^2 (f-u)^2)dx,
\end{equation}
which by choice of $\eta$ turns into the desired
\begin{equation}\nonumber
\int_{B_R}|\nabla u|^2dx\lesssim\int_{B_{2R}}(R^{-2}(u-\bar u)^2+(f-u)^2)dx.
\end{equation}

\medskip
\step{2} Oscillation estimate of the Fourier coefficients.

\noindent Let $\alpha>0$ be the exponent of Lemma~\ref{lem:modif-ener} and let $F_n$
denote the functionals of Step~1 on $H^1(B_{2R})$. In this step we argue that  for all $n\in \N_0$,
\begin{equation}\label{eq:osc-estim-Fourier}
\int_{\R^d} \Big(\osc{A|_{B_2(z)}}{F_{n}(\xi+\nabla \ext\phi_T)}\Big)^2dz\,\lesssim \,  \sup_{z\in \R^d} \Big\{\big(|\frac{z}{R}|+1)^{-2\alpha}\Big(\int_{B_2(z)}|\xi+\nabla \ext\phi_T|^2dx\Big)\Big\}
\end{equation}
Let $F$ and $u$ denote any of the $F_{n}$ and $\frac{u_{n}}{\sqrt{\lambda_n}}$.
Assume first that $A$ is a smooth coefficient field.
For all $z\in \R^d$, let $A_z$ be a smooth coefficient field that coincides with $A$ on $\R^d\setminus B_2(z)$, and denote
by $\ext \phi_{T,z}$ the modified corrector associated with $A_z$.
We first claim that it is enough to prove that for all $\chi\in L^2(\R^d)$,
\begin{multline}\label{eq:osc-estim-Fourier-1}
\Big(\int_{\R^d} \chi(z) \int_{B_{2R}}(\nabla \ext\phi_T(x)-\nabla \ext\phi_{T,z}(x))\cdot \nabla u(x)dxdz\Big)^2
\\
 \lesssim \, \Big(\int_{\R^d}\chi^2dz\Big) \sup_{z\in \R^d} \Big\{\big(|\frac{z}{R}|+1)^{-2\alpha}\int_{B_2(z)}|\xi+\nabla \ext\phi_T|^2dx\Big\}.
\end{multline}
Indeed, since $\chi$ is arbitrary and the RHS does not depend on $\{A_z, z\in \R^d\}$, this implies that 
\begin{multline*}
\int_{\R^d} \sup_{A_z} \Big|\int_{B_{2R}}(\nabla \ext\phi_T(x)-\nabla \ext\phi_{T,z}(x))\cdot \nabla u(x)dx\Big|^2dz
\\
\lesssim \,  \sup_{z\in \R^d} \Big\{\big(|\frac{z}{R}|+1)^{-2\alpha}\int_{B_2(z)}|\xi+\nabla \ext\phi_T|^2dx\Big\},
\end{multline*}
from which \eqref{eq:osc-estim-Fourier} follows by density (to relax the assumption that $A$ be smooth) and the elementary estimate 
$$
\osc{A|_{B_2(z)}}{F(\xi+\nabla \ext\phi_T)}\,\leq\,2\sup_{A_z} |F(\xi+\nabla\ext \phi_T)-F(\xi+\nabla \ext\phi_{T,z})|.
$$
We now prove \eqref{eq:osc-estim-Fourier-1}.
Set $v(x):=\int_{\R^d} \chi(z) (\ext\phi_T(x)-\ext\phi_{T,z}(x))dz$. By Fubini's theorem,
and since $\nabla u$ has $L^2(B_{2R})$-norm unity,
\begin{multline*}
\Big(\int_{{\R^d}} \chi(z) \int_{B_{2R}}(\nabla \ext\phi_T(x)-\nabla \ext\phi_{T,z}(x))\cdot \nabla u(x)dxdz\Big)^2
\,=\,\Big(\int_{B_{2R}}\nabla v(x) \cdot \nabla u(x)dx\Big)^2 
\\
\lesssim \, \int_{B_{2R}}|\nabla v|^2dx \int_{B_{2R}} |\nabla u|^2dx\,=\,\int_{B_{2R}} |\nabla v|^2dx.
\end{multline*}
Since $v$ satisfies
$$
T^{-1}v-\nabla \cdot A\nabla v\,=\,\nabla\cdot \Big(\int_{\R^d}\chi(z)(A-A_z)(\xi+\nabla \ext\phi_{T,z}) dz\Big)
$$
on $\R^d$, we deduce by \eqref{eq:modif-ener-1.2} in Lemma~\ref{lem:modif-ener} that
\begin{eqnarray*}
\int_{B_{2R}} |\nabla v|^2dx&\lesssim& \int_{\R^d} \big(\frac{|x|}{R}+1\big)^{-2\alpha} 
\Big(\int_{\R^d}|\chi(z)||A(x)-A_z(x)||\xi+\nabla \ext\phi_{T,z}(x)| dz\Big)^2dx \\
&=&\int_{\R^d}\int_{\R^d}\int_{\R^d}\big(\frac{|x|}{R}+1\big)^{-2\alpha} |\chi(z)||A(x)-A_z(x)|
\\
&&\quad\times |\xi+\nabla \ext\phi_{T,z}(x)|
|\chi(z')||A(x)-A_{z'}(x)||\xi+\nabla \ext\phi_{T,z'}(x)|dxdzdz'\\
&\leq & \int_{\R^d}\int_{\R^d}\int_{\R^d}\big(\frac{|x|}{R}+1\big)^{-2\alpha} \chi^2(z)|\xi+\nabla \ext\phi_{T,z}(x)|^2
\\
&&\qquad \qquad \times|A(x)-A_z(x)|
|A(x)-A_{z'}(x)|dz'dxdz
\\
&\lesssim & \int_{\R^d}\chi^2(z)\int_{B_2(z)}\big(\frac{|x|}{R}+1\big)^{-2\alpha} |\xi+\nabla \ext\phi_{T,z}(x)|^2dx dz
\\
&\lesssim &\Big(\int_{\R^d}\chi^2dz\Big) \sup_{z\in \R^d} \Big\{ \big(\frac{|z|}{R}+1\big)^{-2\alpha} 
\int_{B_2(z)} |\xi+\nabla \ext\phi_{T,z}|^2dx\Big\}.
\end{eqnarray*}
It remains to show that 
\begin{equation}\label{eq:osc-estim-Fourier-1.1}
\int_{B_2(z)} |\xi+\nabla \ext\phi_{T,z}|^2dx\,\lesssim \, \int_{B_2(z)}|\xi+\nabla \ext\phi_{T}|^2dx.
\end{equation}
Indeed, since $\delta \phi:=\ext\phi_{T,z}-\ext\phi_{T}$ satisfies
$$
T^{-1}\delta \phi-\nabla \cdot A_z\nabla \delta \phi\,=\,\nabla\cdot (A-A_z)(\xi+\nabla \ext\phi_{T}),
$$
an energy estimate yields
$$
\int_{\R^d}|\nabla \delta \phi|^2dx \, \lesssim \, \int_{B_2(z)}|\xi+\nabla \ext\phi_{T}|^2dx,
$$
from which \eqref{eq:osc-estim-Fourier-1.1} follows by the triangle inequality.

\medskip
\step{3} Application of (SG) and elliptic regularity: Proof of 
\begin{multline}\label{eq:qSG-Fourier-coeff}
\expec{(F(\xi+\nabla \ext\phi_T)-\expec{F(\xi+\nabla \ext\phi_T)})^{2q}}^{\frac{1}{q}} 
\\
\lesssim\ R^{\frac{d}{q}-2\alpha} \Big\langle\Big(\int_{B_{R}}(T^{-1}(\xi\cdot x+\ext\phi_T(x))^2+|\xi+\nabla \ext\phi_T(x)|^2+T^{-1}R^2)dx\Big)^q\Big\rangle^{\frac{1}{q}}
\end{multline}
for all $q\ge \frac{d+1}{2\alpha}$, where the multiplicative constant is independent of $T$ and $R$.

\noindent We first apply $q$-(SG) to $F(\xi+\nabla \ext\phi_T)$ and appeal to \eqref{eq:osc-estim-Fourier}
to get
\begin{multline}\label{eq:passage-2-2R-gamma1}
\expec{(F(\xi+\nabla \ext\phi_T)-\expec{F(\xi+\nabla \ext\phi_T)})^{2q}}^{\frac{1}{q}} \\
\lesssim\,  \expec{\Big( \sup_{z\in \R^d} \Big\{\big(|\frac{z}{R}|+1)^{-2\alpha}\int_{B_2(z)}|\xi+\nabla \ext\phi_T|^2dx\Big\}\Big)^q}^{\frac{1}{q}},
\end{multline}
where the multiplicative constant is independent of  $T$ and $R$.
Note that 
\begin{multline*}
\sup_{z\in \R^d} \Big\{\big(|\frac{z}{R}|+1)^{-2\alpha}\int_{B_2(z)}|\xi+\nabla \ext\phi_T|^2dx\Big\}^q
\\
\lesssim \, \int_{\R^d} \big(|\frac{z}{R}|+1\big)^{-2q\alpha}\Big(\int_{B_3(z)}|\xi+\nabla \ext\phi_T|^2dx\Big)^q,
\end{multline*}
where the multiplicative constant only depends on $d$.
Hence, by stationarity of the modified corrector and the estimate 
$$
 \int_{\R^d}\big(|\frac{x}{R}|+1\big)^{-2q\alpha}dz \, \lesssim \,R^d,
$$
which holds for all $q\ge \frac{d+1}{2\alpha}$, \eqref{eq:passage-2-2R-gamma1} turns into
\begin{equation*}
\expec{(F(\xi+\nabla \ext\phi_T)-\expec{F(\xi+\nabla \ext\phi_T)})^{2q}}^{\frac{1}{q}} \,
\lesssim\,R^{\frac{d}{q}} \Big\langle\Big(\int_{B_3}|\xi+\nabla \ext\phi_T|^2dx\Big)^q\Big\rangle^{\frac{1}{q}}.
\end{equation*}
We then appeal to \eqref{eq:modif-ener-1.3} in Lemma~\ref{lem:modif-ener}, that shows that for all $R\ge 6$
\begin{equation*}\label{eq:passage-2-2R-gamma}
\int_{B_3}|\xi+\nabla \ext\phi_T|^2dx \,\lesssim\, R^{-2\alpha}\int_{B_{R}}(T^{-1}(\xi\cdot x+\ext\phi_T(x))^2+|\xi+\nabla \ext\phi_T(x)|^2+T^{-1}R^2)dx.
\end{equation*}

\medskip
\step{4} Buckling and proof of \eqref{eq:estim-lem-gradphi_estim}.

\noindent 
By stationarity, there exists $C$ depending only on $d$ such that for all $q\ge 1$,
$$
\expec{\Big(\int_{B_{2R}} |\xi+\nabla \ext\phi_T|^2dx\Big)^q}^{\frac{1}{q}}\,\leq \, C\expec{\Big(\int_{B_{R}} |\xi+\nabla \ext\phi_T|^2dx\Big)^q}^{\frac{1}{q}}.
$$
Hence, from the first step for $u(x)=\xi\cdot x+\ext\phi_T(x)$ and $f(x)=-\xi\cdot x$, we learn by \eqref{eq:co-u.24} and the triangle inequality that 
for some $\delta>0$ small enough there exist some constant $C<\infty$ and $N\in \N$
such that for all $R>0$, all $T>0$ and all $q\ge 1$,
\begin{multline}\label{eq:buckling-nabla-phiT-1}
\expec{\Big(\int_{B_{R}} (T^{-1}(\xi\cdot x+ \ext\phi_T(x))^2+|\xi+\nabla \ext\phi_T(x)|^2)dx\Big)^q}^{\frac{1}{q}}
\\\leq \, C \max_{n\in \{0,\dots,N-1\}} \expec{(F_{n}(\xi+\nabla \ext\phi_T))^{2q}}^{\frac{1}{q}}+CT^{-1}R^{d+2}.
\end{multline}
In the rest of this step, $C$ may change from line to line but remains independent of $R$ and $T$.
Let $F$ denote any of the $F_n$. By the triangle inequality followed by Jensen's inequality, 
$$
\expec{(F(\xi+\nabla \ext\phi_T))^{2q}}^{\frac{1}{q}}\,\leq\, \expec{(F(\xi+\nabla \ext\phi_T)-\expec{F(\xi+\nabla \ext\phi_T)})^{2q}}^{\frac{1}{q}}+\expec{(F(\xi+\nabla \ext\phi_T))^{2}},
$$
so that the combination of \eqref{eq:buckling-nabla-phiT-1}, \eqref{eq:co-u.17} in Step~1, and of \eqref{eq:qSG-Fourier-coeff} in Step~3, yields for $q\ge \frac{d+1}{2\alpha}$ and $R\ge 6$,
\begin{multline}\label{eq:buckling-nabla-phiT-2}
\expec{\Big(\int_{B_{R}} (T^{-1}(\xi\cdot x+ \ext\phi_T(x))^2+|\xi+\nabla \ext\phi_T(x)|^2)dx\Big)^q}^{\frac{1}{q}}
\\ \leq\, CR^{\frac{d}{q}-2\alpha} \expec{\Big(\int_{B_{2R}} (T^{-1}(\xi\cdot x+ \ext\phi_T(x))^2+|\xi+\nabla \ext\phi_T(x)|^2)dx\Big)^q}^{\frac{1}{q}}\\
+\expec{\int_{B_{2R}}|\xi+\nabla \ext\phi_T|^2dx} + CT^{-1}R^{d+2}.
\end{multline}
By stationarity, there exists $C<\infty$ depending only on $d$ such that for all $q\ge 1$,
\begin{multline*}
\expec{\Big(\int_{B_{2R}} (T^{-1}(\xi\cdot x+ \ext\phi_T(x))^2+|\xi+\nabla \ext\phi_T(x)|^2)dx\Big)^q}^{\frac{1}{q}}
\\
\leq C \expec{\Big(\int_{B_{R}} (T^{-1}(\xi\cdot x+ \ext\phi_T(x))^2+|\xi+\nabla \ext\phi_T(x)|^2)dx\Big)^q}^{\frac{1}{q}},
\end{multline*}
to the effect that for all $q\ge \frac{d+1}{2\alpha}$ we can absorb the first RHS term of \eqref{eq:buckling-nabla-phiT-2}
into the LHS for $R$ large enough.
This yields by the energy estimate of Lemma~\ref{lem:app-corr}, the triangle inequality, and since $\sqrt{T}\ge R$ (as required in Step~1),
\begin{equation*}
\expec{\Big(\int_{B_{R}} |\nabla \ext\phi_T|^2dx\Big)^q}^{\frac{1}{q}}
\, \lesssim \, \expec{\int_{B_{2R}}|\xi+\nabla \ext\phi_T|^2dx} +R^d+T^{-1}R^{d+2} \,\lesssim \, R^d
\end{equation*}
for all $q\ge \frac{d+1}{2\alpha}$ (and therefore all $q\ge 1$ by Jensen's inequality) and $T\gg 1$ large enough.


\subsection{Proof of Theorem~\ref{th:main-1}}

\noindent Let us denote the spatial average of a function $h:\R^d\to \R$ with the averaging function $\eta_L$ by
\begin{equation*}
\moy{h}:=\int_{\R^d}h(x)\eta_L(x)dx,
\end{equation*}
where we recall that $\eta_L$ satisfies
\begin{equation}\label{eq:eta}
\eta_L:\R^d \to \R_+, \quad \supp{\eta_L} \subset B_L, \quad \int_{\R^d} \eta_L(x)\,dx=1, \quad |\nabla \eta_L|\lesssim L^{-d-1}. 
\end{equation}
The claim of the theorem is
\begin{equation}\label{eq:pro-th1-restatement}
\var{\moy{T^{-1}\ext\phi_T'\ext\phi_T+(\xi'+\nabla \ext \phi_T')\cdot A(\xi+\nabla\ext \phi_T)}}\,\lesssim \,
\left\{
\begin{array}{ll}
d=2: & L^{-2} \ln (2+\frac{\sqrt{T}}{L}),\\
d>2: & L^{-d},
\end{array}
\right.
\end{equation}
where $\ext \phi_T,\ext\phi_T'$ are the modified corrector and modified adjoint corrector associated with $A$ through Lemma~\ref{lem:def-phiT-conv} (with $A^*$ in place of $A$ for the adjoint corrector).

\medskip

\noindent This proof is an adaptation and simplification of the corresponding proof in the discrete setting, where we replace
convolution estimates by the triangle inequality combined with the pointwise annealed estimates of Lemma~\ref{lem:annealed-estim}.
The starting point is (SG) applied to  $$\calE=\moy{T^{-1}\ext\phi_T'\ext\phi_T+(\xi'+\nabla \ext \phi_T')\cdot A(\xi+\nabla\ext \phi_T)},$$
which yields
\begin{eqnarray}
\lefteqn{\var{\moy{T^{-1}\ext\phi_T'\ext\phi_T+(\xi'+\nabla \ext \phi_T')\cdot A(\xi+\nabla\ext \phi_T)}}}\nonumber \\
&\lesssim &\expec{\int_{\R^d} \left(\osc{A|_{B_2(z)}}{\moy{T^{-1}\ext\phi_T'\ext\phi_T+(\xi'+\nabla \ext \phi_T')\cdot A(\xi+\nabla\ext \phi_T)}}\right)^2dz}.\label{eq:main-var-start}
\end{eqnarray}

\medskip

\step{1} Sensitivity estimate for the averaged energy density:
\begin{equation}
{\osc{A|_{B_2(z)}}{\calE(A)} } \,\lesssim\,
L^{-(d+1)}\int_{B_L} Y_1(z,x) (Y_2(z)+Y_2(x))dx + (\sup_{B_2(z)}\eta_L)Y_2(z),
\label{eq:theo-step2-0sens}
\end{equation}
where for all $A\in \Omega$, $\calE(A)$ denotes the averaged energy
\begin{eqnarray*}
\calE(A)&:=&\moy{T^{-1}\ext\phi_T'\ext\phi_T+(\xi'+\nabla \ext \phi_T')\cdot A(\xi+\nabla\ext \phi_T)} ,
\end{eqnarray*}
and $Y_1$ and $Y_2$ are stationary random fields given by
\begin{multline*}
Y_1(x,z)\,:=\,\min \left\{
\left(\int_{B_1(x)}\int_{B_2(z)}|\nabla_yG_T(y,x')|^2dydx'\right)^{\frac{1}{2}}, 1 \right\}
\\+
\min \left\{
\left(\int_{B_1(x)}\int_{B_2(z)}|\nabla_yG_T'(y,x')|^2dydx'\right)^{\frac{1}{2}}, 1 \right\}
\end{multline*}
and
$$
Y_2(x)\,:=\,\int_{B_6(x)}|\nabla \ext\phi_T|^2dx'+\int_{B_6(x)}|\nabla \ext\phi_T'|^2dx'+1.
$$

\noindent Let $\tilde A$ coincide with $A$ outside $B_2(z)$, $z\in \R^d$.
We denote by  $\ext{\tilde \phi}_T$ and $\ext{\tilde \phi}_T'$ the modified corrector and adjoint corrector associated with $\tilde A$ so that $\calE(\tilde A)$ is given by 
\begin{eqnarray*}
\calE(\tilde A)&:=&\moy{T^{-1}\ext{\tilde \phi}_T'\ext{\tilde \phi}_T+(\xi'+\nabla \ext{\tilde \phi}_T')\cdot \tilde A(\xi+\nabla \ext{\tilde \phi}_T)} .
\end{eqnarray*}
We first derive a representation formula for the difference $\calE(A)-\calE(\tilde A)$: 
\begin{multline}\label{eq:theo-step1-1}
\calE(\tilde A)-\calE(A)\,=\,-\int_{\R^d}(\ext{\tilde \phi}_T'- \ext\phi_T')\nabla \eta_L\cdot \tilde A(\xi+\nabla \ext{\tilde \phi}_T)dx \\
+\int_{\R^d}(\ext \phi_T-\ext{\tilde \phi}_T)\nabla \eta_L\cdot A^*(\xi'+\nabla \ext \phi_T')dx\\+\int_{\R^d}(\xi'+\nabla \ext{\phi}_T')\cdot (\tilde A-A)(\xi+\nabla \ext{\tilde \phi}_T)\eta_Ldx.
\end{multline}
An elementary calculation yields
\begin{multline*}\label{eq:theo-step1-2}
{\calE(\tilde A)-\calE(A)}\,=\,T^{-1}\int_{\R^d}(\ext{\tilde \phi}_T'-\ext \phi_T')\ext{\tilde \phi}_T\eta_Ldx +\int_{\R^d}\nabla (\ext{\tilde \phi}_T'-\ext \phi_T')\cdot \tilde A(\xi+\nabla \ext{\tilde \phi}_T)\eta_Ldx \\
-T^{-1}\int_{\R^d}(\ext \phi_T-\ext{\tilde \phi}_T)\ext \phi_T'\eta_Ldx -\int_{\R^d}\nabla(\ext  \phi_T-\ext{\tilde  \phi}_T)\cdot A^*(\xi'+\nabla\ext \phi_T')\eta_Ldx \\
+\int_{\R^d}(\xi'+\nabla \ext{\phi}_T')\cdot(\tilde A-A)(\xi+\nabla\ext{\tilde  \phi}_T)\eta_Ldx.
\end{multline*}
This identity, combined with the weak form of the modified corrector and adjoint corrector equations~\eqref{eq:app-corr-distribution} for $\ext{\tilde \phi}_T$ and $\ext\phi_T'$ and test-functions $\eta_L(\ext{\tilde \phi}_T'-\ext \phi_T')$ and $\eta_L(\ext{\tilde \phi}_T-\ext \phi_T)$, turns into \eqref{eq:theo-step1-1}.

\medskip

\noindent
Since $A$ and $\tilde A$ coincide outside $B_2(z)$, one may bound $|\ext{\tilde \phi}_T(x)- \ext\phi_T(x)|$ and $|\ext{\tilde \phi}_T'(x)- \ext\phi_T'(x)|$ by the oscillations over $A|_{B_2(z)}$ of $\ext\phi_T(x)$ and $\ext \phi_T'(x)$, respectively, so that \eqref{eq:theo-step1-1} yields
\begin{multline*}
|\calE(\tilde A)-\calE(A)|\\\,\lesssim \,\int_{\R^d} \Big(\osc{A|_{B_2(z)}}{\ext \phi_T'}\Big)|\nabla \eta_L|(|\nabla \ext{\tilde \phi}_T|+1) dx
+\int_{\R^d} \Big(\osc{A|_{B_2(z)}}{\ext \phi_T}\Big)|\nabla \eta_L|(|\nabla \ext{\phi}_T'|+1) dx
\\+\int_{B_2(z)}(|\nabla \ext{ \phi}_T'|+1)(|\nabla \ext{\tilde\phi}_T|+1)\eta_Ldx.
\end{multline*}
Before we can take the supremum over $A$ and $\tilde A$ and use estimates~\eqref{eq:diff-phi-1} and~\eqref{eq:sup-phi-1} in Lemma~\ref{lem:diff-phi} (and the corresponding estimates for the adjoint correctors), we have to rewrite the RHS in terms of local square averages
of $\nabla \ext \phi_T'$ and $\nabla \ext{\tilde \phi}_T$.
To this purpose we introduce a new variable $y$ in the first RHS term via $\int_{\R^d}dx \lesssim \int_{\R^d} dx \int_{B_1(x)} dy$.
We then use Cauchy-Schwarz' inequality and take the supremum over $A|_{B_2(z)}$ and $\tilde A|_{B_2(z)}$. Since the RHS does not depend on $A|_{B_2(z)}$ and
$\tilde A|_{B_2(z)}$, it controls the oscillation of $\calE( A)$ with respect to $A|_{B_2(z)}$, and we have
\begin{multline*}
\osc{A|_{B_2(z)}}{\calE(A)} \\\,\lesssim \,\int_{\R^d} \left(\int_{B_1(x)}\Big(\osc{A|_{B_2(z)}}{\ext \phi_T'}\Big)^2 dy\right)^{\frac{1}{2}}(\sup_{B_1(x)}|\nabla \eta_L|)  \left(\sup_{\dps A|_{B_2(z)}}\int_{B_1(x)}|\nabla  \ext\phi_T|^2dy +1\right)^{\frac{1}{2}}dx 
 \\
+\int_{\R^d} \left(\int_{B_1(x)}\Big(\osc{A|_{B_2(z)}}{\ext \phi_T}\Big)^2 dy\right)^{\frac{1}{2}}(\sup_{B_1(x)}|\nabla \eta_L|)  \left(\sup_{\dps A|_{B_2(z)}}\int_{B_1(x)}|\nabla  \ext\phi_T'|^2dy +1\right)^{\frac{1}{2}}dx
\\ 
+(\sup_{B_2(z)}\eta_L) \, \left(\sup_{\dps A|_{B_2(z)}}\int_{B_2(z)}|\nabla  \ext\phi_T|^2dy+\sup_{\dps A|_{B_2(z)}}\int_{B_2(z)}|\nabla  \ext\phi_T'|^2dy+1\right).
\end{multline*}
An application of estimates~\eqref{eq:diff-phi-1} and~\eqref{eq:sup-phi-1} in Lemma~\ref{lem:diff-phi} (and the corresponding estimates for $\phi_T'$) with $R=2$  yields \eqref{eq:theo-step2-0sens} by Young's inequality and the properties of~$\eta$.

\medskip

\step{2} Proof of \eqref{eq:estim-var-hom}.

\noindent We apply the spectral gap estimate to $\calE$, use the oscillation estimate \eqref{eq:theo-step2-0},
and expand the square:
\begin{eqnarray}
\lefteqn{\var{\calE}}
\nonumber \\
 &\lesssim & \int_{\R^d} \expec{\big(\osc{A|_{B_2(z)}}{\calE(A)} \big)^2}dz
\nonumber\\
&\lesssim &  \int_{\R^d}\expec{\Big( L^{-(d+1)}\int_{B_L} Y_1(z,x) (Y_2(z)+Y_2(x))dx + (\sup_{B_2(z)}\eta_L) Y_2(z)\Big)^2}dx
\nonumber\\
&\lesssim & L^{-2(d+1)}\int_{B_L}\int_{B_L}\int_{\R^d} \expec{Y_1(z,x) (Y_2(z)+Y_2(x))Y_1(z,x') (Y_2(z)+Y_2(x'))}dzdxdx'
\nonumber\\
&&+\int_{\R^d}(\sup_{B_2(z)}\eta_L)^2\expec{Y_2^2(z)}dz.
\label{eq:theo-step2-0.1}
\end{eqnarray}
To estimate the RHS of \eqref{eq:theo-step2-0.1} we appeal to \eqref{eq:estim-prop2} in Proposition~\ref{prop:main-1} and Lemma~\ref{lem:annealed-estim}, which 
imply that for $\bar p$ as in Lemma~\ref{lem:annealed-estim} and all $q\ge 1$,
\begin{eqnarray}
\expec{|Y_1(z,x)|^{2\bar{p}}}^{\frac{1}{2\bar{p}}}&\lesssim & \frac{1}{1+|x-z|^{d-1}}\exp(-c\frac{|x-z|}{\sqrt{T}}),\label{eq:theo-step2-1}\\
\expec{|Y_2|^{q}}^{\frac{1}{q}}&\lesssim & 1.\label{eq:theo-step2-2}
\end{eqnarray}
Using \eqref{eq:theo-step2-2} for $q=4$ on the second RHS term and H\"older's estimate in probability with exponents
$(2\bar p, \frac{2\bar p}{\bar p-1},2\bar p, \frac{2\bar p}{\bar p-1})$ on the first RHS term followed by \eqref{eq:theo-step2-1} and \eqref{eq:theo-step2-2} for $q=4\frac{\bar p}{\bar p-1}$ then yields
\begin{eqnarray}
{\var{\calE}}
&\lesssim & L^{-2(d+1)} \int_{B_L} \int_{B_L}\int_{\R^d}\frac{1}{1+|x-z|^{d-1}} \exp(-c\frac{|x-z|}{\sqrt{T}})
\nonumber \\
&& \qquad\qquad \qquad\qquad \qquad \times\frac{1}{1+|x'-z|^{d-1}} \exp(-c\frac{|x'-z|}{\sqrt{T}})dzdxdx'
\nonumber\\
&&+\int_{\R^d}(\sup_{B_2(z)}\eta_L)^2dz.
\label{eq:theo-step2-0}
\end{eqnarray}
By definition of $\eta_L$, the second RHS term scales as $L^{-d}$.
For the first RHS term, we treat the cases $d=2$ and $d>2$ differently, and start with $d>2$.
In this case, we may discard the exponential cut-off and a direct calculation yields
\begin{equation*}
\int_{\R^d}\frac{1}{1+|x-z|^{d-1}} \frac{1}{1+|x'-z|^{d-1}} dz
\\\lesssim \, \frac{1}{1+|x-x'|^{d-2}},
\end{equation*}
whereas 
\begin{equation*}
\int_{B_L} \int_{B_L}\frac{1}{1+|x-x'|^{d-2}}dxdx' \,\lesssim\, L^{d+2},
\end{equation*}
so that the claim \eqref{eq:estim-var-hom} follows for $d>2$.

\noindent For $d=2$, we split the integral over $z$ into two parts: the integral over $B_{2L}$ and the integral over $\R^d\setminus B_{2L}$. On $B_{2L}$ we discard the exponential cut-off: 
\begin{equation*}
\int_{B_{2L}} \int_{B_L} \int_{B_L} \frac{1}{1+|x-z|} 
\frac{1}{1+|x'-z|} dxdx'dz \,\lesssim \, L^2  \int_{B_{2L}} \int_{B_{3L}} \frac{1}{1+|x|} 
\frac{1}{1+|x'|} dxdx'\,\lesssim \, L^{4},
\end{equation*}
whereas on $\R^d\setminus B_{2L}$ we take advantage of the exponential cut-off:
\begin{multline*}
 \int_{\R^d\setminus B_{2L}}\int_{B_L} \int_{B_L}\frac{1}{1+|x-z|} \exp(-c\frac{|x-z|}{\sqrt{T}})\frac{1}{1+|x'-z|} \exp(-c\frac{|x'-z|}{\sqrt{T}})dzdxdx'
\\
\lesssim \,\int_{\R^d\setminus B_{L}} L^{4} \frac{1}{1+|z|^{2}}\exp(-c\frac{|z|}{\sqrt{T}})dz
\,\lesssim\, L^4 \ln (1+\frac{\sqrt{T}}{L}),
\end{multline*}
and the claim \eqref{eq:estim-var-hom} follows for $d=2$.

\medskip

\noindent To extend the result to the corrector field itself for $d>2$, we rely on the same soft arguments
as in the discrete case for the limit $T\uparrow \infty$,
and refer the reader to \cite[Proof of Theorem~2.1, Step~8]{Gloria-Otto-09}.

\subsection{Proof of Proposition~\ref{cor:main-2}}

We divide the proof into six steps. 
In the first step we give some preliminary results on the function $\psi_T$ of Lemma~\ref{lem:diff-psi}, which allow us in the second step to reduce the claim of Proposition~\ref{cor:main-2} to an estimate of $\var{\psi_T}$.
The remaining four steps are dedicated to the proof of that estimate.

\medskip
\step{1} Preliminary results.

\noindent By differentiating \eqref{eq:app-corr-eq-proba} wrt $T$ in Remark~\ref{rem:app-corr-eq-proba}, $\psi_T$ solves: For all $\zeta \in \calH^1$, 
\begin{equation}\label{eq:psiT-probab}
\expec{T^{-1} \psi_T \zeta + \DD \zeta \cdot A(0) \DD \psi_T}\,=\,\expec{\phi_T\zeta}.
\end{equation}
Taking $\zeta=\psi_T$ yields the a priori estimate
\begin{equation}\label{eq:fo-gradTpsi1}
T^{-1}\expec{\psi_T^2}+\expec{|\DD \psi_T|^2}\,\lesssim\,\expec{\phi_T\psi_T}.
\end{equation}

\noindent Next, we prove the following formula for the derivative of $\expec{\phi_T\psi_T}$ with respect to $T$:
\begin{equation}\label{eq:fo-gradTpsi2}
|\partial_T \expec{\phi_T\psi_T}|\,=\,|T^{-2}(\var{\psi_T}+2\expec{\psi_T^*\psi_T})|\,\leq \,T^{-2}(2\var{\psi_T}+\var{\psi_T^*}) ,
\end{equation}
where $\psi_T^*$ is the unique weak solution in $\calH^1$ of 
$$
T^{-1}\psi_T^*-\DD\cdot A^*(0)\DD\psi_T^*\,=\,\phi_T.
$$
To this aim we differentiate \eqref{eq:psiT-probab} in its pointwise form with respect to $T$,
$$
T^{-1}\partial_T \psi_T-\DD\cdot A(0)\DD \partial_T\psi_T\,=\,\partial_T\phi_T+T^{-2}\psi_T\,=\,2T^{-2}\psi_T,
$$
which we rewrite as $\partial_T \psi_T\,=\,2T^{-2} (T^{-1}-\DD\cdot A(0)\DD)^{-1} \psi_T$.
Likewise, we  write $\psi_T^*\,=\, (T^{-1}-\DD\cdot A^*(0)\DD)^{-1} \phi_T$. This implies \eqref{eq:fo-gradTpsi2}
 as follows:
\begin{eqnarray*}
\partial_T \expec{\phi_T\psi_T}&=&\expec{\partial_T\phi_T\psi_T}+\expec{\phi_T\partial_T\psi_T}\\
&=&T^{-2}\var{\psi_T}+2T^{-2}\expec{\phi_T  (T^{-1}-\DD\cdot A(0)\DD)^{-1} \psi_T} \\
&=&T^{-2}\var{\psi_T}+2T^{-2}\expec{\big((T^{-1}-\DD\cdot A^*(0)\DD)^{-1}\phi_T\big) \psi_T} \\
&=& T^{-2}(\var{\psi_T}+2\expec{\psi_T^*\psi_T}).
\end{eqnarray*}
Note that the sensitivity estimates for $\psi_T^*$ are identical to the sensitivity estimates for $\psi_T$ in Lemma~\ref{lem:diff-psi} since the distribution of $A^*$ satisfies the same assumption as the one of $A$
(because transposition is a linear local operation).

\medskip
\step{2} Reduction of the proof of Proposition~\ref{cor:main-2} to the proof of 
\begin{equation}\label{eq:estim-var-psiT-a}
\var{\psi_T},\var{\psi_T^*} \,\lesssim \,
\left\{
\begin{array}{rcl}
2\leq d<6&:&\sqrt{T}^{6-d}, \\
d=6&:&\ln \sqrt{T},\\
d>6&:&1.
\end{array}
\right.
\end{equation}
Since $\nabla \ext \phi(0)$ is the weak limit in $\calH$ of $\DD\phi_T$, we have
by lower semi-continuity of the norm, the triangle inequality, the definition of $\psi_T$, and \eqref{eq:fo-gradTpsi1}:
\begin{eqnarray*}
\expec{|\DD\phi_T-\nabla \ext \phi(0)|^2}^{\frac{1}{2}} &\leq &\liminf_{t\uparrow \infty} \expec{\Big|\int_T^t (\partial_\tau \nabla \phi_\tau)d\tau\Big|^2}^{\frac{1}{2}}\\
&\leq &  \int_T^\infty \expec{|\nabla \partial_\tau \phi_\tau|^2}^{\frac{1}{2}}d\tau\\
&=& \int_T^\infty \tau^{-2} \expec{|\nabla \psi_\tau|^2}^{\frac{1}{2}}d\tau \\
&\stackrel{\eqref{eq:fo-gradTpsi1}}{\lesssim} & \int_T^\infty \tau^{-2} \expec{\phi_\tau \psi_\tau}^{\frac{1}{2}}d\tau .
\end{eqnarray*}
To prove \eqref{eq:prop1} it is therefore enough to show that
\begin{equation}\label{eq:fo-gradTpsi3}
0\,\leq\,\expec{\phi_T \psi_T}\,\lesssim\,
\left\{
\begin{array}{rcl}
d=2&:&T, \\
d=3&:&\sqrt{T},\\
d=4&:&\ln T,\\
d>4&:&1.
\end{array}
\right.
\end{equation}
By \eqref{eq:fo-gradTpsi2}, and Young's and Cauchy-Schwarz' inequalities,
for all $T_0\lesssim 1$ and $T\ge T_0$,
\begin{eqnarray*}
\expec{\phi_T \psi_T}&=&\int_{T_0}^T \partial_\tau \expec{\phi_\tau\psi_\tau}d\tau +\expec{\phi_{T_0}\psi_{T_0}}\\
&\leq&\int_{T_0}^T \tau^{-2} (2\var{\psi_\tau}+\var{\psi_\tau^*})d\tau + \expec{\phi_{T_0}^2}^{\frac{1}{2}}\expec{\psi_{T_0}^2}^{\frac{1}{2}},
\end{eqnarray*}
so that \eqref{eq:fo-gradTpsi3} follows from \eqref{eq:estim-var-psiT-a} to bound the integral term, and from \eqref{eq:fo-gradTpsi1}, Cauchy-Schwarz' inequality, and Proposition~\ref{prop:main-1} with $T=T_0\lesssim 1$ to bound the second term.

\medskip

\noindent The rest of the proof is dedicated to the proof of \eqref{eq:estim-var-psiT-a}.
Since the proofs of the estimates of $\var{\psi_T}$ and $\var{\psi_T^*}$ are similar, we only treat the former.

\medskip
\step{3} Proof of
\begin{multline}\label{eq:var-estim-psiT}
\var{\psi_T} \,\lesssim \, \expec{\int_{\R^d}{\ext h}_{T}^2(z,0) \left(\int_{B_6(z)}|\nabla \ext\psi_T(z')|^2dz'+\nu_d(T)\Big(\int_{B_{18}(z)}|\nabla \ext\phi_T(z')|^2dz'+1\Big) \right) dz  }\\
+\expec{\int_{\R^d} \Big(\int_{B_{6}(z)}|\nabla \ext\phi_T(z')|^2dz'+1\Big) \Big(\int_{\R^d}g_{T}(y){\ext h}_{T}(z,y)dy\Big)^2dz}, 
\end{multline}
where ${\ext h}_T$ and $g_T$ are as in \eqref{eq:def-frak-h} and \eqref{eq:ptwise-decay-estim} for $R=2$, respectively, and $\nu_d(T)$ is given by \eqref{eq:def-nud}.

\noindent Since $\psi_T=T^2\partial_T \phi_T$,
one may apply (SG) to $\psi_T$. The claim then follows from \eqref{eq:diff-psi-1} in Lemma~\ref{lem:diff-psi} 
with $R=2$, and Young's inequality.

\medskip
\noindent The first term of the RHS is a nonlinear term since it involves $\ext\psi_T$, whereas the second term is linear.
We estimate these terms separately in Steps~4 and~5.

\medskip
\step{4} Suboptimal estimate of the nonlinear term:
\begin{multline}
\expec{\int_{\R^d}{\ext h}_{T}^2(z,0) \left(\int_{B_{6}(z)}|\nabla \ext\psi_T(z')|^2dz'+\nu_d(T)\Big(\int_{B_{18}(z)}|\nabla \ext\phi_T(z')|^2dz'+1 \Big)\right) dz }\\
\lesssim \,  
\left\{
\begin{array}{rl}
d=2:&T^{2-2\alpha}\ln T,\\
d=3:&T ,\\
d=4:&\ln^2 T ,\\
d>4:&1 ,
\end{array}
\right\}
+
\expec{|\DD \psi_T|^2} \times 
\left\{
\begin{array}{rl}
d=2:&T^{1-2\alpha},\\
d=3:&\sqrt{T} ,\\
d=4:&\ln T ,\\
d>4:&1 ,
\end{array}
\right\}
\label{eq:ap-2}
\end{multline}
where $\alpha>0$ is the H\"older exponent of Lemma~\ref{lem:ptwise-estim}.
Indeed, for $|z|>3$ we bound ${\ext h}_{T}(z,0)$ by Lemma~\ref{lem:ptwise-estim}, whereas  ${\ext h}_{T}(z,0)$
is of order one for $|z|\leq 3$ by  \eqref{eq:def-frak-h}, so that
\begin{equation*}
{\ext h}_{T}(z,0) \,\lesssim\, 
 \left\{ 
\begin{array}{lll}
d=2&:&\min \left\{|z|^{-\alpha}\exp(-c\frac{|z|}{\sqrt{T}} ),1\right\},\\
d>2&:&\min \left\{{|z|^{2-d}}\exp(-c\frac{|z|}{\sqrt{T}} ),1\right\}.
\end{array}
\right\}
\end{equation*}
Since this estimate is deterministic, one may take it out of the expectation in
the LHS of \eqref{eq:ap-2}.
By stationarity and Lemma~\ref{lem:app-corr},
\begin{multline*}
\expec{\int_{B_{6}(z)}|\nabla \ext\psi_T(z')|^2dz'+\nu_d(T)\Big(\int_{B_{18}(z)}|\nabla \ext\phi_T(z')|^2dz'+1\Big)}\\ \lesssim \, 
\expec{|\DD \psi_T|^2}+\nu_d(T)\Big(\expec{|\DD \phi_T|^2}+1\Big) \,\lesssim \, \expec{|\DD \psi_T|^2}+\nu_d(T).
\end{multline*}
Estimate \eqref{eq:ap-2} thus follows by integrating over $z$ and by the definition~\eqref{eq:def-nud} of $\nu_d(T)$.

\medskip
\step{5} Estimate of the linear term: 
\begin{multline}
\expec{\int_{\R^d} \Big(\int_{B_{6}(z)}|\nabla \ext\phi_T(z')|^2dz'+1\Big) \Big(\int_{\R^d}g_{T}(y){\ext h}_{T}(z,y)dy\Big)^2dz}
\\
\lesssim 
\left\{
\begin{array}{rl}
2\leq d<6:&\sqrt{T}^{6-d}, \\
d=6:&\ln \sqrt{T} ,\\
d>6:&1 .
\end{array}
\right.
\label{eq:ap-3}
\end{multline}
By the triangle inequality in probability,
\begin{multline*}
\expec{\int_{\R^d} \Big(\int_{B_{6}(z)}|\nabla \ext\phi_T(z')|^2dz'+1\Big) \Big(\int_{\R^d}g_{T}(y){\ext h}_{T}(z,y)dy\Big)^2dz}
\\
=\,\int_{\R^d} \expec{\Big(\int_{\R^d}g_{T}(y){\ext h}_{T}(z,y)\Big(\int_{B_{6}(z)}|\nabla \ext\phi_T(z')|^2dz'+1\Big)^{\frac{1}{2}}dy\Big)^2}dz
\\
\leq\, \int_{\R^d}\bigg(\int_{\R^d}g_{T}(y)\expec{\Big(\int_{B_{6}(z)}|\nabla \ext\phi_T(z')|^2dz'+1\Big){\ext h}^2_{T}(z,y)}^\frac{1}{2}dy\bigg)^2dz.
\end{multline*}
Let $\bar p$ be the Meyers exponent of Lemma~\ref{lem:annealed-estim}.
H\"older's inequality in probability with exponents $(\frac{\bar p}{\bar p-1},\bar p)$, Lemma~\ref{lem:annealed-estim}, Proposition~\ref{prop:main-1}, and the definition of $g_T$ then yield
\begin{multline*}
\expec{\int_{\R^d} \Big(\int_{B_{6}(z)}|\nabla \ext\phi_T(z')|^2dz'+1\Big) \Big(\int_{\R^d}g_{T}(y){\ext h}_{T}(z,y)dy\Big)^2dz}
\\
\leq\, \int_{\R^d}\bigg(\int_{\R^d}
\left\{\begin{array}{ccc}
\ln(2+\frac{\sqrt{T}}{|y|})&\mbox{for}&d=2\\
|y|^{2-d}&\mbox{for}&d>2
\end{array}\right\}\frac{1}{1+|y-z|^{d-1}}\exp(-c\frac{|y|+|y-z|}{\sqrt{T}})dy\bigg)^2dz.
\end{multline*}
For $2\le d \le 3$ we use the exponential cut-off both in the inner and outer integrals (dimension $d=3$ is critical for the inner integral),
for $3< d\le 6$ we use the exponential cut-off for the outer integral only (dimension $d=6$ is critical for the outer integral), 
and for $d>6$ one may discard the exponential cut-off.
We start with $d>3$:
\begin{multline*}
\expec{\int_{\R^d} \Big(\int_{B_{6}(z)}|\nabla \ext\phi_T(z')|^2dz'+1\Big) \Big(\int_{\R^d}g_{T}(y){\ext h}_{T}(z,y)dy\Big)^2dz}
\\
\lesssim\, \int_{\R^d} \frac{1}{1+|z|^{2(d-3)}}\exp(-c\frac{|z|}{\sqrt{T}})dz\,\lesssim 
\left\{
\begin{array}{rl}
3< d<6:&\sqrt{T}^{6-d}, \\
d=6:&\ln \sqrt{T} ,\\
d>6:&1 .
\end{array}
\right.
\end{multline*}
For $d=2$, the inner integral scales as $\sqrt{T}\exp(-c\frac{|z|}{\sqrt{T}})$, and the claim \eqref{eq:ap-3} follows.
For $d=3$, we rewrite the inner integrand using the exponential cut-off (up to changing the value of $c$) in the form
\begin{multline*}
\expec{\int_{\R^3} \Big(\int_{B_{6}(z)}|\nabla \ext\phi_T(z')|^2dz'+1\Big) \Big(\int_{\R^3}g_{T}(y){\ext h}_{T}(z,y)dy\Big)^2dz}
\\
\lesssim \, \int_{\R^3}\bigg(\int_{\R^3}\sqrt{T}|y|^{-2}\frac{1}{1+|y-z|^{2}}\exp(-c\frac{|y|+|y-z|}{\sqrt{T}})dy\bigg)^2dz
\\
\lesssim \, \int_{\R^3} T \frac{1}{1+|z|^{2}}\exp(-c\frac{|z|}{\sqrt{T}})dz\,\lesssim \, \sqrt{T}^3,
\end{multline*}
that is, \eqref{eq:ap-3}.

\medskip

\step{6} Nonlinear estimate and buckling.

\noindent The combination of \eqref{eq:var-estim-psiT}  with \eqref{eq:ap-2} \& \eqref{eq:ap-3}  yields
\begin{equation}\label{eq:pr-th2-step7-1}
\var{\psi_T} \,\lesssim \,
\expec{|\DD \psi_T|^2}\times \left\{
\begin{array}{rl}
d=2:&T^{1-2\alpha}\\
d=3:&\sqrt{T}\\
d=4:&\ln T\\
d>4:&1
\end{array}
\right\}
+
\left\{
\begin{array}{rl}
2\leq d<6:&\sqrt{T}^{6-d}\\
d=6:&\ln \sqrt{T}\\
d>6:&1
\end{array}
\right\}.
\end{equation}
We then appeal to the following nonlinear estimate, which follows from \eqref{eq:fo-gradTpsi1}, Cauchy-Schwarz' inequality
and Proposition~\ref{prop:main-1}:
\begin{equation}\label{eq:syst-subopt}
\expec{|\DD \psi_T|^2}\,\lesssim \,
\var{\psi_T}^{\frac{1}{2}} \left\{
\begin{array}{rcl}
d=2&:&(\ln T)^{\frac{1}{2}},\\
d>2&:&1.
\end{array}
\right.
\end{equation}

\noindent Combined with the nonlinear estimate~\eqref{eq:syst-subopt}, \eqref{eq:pr-th2-step7-1} thus turns into
\begin{equation*}
\var{\psi_T} \,\lesssim \,
\var{\psi_T}^{\frac{1}{2}}\times \left\{
\begin{array}{rl}
d=2:&T^{1-2\alpha}(\ln T)^{\frac{1}{2}}\\
d=3:&\sqrt{T}\\
d=4:&\ln T\\
d>4:&1
\end{array}
\right\}
+
\left\{
\begin{array}{rl}
2\leq d<6:&\sqrt{T}^{6-d}\\
d=6:&\ln \sqrt{T}\\
d>6:&1
\end{array}
\right\},
\end{equation*}
which yields the desired estimate \eqref{eq:estim-var-psiT-a} for $\psi_T$ by Young's inequality.


\subsection{Proof of Theorem~\ref{th:main-2}}

Theorem~\ref{th:main-2} follows from
the identities
\begin{eqnarray}
\expec{(\DD \phi_T'-\nabla \ext \phi'(0))\cdot A(0)(\xi+\nabla \ext \phi(0))}&=&0,\label{eq:pr-Th2-1}\\
\expec{(\DD \phi_T-\nabla \ext \phi(0))\cdot A^*(0)(\xi'+\nabla \ext \phi'(0))}&=&0,\label{eq:pr-Th2-2}
\end{eqnarray}
the calculation
\begin{eqnarray*}
\lefteqn{\xi'\cdot A_T\xi-\xi'\cdot A_\ho \xi}\\
&=& \expec{(\xi'+\DD \phi_T')\cdot A(0)(\xi+\DD \phi_T)}
-\expec{(\xi'+\nabla \ext \phi'(0))\cdot A(0)(\xi+\nabla \ext \phi(0))} \\
&=&\expec{(\DD \phi_T'-\nabla \ext \phi'(0))\cdot A(0)(\xi+\nabla \ext \phi(0))}
+\expec{(\xi'+\DD \phi_T')\cdot A(0)(\DD \phi_T-\nabla \ext \phi(0))}\\
&\stackrel{\eqref{eq:pr-Th2-1},\eqref{eq:pr-Th2-2}}{=}&-\expec{(\DD \phi_T-\nabla \ext \phi(0))\cdot A^*(0)(\xi'+\nabla \ext \phi'(0))}
+\expec{(\xi'+\DD \phi_T')\cdot A(0)(\DD \phi_T-\nabla \ext \phi(0))}\\
&=&\expec{(\DD \phi_T'-\nabla \ext \phi'(0))\cdot A(0)(\DD \phi_T-\nabla \ext \phi(0))},
\end{eqnarray*}
Cauchy-Schwarz' inequality, and Proposition~\ref{cor:main-2} (which holds both for $\phi_T$ and $\phi_T'$).

\subsection{Proof of Corollary~\ref{coro:spectr-bottom}}

By the estimate \eqref{eq:estim-var-psiT-a} of $\var{\psi_T}$ in the proof of Proposition~\ref{cor:main-2}, it is enough
to prove that for all $T\geq 1$
\begin{equation}\label{eq:pr-estim-bottom-1}
\expec{\mathfrak{d}G(d\lambda)\mathfrak{d}}([0,T^{-1}])\,\lesssim \, T^{-4} \var{\psi_T}.
\end{equation}
Indeed, since for all $\lambda \leq T^{-1}$,
$$
\frac{T^{-4}}{(T^{-1}+\lambda)^4} \,\gtrsim\, 1,
$$
we have
\begin{equation}\label{eq:pr-estim-bottom-1b}
\expec{\mathfrak{d}G(d\lambda)\mathfrak{d}}([0,T^{-1}])\,=\, \int_0^{T^{-1}}d e_{\mathfrak{d}}(\lambda)\,\lesssim \,T^{-4}\int_0^{+\infty} \frac{1}{(T^{-1}+\lambda)^4}d e_{\mathfrak{d}}(\lambda).
\end{equation}
Since $\psi_T=(T^{-1}+\calL)^{-2}\mathfrak{d}$, 
we recognise in the integral of the RHS of \eqref{eq:pr-estim-bottom-1b} 
the spectral representation of $\expec{\psi_T^2}=\var{\psi_T}$, which proves \eqref{eq:pr-estim-bottom-1}.


\section{Proofs of the spectral gap estimates}

\subsection{Proof of Lemma~\ref{lem:SG-ergo}}

We shall prove ergodicity in the following form: For all $X\in L^1(\Omega)$, 
we have
\begin{equation}\label{eq:SG-ergo-0}
\lim_{R\uparrow \infty} \expec{\Big|\fint_{B_R} \overline{X}(y)dy-\expec{X}\Big|}\,=\,0,
\end{equation}
where $\ext X$ is the stationary extension of $X$.
We divide the proof into two steps. We first show by approximation that it is enough to prove \eqref{eq:SG-ergo-0} for
bounded random fields $X\in L^\infty(\Omega)$ which only depend on $A$ through its restriction $A|_V$ on some bounded set $V$.
We then show that for such random fields, \eqref{eq:SG-ergo-0} follows from (SG).

\medskip

\step{1} Approximation argument. 

\noindent 
Since the map $X\mapsto \langle|\fint_{B_R} \overline{X}(y)dy-\langle X\rangle|\rangle$ is Lipschitz continuous on $L^1(\Omega)$ uniformly in $R$, it is enough to establish \eqref{eq:SG-ergo-0} on an $L^1(\Omega)$-dense
subset of $X$'s. By definition of measurability, we thus may restrict ourselves to $X$'s that
depend on $A$ only through its restriction $A_{|V}$ on some bounded set $V$. Moreover, a simple
truncation argument shows that any $X\in L^1(\Omega)$ can be approximated in $L^1(\Omega)$ by
$\tilde X\in L^\infty(\Omega)$. Hence we may restrict ourselves to $X\in L^\infty(\Omega)$
that depend on $A$ only through its restriction on a some ball $B_L$.

\medskip

\step{2} Proof that (SG) implies \eqref{eq:SG-ergo-0}.

\noindent By Step~1, it is enough to prove \eqref{eq:SG-ergo-0} for bounded random fields $X\in L^\infty(\Omega)$ which only depend on $A$ through its restriction to a bounded set $B_L$.
In that case, by stationarity of $\ext X$ and since $X$ does not depend on $A|_{\R^d\setminus B_L}$,
\begin{eqnarray*}
\var{\fint_{B_R} \ext X(y)dy} &\leq & \frac{1}{\rho}\int_{\R^d} \expec{\Big(\osc{A|_{B_\ell(x)}}{\fint_{B_R}\ext X(y)dy}\Big)^2}dx\\
&\leq &\frac{1}{\rho} \fint_{B_R}\fint_{B_R} \int_{\R^d} \expec{\osc{A|_{B_\ell(x)}}{\ext X(y)}\osc{A|_{B_\ell(x)}}{\ext X(y')}}dxdydy' \\
&\lesssim & \|X\|_{L^\infty(\Omega)}^2 \fint_{B_R}\fint_{B_R} \int_{\R^d} \mathds{1}_{|x-y|\leq L+\ell}\mathds{1}_{|x-y'|\leq L+\ell}dxdydy'\\
&\lesssim & R^{-d} (L+\ell)^{2d} \|X\|_{L^\infty(\Omega)}^2,
\end{eqnarray*}
so that by Cauchy-Schwarz' inequality and stationarity of $\ext X$,
$$
\expec{\Big|\fint_{B_R} \overline{X}(y)dy-\expec{X}\Big|}\,\leq \,\var{\fint_{B_R} \ext X(y)dy} \,\stackrel{R\uparrow \infty}{\longrightarrow} 0.
$$

\subsection{Proof of Corollary~\ref{coro:var-estim}}

We assume w.~l.~o.~g. that $\expec{X}=0$, and divide the proof into three steps.

\medskip
\step{1} Proxy for the Leibniz rule: For any function $\zeta$ and all $q\geq 1$,
\begin{equation}\label{eq:Leib-osc}
\osc{}{|\zeta|^{q}}\,\lesssim \, |\zeta|^{q-1}\osc{}{\zeta}+\left(\osc{}{\zeta}\right)^{q}.
\end{equation}
This follows from Young's inequality and the two elementary estimates
\begin{eqnarray*}
\osc{}{|\zeta|^{q}}&\lesssim &\Big(\sup |\zeta|^{q-1} \Big)\osc{}{\zeta},\\
\sup |\zeta| &\leq&|\zeta|+\osc{}{\zeta}.
\end{eqnarray*}

\medskip

\step{2} Proof that for all $q\ge 1$,
\begin{equation}\label{eq:pr-coro-sg-2-1}
\expec{X^{2q}}^{\frac{1}{q}}\,\lesssim \, \expec{X^2}+\expec{\Big( \int_{\R^d} \Big(\osc{A|_{B_{2\ell}(z)}}{X}\Big)^2dz\Big)^q}^{\frac{1}{q}}.
\end{equation}
By definition of the oscillation we have $\osc{A|_{B_{2\ell}(z)}}{X}\ge \osc{A|_{B_{\ell}(0)}}{X}$ for all $z\in B_\ell(0)$
so that $\int_{\mathbb{R}^d}(\osc{A|_{B_{2\ell}(z)}}{X})^2dz\gtrsim(\osc{A|_{B_\ell(0)}}{X})^2$. Since the origin
plays no special role, this can be rewritten as
\begin{equation}\nonumber
\sup_{z}\Big(\osc{A|_{B_\ell(z)}}{X}\Big)^2\lesssim \int_{\mathbb{R}^d}\Big(\osc{A|_{B_{2\ell}(z)}}{X}\Big)^2dz.
\end{equation}
This immediately implies for any $q\ge 1$
\begin{equation}\label{eq:disc-lq-l1}
\int_{\mathbb{R}^d}\Big(\osc{A|_{B_\ell(z)}}{X}\Big)^{2q}dz\lesssim\left(\int_{\mathbb{R}^d}\Big(\osc{A|_{B_{2\ell}(z)}}{X}\Big)^2dz\right)^q.
\end{equation}
We then apply (SG) to $|X|^q$:
\begin{eqnarray*}
\expec{X^{2q}} -\expec{|X|^q}^2\,=\,\var{|X|^q}
&\lesssim & \expec{ \int_{\R^d} \Big( \osc{A|_{B_{\ell}(z)}}{|X|^q}\Big)^2dz} .
\end{eqnarray*}

\medskip

\noindent By the Leibniz rule \eqref{eq:Leib-osc} this implies
\begin{equation}\label{eq:pr-coro-sg-2-2}
\expec{X^{2q}}\,\lesssim \, \expec{|X|^q}^2+\expec{ \int_{\R^d} X^{2(q-1)} \Big(\osc{A|_{B_{\ell}(z)}}{X}\Big)^2} + \expec{ \int_{\R^d} \Big(\osc{A|_{B_{\ell}(z)}}{X}\Big)^{2q} }.
\end{equation}
We treat the three terms of the RHS separately.
For the third term we appeal to \eqref{eq:disc-lq-l1}.
For the second term, we use H\"older's and Young's inequalities both with exponents $(\frac{q}{q-1},q)$,
which yields for all $C>0$
\begin{eqnarray}
\lefteqn{\expec{X^{2(q-1)} \int_{\R^d}   \Big(\osc{A|_{B_{\ell}(z)}}{X}\Big)^2} }\nonumber
\\
&\leq & \expec{X^{2q}}^{\frac{q-1}{q}} \expec{\Big(\int_{\R^d}  \Big(\osc{A|_{B_{\ell}(z)}}{X}\Big)^2\Big)^q}^{\frac{1}{q}}\nonumber \\
&\leq & \frac{q-1}{Cq}\expec{X^{2q}}+\frac{C^{q-1}}{q}\expec{\Big(\int_{\R^d}  \Big(\osc{A|_{B_{\ell}(z)}}{X}\Big)^2\Big)^q}.
\label{eq:pr-coro-sg-2-4}
\end{eqnarray}
\noindent For the first term of the RHS of \eqref{eq:pr-coro-sg-2-2} it is enough to treat the case $q>2$ since for $q\leq 2$, it is controlled by the $q^{th}$ power of the RHS of \eqref{eq:pr-coro-sg-2-1}.
We then apply H\"older's inequality with exponents $(2\frac{q-1}{q-2},2\frac{q-1}{q})$ to
$\expec{|X|^{q}}=\expec{|X|^{q\frac{q-2}{q-1}}|X|^{\frac{q}{q-1}}}$:
This yields for all $C>0$ using then Young's inequality 
\begin{eqnarray}
 \expec{|X|^q}^2 &\leq & \expec{X^{2q}}^{\frac{q-2}{q-1}}\expec{X^2}^{\frac{q}{q-1}} \nonumber \\
&\leq & \frac{q-2}{C(q-1)}\expec{X^{2q}}+\frac{C^{q-2}}{q-1}\expec{X^2}^q.
\label{eq:pr-coro-sg-2-5}
\end{eqnarray}
The combination of \eqref{eq:pr-coro-sg-2-2}---\eqref{eq:pr-coro-sg-2-5} with \eqref{eq:disc-lq-l1} and Young's inequality
yields  \eqref{eq:pr-coro-sg-2-1}.

\medskip
\step{3} Conclusion.

\noindent The spectral gap estimate applied to $X$
\begin{equation*}
\expec{X^{2}} \,\lesssim \, \expec{ \int_{\R^d} \Big( \osc{A|_{B_\ell(z)}}{X}\Big)^2dz},
\end{equation*}
combined with Jensen's inequality in probability yields
\begin{equation*}
\expec{X^{2}} \,\lesssim \, \expec{ \Big(\int_{\R^d} \Big( \osc{A|_{B_\ell(z)}}{X}\Big)^2dz\Big)^q}^{\frac{1}{q}},
\end{equation*}
so that the claim follows from \eqref{eq:pr-coro-sg-2-1}.


\section{Proofs of the sensitivity estimates}

\noindent The sensitivity estimates do not require the coefficients $A\in \Omega$ to be smooth.
It is however convenient to first prove these estimates under that additional assumption.
These estimates are then recovered for general coefficients by density.
Indeed, by elementary $L^2$-theory, if the coefficients $A$ are approximated by a sequence of smooth coefficients $A_k$ in $L^1_\loc(\R^d)$, then $\ext\phi_T(\cdot;A_k)$
converges in $H^1_\loc(\R^d)$ to $\phi_T(\cdot;A)$, and for all $x$ the Green function $y\mapsto G_T(y,x;A_k)$ converges in $H^1_\loc(\R^d\setminus B_r(x))$ for all $r>0$.
This is enough to prove the convergence of the RHS of the oscillation estimates \eqref{eq:diff-phi-1} and \eqref{eq:diff-psi-1}.
For the LHS we use in addition that $\ext\phi_T$ and $\ext\psi_T$ are H\"older continuous uniformly in space and with respect to $A$, so that $L^2_\loc(\R^d)$ convergence implies pointwise convergence.
The  H\"older continuity of $\ext\phi_T$ is a consequence of the De Giorgi-Nash-Moser theory, while the \emph{uniform} H\"older continuity in addition relies on the uniform $L^2$-bound \eqref{eq:uniqueness-cond-refined} of Lemma~\ref{lem:def-phiT-conv}. A similar argument holds for $\ext \psi_T$.

\subsection{Proof of Lemma~\ref{lem:diff-phi}}

We let $A_1, A_2\in \Omega$ be  smooth and coincide outside $B_R(z)$, $z\in \R^d$, with some $A\in \Omega$.
For convenience we denote by $\ext \phi_1$ and $\ext \phi_2$, and $G_1$ and $G_2$ the associated 
modified correctors for $\xi\in \R^d$, $|\xi|=1$, and Green functions for $T>0$.

\medskip
\step{1} Preliminaries.

\noindent By definition, $\ext{ \phi}_1$ and $\ext\phi_2$ are smooth and $\ext{\phi}_1-\ext\phi_2$ is a classical solution of 
\begin{equation}\label{eq:pr-diff-phi-2}
T^{-1}(\ext\phi_1-\ext{\phi}_2)-\nabla \cdot (A_1 \nabla (\ext\phi_1-\ext{\phi}_2))\,=\, \nabla \cdot ((A_1-A_2)(\xi+\nabla \ext{ \phi}_2)).
\end{equation}
Since $A_1$ and $A_2$ coincide outside $B_R$, the RHS of \eqref{eq:pr-diff-phi-2} has compact support so that 
$\ext\phi_1-\ext{ \phi}_2\in H^1(\R^d)$.
Since all the quantities are smooth and $x\mapsto  G_1(x,y) \in W^{1,d/(d-1+\e)}(\R^d)$ for all $0<\e\leq 1$, $\ext\phi_1-\ext{ \phi}_2$ satisfies the Green representation formula
\begin{equation}\label{eq:pr-osc-phi1}
\ext\phi_1(x)-\ext{\phi}_2(x)\,=\,-\int_{\R^d}\nabla_y G_1(y,x)\cdot (A_1(y)-A_2(y))(\xi+\nabla \ext{ \phi}_2(y))dy.
\end{equation}
\medskip

\noindent The second ingredient to the proof is estimate \eqref{eq:sup-phi-1}, which we prove now.
Since $\ext{\phi}_1-\ext\phi_2\in H^1(\R^d)$, an a priori estimate based on \eqref{eq:pr-diff-phi-2}
yields
\begin{equation*}
\int_{\R^d}|\nabla \ext\phi_1(y)-\nabla \ext{ \phi}_2(y)|^2dy \,\lesssim \,\int_{B_R(z)}|\xi+\nabla \ext{\phi}_2(y)|^2dy.
\end{equation*}
This shows by the triangle inequality that for all $x\in\R^d$
\begin{equation*}
\int_{B_R(x)}|\xi+\nabla  \ext\phi_1(y)|^2dy \,\lesssim \,\int_{B_R(x)}|\xi+\nabla \ext{\phi}_2(y)|^2dy+\int_{B_R(z)}|\xi+\nabla \ext{ \phi}_2(y)|^2dy,
\end{equation*}
which yields the claim.

\medskip
\step{2} Proof of \eqref{eq:diff-phi-1} for $|z-x|\geq 2R$.

\noindent The starting point is the Green representation formula \eqref{eq:pr-osc-phi1}, which yields by Cauchy-Schwarz' inequality:
\begin{equation}
|\ext\phi_1(x)-\ext{ \phi}_2(x)|
\,\lesssim \, \left(\int_{B_R(z)} |\nabla_y  G_1(y,x)|^2dy\right)^{\frac{1}{2}}  \left(\int_{B_R(z)} |\xi+\nabla  \ext{\phi}_2(y)|^2dy\right)^{\frac{1}{2}} . \label{eq:pr-diff-phi-1}
\end{equation}
In order to conclude, we need to take the supremum over all the smooth coefficients $A_1,A_2$ such that ${A_1}|_{B_R(z)}={A_2}|_{B_R(z)}={A}|_{B_R(z)}$ in the RHS of \eqref{eq:pr-diff-phi-1}.
For the first term, which depends on $A_1$ but not on $A_2$,  we appeal to Lemma~\ref{lem:diff-Green}.
For the second term, we use
\eqref{eq:sup-phi-1}.
Note that  the RHS of \eqref{eq:pr-diff-phi-1} is finite only for $|z-x|>R$, and that $\int_{B_R(z)} |\nabla_y  G_1(y,x)|^2dy \, \lesssim \,1$ for all $|z-x|\geq 2R$
by property~\eqref{gu.11} of Definition~\ref{def:Green}.

\medskip 
\step{3} Proof of \eqref{eq:diff-phi-1} for $|z-x|< 2R$.

\noindent By definition of the oscillation $\osc{A|_{B_R(z)}}{}$,
using the triangle inequality and thus just at the expense of a factor of two,
we may make any restrictions on one of the two coefficient fields $A_1$ and $A_2$, 
say on $A_2$, provided it does not violate its smooth connection to $A$ outside of $B_R(z)$.
For our purpose, it is convenient to have {\it quantitative} smoothness of $A_2$ near $z$ in form of
\begin{equation}\label{eq:choiceA}
{A_2}|_{\R^d\setminus B_{R}(z)}\,=\,A|_{\R^d\setminus B_{R}(z)},  \quad {A_2}|_{B_{\frac{R}{2}}(z)}\,=\,\Id.
\end{equation}
As mentioned above, this can be obtained by setting $A_2=(1-\eta)A+\eta\Id$,
where $\eta$ is a smooth cut-off function for $B_{\frac{R}{2}}(z)$ in $B_R(z)$.

\medskip
\noindent We turn now to the proof of \eqref{eq:diff-phi-1}.
It is enough to prove that for all $R\lesssim 1$ and all $|z-x|\leq \frac{2}{3}R$,
$$
\osc{A_1|_{B_R(z)}}{\ext\phi_1(x)}\,\lesssim \, \left(\int_{B_{R}(z)}|\nabla \ext \phi_1(y)|^2dy+1\right)^{\frac{1}{2}},
$$
then to replace $R$ by $3R$ in this estimate, and to use that $$\dps\osc{A_1|_{B_R(z)}}{\ext\phi_1(x)}\,\leq \,\dps\osc{A_1|_{B_{3R}(z)}}{\ext\phi_1(x)}.$$

\noindent Due to the singularity of the Green function at $x=y$, the estimate \eqref{eq:pr-diff-phi-1} of Step~2 cannot be used for $|z-x|\leq R$.
Instead of using Cauchy-Schwarz' inequality, we cut the integral into two parts $B_{\frac{R}{4}}(z)$ and $B_R(z)\setminus B_{\frac{R}{4}}(z)$, and use H\"older's inequality with exponents $(p,q)$ for some $1<p<\frac{d}{d-1}$  on the first term, and Cauchy-Schwarz' inequality on the second term:
\begin{eqnarray}
\lefteqn{|\ext\phi_1(x)-\ext{\phi}_2(x)|}\nonumber \\
&\lesssim &\left( \int_{B_{\frac{R}{4}}(z)} |\nabla_y   G_1(y,x)|^pdy\right)^{\frac{1}{p}} \left(\int_{B_{\frac{R}{4}}(z)} |\xi+\nabla \ext{\phi}_2(y)|^qdy\right)^{\frac{1}{q}}
\nonumber\\
& &+  \left(\int_{B_R(z)\setminus B_{\frac{R}{4}}(z)} |\nabla_y  G_1(y,x)|^2dy\right)^{\frac{1}{2}}  
 \left(\int_{B_R(z)\setminus B_{\frac{R}{4}}(z)} |\xi+\nabla  \ext{ \phi}_2(y)|^2dy\right)^{\frac{1}{2}} . \label{eq:pr-lem-osc1}
\end{eqnarray}
We first treat the first summand on the RHS: The first factor is bounded uniformly in $A_1$ since $\nabla G_1 $ is bounded in $L^{p}(\R^d)$ uniformly with respect to $A_1$ and $T>0$ (as a consequence of property~\eqref{gu.11} in Definition~\ref{def:Green} and a dyadic decomposition of $B_{\frac{R}{4}}$). For the second factor, we note that $\ext{\phi}_2$ satisfies 
$$
T^{-1} \ext{\phi}_2-\triangle \ext{\phi}_2\,=\,0
$$
in $B_{\frac{R}{2}}(z)$ since $A_2|_{B_{\frac{R}{2}}(z)}=\Id$, so that for all $i\in \{1,\dots,d\}$, $\partial_{x_i} \ext{\phi}_2$
satisfies 
$$
T^{-1}  \partial_{x_i} \ext{\phi}_2-\triangle \partial_{x_i}\ext{\phi}_2\,=\,0
$$
in $B_{\frac{R}{2}}(z)$.
Hence, by classical interior elliptic regularity (see for instance \cite[Theorem~2, Sec.~6.3]{Evans-10}), for all $k\in \N$,
\begin{equation}\label{eq:interior-reg-phiT}
\|\nabla \ext{\phi}_2\|_{H^k(B_{\frac{R}{4}}(z))} \,\lesssim \, \|\nabla \ext{\phi}_2\|_{L^2(B_{\frac{R}{2}}(z))},
\end{equation}
where the multiplicative constant depends on $k$ and $R$.
This yields by Sobolev embedding
\begin{equation*}
\left(\int_{B_{\frac{R}{4}}(z)} |\nabla \ext{\phi}_2(y)|^qdy\right)^{\frac{1}{q}} \,\lesssim \,\left(\int_{B_{\frac{R}{2}}(z)} |\nabla \ext{ \phi}_2(y)|^2dy\right)^{\frac{1}{2}},
\end{equation*}
so that the first summand on the RHS of \eqref{eq:pr-lem-osc1} is estimated by 
\begin{multline*}
\left( \int_{B_{\frac{R}{4}}(z)} |\nabla_y   G_1(y,x)|^pdy\right)^{\frac{1}{p}} \left(\int_{B_{\frac{R}{4}}(z)} |\xi+\nabla \ext{\phi}_2(y)|^qdy\right)^{\frac{1}{q}}
\\ \,\lesssim \, \left(\int_{B_{R}(z)}|\nabla \ext{\phi}_2(y)|^2dy+1\right)^{\frac{1}{2}}.
\end{multline*}
For the second summand of the RHS of \eqref{eq:pr-lem-osc1}, the first factor is of order 1 by property~\eqref{gu.11} in Definition~\ref{def:Green}, so that 
\begin{multline*}
\left(\int_{B_R(z)\setminus B_{\frac{R}{4}}(z)} |\nabla_y  G_1(y,x)|^2dy\right)^{\frac{1}{2}}  
 \left(\int_{B_R(z)\setminus B_{\frac{R}{4}}(z)} |\xi+\nabla  \ext{\phi}_2(y)|^2dy\right)^{\frac{1}{2}}\\\,\lesssim \, \left(\int_{B_{R}(z)}|\nabla \ext{\phi}_2(y)|^2dy+1\right)^{\frac{1}{2}}.
\end{multline*}
The claim then follows from taking the supremum in $A_2|_{B_{R}(z)}$ of these two estimates using \eqref{eq:sup-phi-1}.

\subsection{Proof of Lemma~\ref{lem:diff-psi}}

The proof has the same structure as the proof of Lemma~\ref{lem:diff-phi}.
We let $A_1, A_2\in \Omega$ be  smooth and coincide outside $B_R(z)$, $z\in \R^d$, with some $A\in \Omega$.
For convenience we denote by $\ext \phi_1$ and $\ext \phi_2$, and $\ext\psi_1$ and $\ext{\psi}_2$ the associated 
modified correctors for $\xi\in \R^d$, $|\xi|=1$ and  the functions given by \eqref{eq:diff-psi}.
We also denote by $G_1$ the Green function associated with $A_1$ and a zero-order term of magnitude $T$.

\medskip
\step{1} Preliminaries and proof of \eqref{eq:sup-psi-1}.

\noindent 
Since $\delta \ext \phi:=\ext \phi_1-\ext{ \phi}_2$ is smooth and in $H^{1}(\R^d)$ according to Step~1 in the proof of Lemma~\ref{lem:diff-phi},
and $\ext\psi_2$ is smooth, the function $\delta \ext\psi:= \ext\psi_1-\ext{\psi}_2$ is a classical solution of
\begin{equation}\label{eq:pr-diff-psi-1.0}
T^{-1}\delta \ext\psi-\nabla \cdot  A_1 \nabla \delta \ext\psi \,=\,\delta \ext\phi-\nabla \cdot (A_1-A_2)\nabla  \ext{ \psi}_2,
\end{equation}
and is in $H^1(\R^d)$.
Hence  the Green representation formula holds: For all $x\in \R^d$
\begin{equation}\label{eq:GRF-psi}
\delta \ext\psi(x)\,=\,\int_{\R^d}\nabla_y G_1(y,x) \cdot (A_1-A_2)(y)\nabla \ext{\psi}_2(y)dy+\int_{\R^d} G_{1}(y,x)\delta \ext\phi(y)dy.
\end{equation}

\medskip
\noindent We first establish \eqref{eq:sup-psi-1}. We test the following equivalent form of \eqref{eq:pr-diff-psi-1.0} 
\begin{equation*}
T^{-1}\delta \ext\psi-\nabla \cdot  A_2 \nabla \delta \ext\psi \,=\,\delta \ext\phi-\nabla \cdot (A_2-A_1)\nabla  \ext{\psi}_1
\end{equation*}
with $\delta \ext\psi$, which yields the a priori estimate
\begin{equation}\label{hardy}
T^{-1}\int_{\R^d} (\delta \ext\psi)^2dx+\int_{\R^d} |\nabla\delta \ext\psi|^2dx
\,\lesssim \, \int_{\R^d} |\delta \ext\phi\delta\ext \psi| dx+\int_{B_R(z)} |\nabla \ext{\psi}_1|^2dx .
\end{equation}
We then appeal to \eqref{eq:diff-phi-1} in Lemma~\ref{lem:diff-phi} to bound the first term of the RHS by
$$
\int_{\R^d} |\delta \ext\phi\delta \ext\psi|dx \,\leq \, \int_{\R^d} |\delta \ext\psi(x)| {\ext{h}}_{1}(z,x) dx
\left(\int_{B_{3R}(z)}|\nabla \ext{\phi}_1|^2dy+1\right)^{\frac{1}{2}},
$$
where ${\ext{h}}_{1}$ is given by \eqref{eq:def-frak-h} (with Green's function $G_T(\cdot,\cdot;A_1)$).
In dimension $d=2$, we use Cauchy-Schwarz' inequality, and obtain by integrating $ {\ext{h}}_{1}$ on dyadic annuli (and using~\eqref{gu.11} in Definition~\ref{def:Green}):
\begin{equation*}
\int_{\R^d} |\delta \ext\phi\delta \ext\psi|dx \,\leq \, \sqrt{T}\left(T^{-1}\int_{\R^d} \delta \ext\psi^2dx\right)^{\frac{1}{2}} (\ln T)^{\frac{1}{2}} \left(\int_{B_{3R}(z)}|\nabla \ext{\phi}_1|^2dy+1\right)^{\frac{1}{2}}.
\end{equation*}
Using Young's inequality and absorbing the $L^2$-norm of $\delta \ext\psi$ into the LHS of \eqref{hardy} yield
\begin{equation*}
\int_{\R^d} |\nabla\delta \ext\psi|^2dx \,\lesssim \,  \int_{B_R(z)} |\nabla \ext{\psi}_1|^2dy+T\ln T \Big(\int_{B_{3R}(z)}|\nabla \ext{\phi}_1|^2dy+1\Big),
\end{equation*}
from which \eqref{eq:sup-psi-1} follows for $d=2$, using in addition~\eqref{eq:sup-phi-1} in Lemma~\eqref{lem:diff-phi}.

\noindent For $d>2$, we use Cauchy-Schwarz' inequality with weight
\begin{multline*}
\int_{\R^d} |\delta \ext\phi\delta \ext\psi|dx \,\leq \, \left(\int_{\R^d} \frac{1}{|z-x|^2+1}(\delta \ext\psi(x))^2dx\right)^{\frac{1}{2}} \\
\times \left(\int_{\R^d} (|z-x|^2+1){\ext{h}}_{1}^2(z,x) dx\right)^{\frac{1}{2}}  \left(\int_{B_{3R}(z)}|\nabla \ext{\phi}_1|^2dy+1\right)^{\frac{1}{2}}.
\end{multline*}
On the first factor, we apply Hardy's inequality in the form
$$
\int_{\R^d} \frac{1}{|z-x|^2+1}(\delta \ext\psi(x))^2dx \,\lesssim \, \int_{\R^d} |\nabla \delta \ext\psi|^2dx.
$$
For the second factor, we appeal to \eqref{gu.11} in Definition~\ref{def:Green} for ${\ext{h}}_{1}$ when integrated on dyadic annuli. This yields uniformly w.~r.~t. $z\in \R^d$
$$
\int_{\R^d} (|z-x|^2+1){\ext{h}}_{1}^2(z,x) dx\,\lesssim \,\left\{
\begin{array}{ll}
d=3:& \sqrt{T},\\
d=4:&\ln T,\\
d>4:&1.
\end{array}
\right.
$$
This implies the desired estimate \eqref{eq:sup-psi-1} for $d>2$ by Young's inequality
and \eqref{eq:sup-phi-1} in Lemma~\ref{lem:diff-phi}.

\medskip
\step{2} Proof of \eqref{eq:diff-psi-1} for $|z-x|\geq 2R$.

\noindent The starting point is the Green representation formula \eqref{eq:GRF-psi}.
By Cauchy-Schwarz' inequality, we bound the first term of the RHS by
\begin{multline*}
\left|\int_{\R^d}\nabla_y  G_1(y,x) \cdot (A_1-A_2)(y)\nabla \ext{ \psi}_2(y)dy\right|\\
\,\lesssim \, \left(\int_{B_R(z)}|\nabla_y  G_1(y,x)|^2dy\right)^{\frac{1}{2}} \left(\int_{B_R(z)}|\nabla \ext{\psi}_2|^2dy\right)^{\frac{1}{2}}.
\end{multline*}
We then take the supremum in $A_1$ and $A_2$ using Lemma~\ref{lem:diff-Green} and estimate \eqref{eq:sup-psi-1}, respectively.
This yields
\begin{multline}\label{eq:psi-off-diag1}
\sup_{A_1,A_2}\left|\int_{\R^d}\nabla_y  G_1(y,x) \cdot (A_1-A_2)(y)\nabla \ext{ \psi}_2(y)dy\right|\\
\,\lesssim \, \left(\int_{B_R(z)}|\nabla_y  G_{T}(y,x)|^2dy\right)^{\frac{1}{2}}\left( \int_{B_R(z)}|\nabla \ext\psi_T|^2dy+\nu_d(T)\Big(\int_{B_{3R}(z)}|\nabla \ext\phi_T|^2dy+1\Big)\right)^{\frac{1}{2}}.
\end{multline}
Note that the RHS of \eqref{eq:psi-off-diag1} is only finite for $|z-x|>R$ and that  $\int_{B_R(z)}|\nabla_y  G_{T}(y,x)|^2dy\lesssim 1$
for all for $|z-x|\geq 2R$ by property~\eqref{gu.11} in Definition~\ref{def:Green} and a dyadic decomposition of $B_R(z)$.

\medskip
\noindent For the second term of the RHS of \eqref{eq:GRF-psi}, we bound the Green function pointwise by $g_{T}$, cf. property \eqref{eq:ptwise-decay-estim} in Definition~\ref{def:Green}, and use the oscillation estimate \eqref{eq:diff-phi-1} to bound  $\delta \ext\phi$
\begin{equation}\label{eq:psi-off-diag2}
\sup_{A_1,A_2}\int_{\R^d} G_1(y,x) |\delta \ext\phi(y)| dy 
\,\lesssim \, \int_{\R^d} g_{T}(x-y) {\ext h}_{T}(z,y) dy\left(\int_{B_{3R}(z)}|\nabla \ext\phi_T|^2dy+1\right)^{\frac{1}{2}}.
\end{equation}
Combining the two estimates \eqref{eq:psi-off-diag1} and \eqref{eq:psi-off-diag2} yields  \eqref{eq:diff-psi-1} for $|z-x|\geq 2R$.

\medskip
\step{3} Proof of \eqref{eq:diff-psi-1} for $|z-x|<2R$.

\noindent As in the proof of Lemma~\ref{lem:diff-phi}, it is enough to consider smooth functions $A_2$ of the form
$$
{A}_2|_{\R^d\setminus B_{R}(z)}\,=\,A|_{\R^d\setminus B_{R}(z)},  \quad {A}_2|_{B_{\frac{R}{2}}(z)}\,=\,\Id,
$$
and prove that for all $R>0$ and all $|z-x|\leq \frac{2}{3}R$,
\begin{multline*}
\dps\sup_{A_1,A_2}{\delta\ext\psi(x)}\,\lesssim \, \left(\int_{B_{R}(z)}|\nabla \ext\psi_T|^2dy+ \nu_d(T)\Big(\int_{B_{3R}(z)}|\nabla \ext\phi_T|^2dy +1 \Big)\right)^{\frac{1}{2}}
\\+\left(\int_{B_{R}(z)}|\nabla \ext\phi_T(y)|^2dy+1\right)^{\frac{1}{2}} \int_{\R^d}g_{T}(x-y){\ext h}_{T}(z,y)dy,
\end{multline*}
then to replace $R$ by $3R$ in this estimate, and to use that 
$$\dps\osc{A|_{B_R(z)}}{\ext\psi_T(x)}\,\leq \,\dps\osc{A|_{B_{3R}(z)}}{\ext\psi_T(x)}.$$

\noindent The starting point is again the Green representation formula \eqref{eq:GRF-psi}.
The second term can be dealt with as in Step~2.
For the first term however, due to the singularity of the Green function at $x=y$, we cannot use the Cauchy-Schwarz inequality. Instead, we proceed as 
in the proof of Lemma~\ref{lem:diff-phi}. We split the integrals into two parts, and use H\"older's inequality with exponents $(p,q)$ for some $1<p<\frac{d}{d-1}$  on the first term, and Cauchy-Schwarz' inequality on the second term:
\begin{eqnarray}
\lefteqn{\left|\int_{\R^d}\nabla_y  G_1(y,x) \cdot (A_1-A_2)(y)\nabla \ext{\psi}_2(y)dy\right|}\nonumber \\
&\lesssim &\left( \int_{B_{\frac{R}{8}}(z)} |\nabla_y   G_1(y,x)|^pdy\right)^{\frac{1}{p}} \left(\int_{B_{\frac{R}{8}}(z)} |\nabla \ext{ \psi}_2|^qdy\right)^{\frac{1}{q}}
\label{eq:rhs1}\\
& &+  \left(\int_{B_R(z)\setminus B_{\frac{R}{8}}(z)} |\nabla_y  G_1(y,x)|^2dy\right)^{\frac{1}{2}}  
 \left(\int_{B_R(z)\setminus B_{\frac{R}{8}}(z)} |\nabla  \ext{\psi}_2|^2dy\right)^{\frac{1}{2}} . \label{eq:rhs2}
\end{eqnarray}
We first treat \eqref{eq:rhs1}. The first factor in \eqref{eq:rhs1} is bounded uniformly in $A_1$ since $\nabla G_1 $ is bounded in $L^{p}(\R^d)$ uniformly with respect to $A_1$ and $T>0$, cf. property~\eqref{gu.11} in Definition~\ref{def:Green}. For the second factor, we note that in $B_{\frac{R}{2}}(z)$, $\ext{\psi}_2$ satisfies
\begin{equation*}
T^{-1} \ext{\psi}_2-\triangle \ext{\psi}_2\,=\, \ext{\phi}_2.
\end{equation*}
Hence, for all $i\in \{1,\dots,d\}$,
\begin{equation*}
T^{-1} \partial_{x_i}\ext{\psi}_2-\triangle \partial_{x_i}\ext{\psi}_2\,=\, \partial_{x_i}\ext{\phi}_2,
\end{equation*}
so that by classical interior regularity (see for instance \cite[Theorem~2, Sec.~6.3]{Evans-10}), for all $k\in \N_0$,
\begin{equation}\label{eq:interior-reg-psiT}
\|\nabla \ext{\psi}_2\|_{H^{k+2}(B_{\frac{R}{8}}(z))} \,\lesssim \, \|\nabla \ext{\phi}_2\|_{H^{k}(B_{\frac{R}{4}}(z))} +\|\nabla \ext{ \psi}_2\|_{L^2(B_{\frac{R}{4}}(z))},
\end{equation}
where the multiplicative constant depends on $k$ and $R$.
This yields by Sobolev embedding and the regularity property~\eqref{eq:interior-reg-phiT} in the proof of Lemma~\ref{lem:diff-phi}
\begin{equation*}
\left(\int_{B_{\frac{R}{8}}(z)} |\nabla \ext{ \psi}_2|^qdy\right)^{\frac{1}{q}} \,\lesssim \,\left(\int_{B_{\frac{R}{2}}(z)} |\nabla \ext{ \phi}_2|^2dy\right)^{\frac{1}{2}}+\left(\int_{B_{\frac{R}{4}}(z)} |\nabla \ext{\psi}_2|^2dy\right)^{\frac{1}{2}}.
\end{equation*}
Thus \eqref{eq:rhs1} is bounded as follows
\begin{multline}\label{eq:contrib-1}
\left( \int_{B_{\frac{R}{8}}(z)} |\nabla_y   G_1(y,x)|^pdy\right)^{\frac{1}{p}} \left(\int_{B_{\frac{R}{8}}(z)} |\nabla \ext{\psi}_2|^qdy\right)^{\frac{1}{q}}
\\
\,\lesssim\, \left(\int_{B_{R}(z)} |\nabla  \ext{\phi}_2|^2dy\right)^{\frac{1}{2}}+\left(\int_{B_{R}(z)} |\nabla \ext{\psi}_2|^2dy\right)^{\frac{1}{2}}.
\end{multline}
We now turn to \eqref{eq:rhs2}. We recall that the first factor in \eqref{eq:rhs2} is bounded by 1, cf. property~\eqref{gu.11} in Definition~\ref{def:Green}, so that
\begin{equation}\label{eq:contrib-2}
\left(\int_{B_R(z)\setminus B_{\frac{R}{8}}(z)} |\nabla_y  G_{1}(y,x)|^2dy\right)^{\frac{1}{2}}  \left(\int_{B_R(z)\setminus B_{\frac{R}{8}}(z)} |\nabla  \ext{\psi}_2|^2dy\right)^{\frac{1}{2}}\,\lesssim \, \left(\int_{ B_{R}(z)} |\nabla  \ext{\psi}_2|^2dy\right)^{\frac{1}{2}} .
\end{equation}
Appealing to \eqref{eq:sup-phi-1} and \eqref{eq:sup-psi-1} to estimate the supremum with respect to  $A_2$ of the RHS of \eqref{eq:contrib-1} and \eqref{eq:contrib-2} completes the oscillation estimate for $|z-x|< 2R$.


\appendix

\section{Proofs of the other auxiliary lemmas}

\subsection{Proof of Lemma~\ref{lem:def-phiT-conv}}

We proceed in two steps, first sketch the argument for the existence, and then turn to uniqueness.
W.~l.~o.~g. we may consider $T=1$ by scaling.

\medskip

\step{1} Existence.

\noindent Let $\xi\in \R^d$. To obtain a sequence of approximate solutions $\ext\phi_{R}$, we solve \eqref{eq:app-corr-distribution} on balls $B_R$ with increasing radii and homogeneous Dirichlet boundary conditions.
We test the defining equation for $\ext\phi_{R}$ 
with the function $\eta_z^2 \ext\phi_{R}$ where $\eta_z(x)=\exp(-c |z-x|)$ for arbitrary $z\in \R^d$ and some $c>0$ to be fixed later.
This yields
\begin{multline*}
 \int_{B_R}\eta_z^2 \ext\phi_{R}^2 dx+ \int_{\R^d} \eta_z^2\nabla \ext\phi_{R}\cdot A \nabla \ext \phi_{R} dx\\
\,=\, -2\int_{B_R} \ext\phi_R\eta_z\nabla \eta_z\cdot A\nabla \ext\phi_Rdx-\int_{B_R}\eta_z^2 \nabla \ext\phi_R\cdot A\xi dx -2\int_{B_R}\eta_z\ext\phi_R \nabla \eta_z\cdot A\xi 
dx,
\end{multline*}
which, by the bounds on $A$ and Young's inequality on each term of the RHS with constants $\kappa,2\kappa$ and $\kappa>0$, respectively, turns into
\begin{equation}\label{eq:apriori-kappa}
\int_{B_R}(\eta_z^2 -\frac{2}{\kappa}|\nabla \eta_z|^2)\ext\phi_{R}^2 dx+
\lambda \int_{B_R} \eta_z^2(1-2\frac{\kappa}{\lambda}) |\nabla \ext\phi_{R}|^2 dx
\,\leq \, (\frac{1}{4\kappa}+\kappa)|\xi|^2 \int_{\R^d}\eta_z^2dx.
\end{equation}
Choosing $\kappa=\frac{\lambda}{4}$ and $c=\frac{\sqrt{\lambda}}{4}$
then yields the a priori estimate
\begin{equation*}
 \int_{B_R}\eta_z^2 \ext\phi_{R}^2 dx+ \int_{\R^d} \eta_z^2 |\nabla  \ext\phi_{R}|^2 dx\\
\,\lesssim \, |\xi|^2\int_{\R^d}\eta_z^2dx.
\end{equation*}
By weak compactness, the sequence $\ext\phi_{R}$ weakly converges in $H^1_\loc(\R^d)$ up to extraction to some function $\ext\phi$, 
which is a distributional solution of \eqref{eq:app-corr-distribution}  on $\R^d$. In addition, $\ext\phi$ satisfies the a priori estimate 
\begin{equation}\label{eq:unif-averaged-bound1}
\int_{\R^d}\eta_z^2\ext \phi^2dx + \int_{\R^d} \eta_z^2 |\nabla  \ext\phi|^2 dx\\
\,\lesssim \, |\xi|^2\int_{\R^d}\eta_z^2dx,
\end{equation}
which implies 
\eqref{eq:uniqueness-cond-refined} since its RHS does not depend on $z$.

\medskip

\step{2} Uniqueness.

\noindent Let $\delta \ext\phi$ be such that $\limsup_{t\uparrow \infty} \fint_{B_t} \Big((\delta \ext\phi)^2+|\nabla \delta \ext\phi|^2\Big)dx <\infty$ 
and satisfy \eqref{eq:app-corr-distribution} with $\xi=0$.
Let $\eta_{0}$ be as in Step~1 for $z=0$. We first argue that 
$$
 \int_{\R^d}\eta_0^2 \delta \ext\phi^2 dx+ \int_{\R^d} \eta_0^2|\nabla \delta \ext\phi|^2 dx\,<\,\infty.
 $$
Indeed, by assumption, there exists $C<\infty$ such that $\sup_{t\ge 1}\fint_{B_t} \Big((\delta \ext\phi)^2+|\nabla \delta \ext\phi|^2\Big)dx \leq C$, so that for all $N\in \N$,
\begin{multline*}
 \int_{B_N}\eta_0^2 \delta \ext\phi^2 dx+ \int_{B_N} \eta_0^2|\nabla \delta \ext\phi|^2 dx\,\lesssim\,\sum_{t=1}^N t^d \exp(2c(1-t))\fint_{B_t} ( \delta \ext\phi^2 +|\nabla \delta \ext\phi|^2)dx \\
 \,\leq \, C\sum_{t=1}^\infty t^d \exp(-2ct) <\infty.
 \end{multline*}
We may thus test equation \eqref{eq:app-corr-distribution}  with test function $\eta_{0,R}^2 \delta \ext\phi$, where $\eta_{0,R}=\eta_0\mu_R$ and $\mu_R$ a smooth cut-off function 
on $B_R$. Passing to the limit $R\uparrow \infty$ by dominated convergence leads to the energy estimate
\begin{equation*}
 \int_{\R^d}\eta_0^2 \delta \ext\phi^2 dx+ \int_{\R^d} \eta_0^2\nabla \delta \ext\phi\cdot A \nabla \delta\ext \phi dx
\,=\, -2\int_{\R^d} \eta_0 \delta\ext \phi\nabla \eta_0 \cdot A\nabla \delta\ext\phi dx.
\end{equation*}
Using Young's inequality as for \eqref{eq:apriori-kappa} then yields
\begin{equation*}
\int_{\R^d}(\eta_{0}^2-\frac{2}{\kappa}|\nabla \eta_{0}|^2) (\delta \ext\phi)^2 dx+ \int_{\R^d} \eta_{0}^2(1-2\frac{\kappa}{\lambda})|\nabla \delta\ext\phi|^2 dx
\,\leq\,  0,
\end{equation*}
which, with the choice $\kappa=\frac{\lambda}{4}$ and $c=\frac{\sqrt{\lambda}}{4}$ as in Step~1, establishes uniqueness.

\subsection{Proof of Lemma~\ref{lem:diff-Green}}

By a standard regularization argument, one may assume that $\tilde A$ and $A$ are smooth and coincide outside $B_R(z)$, $z\in \R^d$.
We denote by $G_T,\tilde G_T\in W^{1,1}(\R^d)$ the associated Green functions, $T>0$.
Substracting the defining equations \eqref{eq:Green} with singularity at $y\in \R^d$ for $G_T$ and $\tilde G_T$ then yields
\begin{multline}\label{eq:pr-diff-green-1}
T^{-1}(\tilde G_T(x,y)-G(x,y))-\nabla_x\cdot (\tilde A(x)\nabla_x(\tilde G_T(x,y)-G_T(x,y))) \,\\
=\,\nabla_x\cdot ((\tilde A-A)(x)\nabla_xG_T(x,y))
\end{multline}
in the sense of distributions on $\R^d_x$.
Since $G_T$ and $\tilde G_T$ belong to $C^\infty(\R^d\times \R^d\setminus\{x=y\})$,
the RHS of \eqref{eq:pr-diff-green-1} is smooth with support in $B_R(z)$ provided $|z-y|>R$.
Hence $G_T(\cdot,y)-\tilde G_T(\cdot,y)$ is also a classical solution of \eqref{eq:pr-diff-green-1}
and therefore belongs to $H^1(\R^d)$ since the RHS has compact support.

\medskip
\noindent  The energy estimate  yields
\begin{multline*}
\int_{\R^d} |\nabla_x(\tilde G_T(x,y)-G_T(x,y))|^2dx \\\,\lesssim \, \int_{\R^d}\nabla_x (\tilde G_T(x,y)-G(x,y)) \cdot (\tilde A-A)(x)\nabla_xG_T(x,y)dx.
\end{multline*}
Using that $A$ and $\tilde A$ coincide outside $B_R(z)$ together with the Cauchy-Schwarz inequality, this turns into 
\begin{equation*}
\left(\int_{B_R(z)} |\nabla_x(\tilde G_T(x,y)-G_T(x,y))|^2dx\right)^{\frac{1}{2}}\,\lesssim \, \left(\int_{B_R(z)}|\nabla_xG_T(x,y)|^2dx\right)^{\frac{1}{2}},
\end{equation*}
so that by the triangle inequality
\begin{equation*}
\left(\int_{B_R(z)} |\nabla_x\tilde G_T(x,y)|^2dx\right)^{\frac{1}{2}}\,\lesssim \, \left(\int_{B_R(z)}|\nabla_xG_T(x,y)|^2dx\right)^{\frac{1}{2}},
\end{equation*}
as desired.


\subsection{Properties of the Green functions}

We first address the existence part. By a scaling argument, it is sufficient to
consider $T=1$ (we thus drop the subscript $T=1$ from our notation). 
By a standard approximation argument, it is sufficient to consider
the case of a {\it smooth} uniformly elliptic coefficient field $A$ on a (large) ball $D$.
Let $G(x,y)$ denote the Green function for these data, which is known to exist by
classical theory. For the above properties of the whole-space, non-smooth coefficients
Green function, it is enough to establish the following properties of $G$:
\begin{itemize}
\item Uniform, but qualitative continuity off the diagonal and off the boundary, 
that is, for all $r>0$:
\begin{equation}\label{gu.21}
\begin{array}{c}
\{(x,y)\in D^2\,|\,{\mathrm{dist}}(\{x,y\},\partial D)\ge 2,\,|x-y|>r\}\ni (x,y)\mapsto G(x,y)\\[2ex]
\quad\mbox{has modulus of continuity only depending on}\;d, \lambda, r,
\end{array}
\end{equation}
but not on the modulus of continuity of $A$ nor on $D$. By Arzel\`a-Ascoli's compactness criterion, it is this 
equi-continuity that ensures the continuity \eqref{gu.14} when taking the limit in the approximation argument.
\item Pointwise upper bounds on $G$: For $x,y$ that stay away from
the boundary in the sense of ${\mathrm{dist}}(\{x,y\},\partial D)\ge 1$ we claim
\begin{equation}\label{gu.20}
G(x,y)\lesssim\exp(-c|x-y|)
\left\{\begin{array}{ccc}
\ln(2+\frac{1}{|x-y|})&\mbox{for}&d=2\\
|x-y|^{2-d}&\mbox{for}&d>2
\end{array}\right\}.
\end{equation}
It is obvious that under the locally uniform convergence off the diagonal coming from Arzel\`a-Ascoli's compactness
criterion this turns into \eqref{eq:ptwise-decay-estim} in the limit.
\item Averaged bounds on $\nabla_x G$ and $\nabla_y G$: For ${\mathrm{dist}}(y,\partial D)\ge 1$ we have
\begin{eqnarray}
\left(R^{-d}\int_{D\cap\{R<|x-y|<2R\}}|\nabla_xG(x,y)|^2dx\right)^\frac{1}{2}&\lesssim& \exp(-cR) R^{1-d},\label{gu.18}
\end{eqnarray}
and for ${\mathrm{dist}}(x,\partial D)\ge 1$ we have
\begin{eqnarray}
\left(R^{-d}\int_{D\cap\{R<|y-x|<2R\}}|\nabla_yG(x,y)|^2dy\right)^\frac{1}{2}&\lesssim& \exp(-cR) R^{1-d}.\label{gu.19}
\end{eqnarray}
By lower semi-continuity of these expressions under pointwise convergence of $G$,
\eqref{gu.18} \& \eqref{gu.19} turn into \eqref{gu.11} \& \eqref{gu.12} in the limit.
\item Differential equation:
\begin{eqnarray}
G-\nabla_x\cdot A(x)\nabla_xG=\delta(x-y)&&\mbox{distributionally in}\;D_x,\nonumber\\
G-\nabla_y\cdot A^*(y)\nabla_yG=\delta(y-x)&&\mbox{distributionally in}\;D_y.\nonumber
\end{eqnarray}
Since \eqref{gu.20} and \eqref{gu.18} \& \eqref{gu.19} imply local equi-integrability for
$\mathbb{R}^d\ni x\mapsto (G(x,y),\nabla_x G(x,y))$ and $\mathbb{R}^d\ni y\mapsto (G(x,y),\nabla_y G(x,y))$,
this yields \eqref{eq:Green} \& \eqref{eq:Green-t} in the limit.
\end{itemize}

\medskip

\noindent We now come to a further reduction: Because of the symmetry of our assumptions under exchanging the
roles of $x$ and $y$, we may restrict to the $x$-variable in proving the above estimates.
Because our assumptions are invariant under translation, we may restrict to the case of $y=0$ 
and we may assume
\begin{equation}\label{gu.30}
{\mathrm{dist}}(0,\partial D)\ge 1.
\end{equation}
We thus suppress the $y$-dependence in our notation and just write $G(x)$,
which is characterized as the solution of
\begin{equation}\label{g.1}
G-\nabla\cdot A\nabla G=\delta\quad\mbox{in}\;D\quad\mbox{and}\quad G=0\quad\mbox{on}\;\partial D.
\end{equation}
It will be convenient to separate the near-field behavior dominated by the singularity (i.\ e.\ for $|x|\ll 1$)
from the far-field behavior dominated by the massive term (i.\ e.\ for $|x|\gg 1$):
\begin{itemize}
\item Pointwise upper bounds on $G$: We shall show 
\begin{equation}\label{g.3near}
G(x)\lesssim
\left\{\begin{array}{ccc}
\ln(2+\frac{1}{|x|})&\mbox{for}&d=2\\
|x|^{2-d}&\mbox{for}&d>2
\end{array}\right\}\quad\mbox{for}\;|x|<\frac{2}{3}
\end{equation}
and
\begin{equation}\label{g.3far}
G(x)\lesssim \exp(-c|x|)\quad\mbox{for}\quad|x|\ge\frac{2}{3},\;{\mathrm dist}(x,\partial D)\ge 1.
\end{equation}
This yields \eqref{gu.20} (with a reduced value for the generic $c>0$).
\item Uniform, but qualitative continuity of $G$: We note that by De Giorgi's a priori estimate of 
the H\"older modulus of an $A$-harmonic function,
these quantitative pointwise estimates yield that for all $r>0$
\begin{equation}\nonumber
\begin{array}{c}
\{x\in D\,|\,{\mathrm{dist}}(x,\partial D)\ge 2,\,|x|>r\}\ni x\mapsto G(x)\\[2ex]
\quad\mbox{has modulus of continuity only depending on}\;d, \lambda, r.
\end{array}
\end{equation}
Note that because of the massive term (which however is under good control because of \eqref{g.3near} \& \eqref{g.3far}), 
we need a version of De Giorgi's estimate with a (bounded) right hand side, see for instance \cite[Theorem~4.1]{Han-Lin-97}.
Since a uniform modulus of continuity of $G(x,y)$ in $x$ (for all $y$) and a uniform modulus of continuity in $y$
(for all $x$) implies a uniform modulus of continuity in $(x,y)$, this yields \eqref{gu.21}.
\item Average estimates on $\nabla G$: We shall show
\begin{eqnarray}
\left(R^{-d}\int_{D\cap\{R<|x|<2R\}}|\nabla G|^2dx\right)^\frac{1}{2}&\lesssim& R^{1-d}\quad\mbox{for}\;0<R\le\frac{1}{6},
\label{gu.2}\\
\left(R^{-d}\int_{D\cap\{|x|>R\}}|\nabla G|^2dx\right)^\frac{1}{2}&\lesssim&\exp(-cR)\quad\mbox{for}\;R\ge\frac{1}{6}.
\label{gu.3}
\end{eqnarray}
This implies \eqref{gu.18}.
\end{itemize}

\medskip


\noindent Our argument for \eqref{g.3near}, \eqref{g.3far}, \eqref{gu.2}, and \eqref{gu.3}, is self-contained with the exception of 
De Giorgi's a priori estimate for $A$-subharmonic functions
$u$ (i.\ e.\ satisfying $-\nabla\cdot A\nabla u\le 0$) in some ball $B_R(x)$
\begin{equation}\label{g.2}
u(x)\lesssim R^{-d}\int_{B_{R}(x)}\max\{u,0\}.
\end{equation}
With this key ingredient, we split the proof into several easy steps.

\medskip

\step{1} Near-field estimates on $G$.

\noindent We start by establishing {\it average} near-field estimates on $G$.
We start with the easier case of $d>2$ and shall establish
\begin{equation}\label{g.26}
R^{-d}\int_{D\cap\{|x|\le R\}}Gdx\lesssim R^{2-d}\quad\mbox{for all}\;R>0,
\end{equation}
reproducing the classical argument of Gr\"uter \& Widman \cite[(1.1) Theorem]{Grueter-Widman-82}.
To this purpose, we test \eqref{g.1} with
$\min\{G,M\}$ for an arbitrary $0\le M<\infty$. Using the uniform ellipticity $A\ge\lambda{\mathrm Id}$
we obtain the inequality
\begin{equation}\nonumber
\int_{D}\min\{G,M\}^2dx+\lambda\int_{D}|\nabla\min\{G,M\}|^2dx\le M.
\end{equation}
We throw away the first positive term, which comes from the massive term.
With help of the scale invariant Sobolev's estimate (here we use $d>2$) on $D$ 
(with vanishing boundary data) this yields
\begin{equation}\nonumber
\int_{D}|\min\{G,M\}|^\frac{2d}{d-2}dx\lesssim M^\frac{d}{d-2},
\end{equation}
from which, redefining $\frac{M}{2}$ to be $M$, we deduce the
weak $L^\frac{d}{d-2}$-estimate
\begin{equation}\nonumber
|D\cap \{G>M\}|\lesssim M^{-\frac{d}{d-2}},
\end{equation}
where $|\cdot|$ denotes the $d$-dimensional volume.
We now restrict to the ball of radius $R$:
\begin{equation}\nonumber
|D\cap\{|x|<R\}\cap\{G>M\}|\lesssim \min\{R^{-d},M^{-\frac{d}{d-2}}\}
\end{equation}
and integrate over $M\in(0,\infty)$ to recover the $L^1$-norm:
\begin{eqnarray*}\nonumber
\int_{D\cap\{|x|\le R\}}Gdx&=&\int_0^\infty|D\cap\{|x|<R\}\cap\{G>M\}|dM\\
&\lesssim& \int_0^\infty\min\{R^{-d},M^{-\frac{d}{d-2}}\}dM\\
&\stackrel{M=R^{-(d-2)}\hat M}{=}&
R^2\int_0^\infty\min\{1,\hat M^{-\frac{d}{d-2}}\}d\hat M\sim R^2,
\end{eqnarray*}
which establishes \eqref{g.26}.

\medskip

\noindent The average near-field estimates on $G$ is more subtle for $d=2$;
in fact, one naturally controls only the {\it oscillation} of $G$ in the sense of
\begin{equation}\label{g.20}
\left(R^{-2}\inf_{c\in\mathbb{R}}\int_{|x|<R}(G-c)^2dx\right)^\frac{1}{2}\lesssim 1
\quad\mbox{for all}\;0<R\le 1.
\end{equation}
We note that because of \eqref{gu.30} and $R\le 1$, we have $\{|x|<R\}\subset D$.
The argument for \eqref{g.20} mimics \cite[Lemma~10]{Lamacz-Neukamm-Otto-13}
which is a simplification of \cite[Lemma~2.8]{Gloria-Otto-09}, 
which itself was a quantification of \cite[Lemma~2.5]{Dolzmann-Muller-95}.
Let $c_R$ denote the {\it median} of $G$ over $\{|x|\le R\}$. Following the
argument for $d>2$, we test \eqref{g.1} with the truncated
$G-c_R$, that is $\max\{\min\{G-c_R,M\},-M\}$ for some arbitrary $0\le M<\infty$. 
Since the test function has no sign, the massive term now gets into our way. 
However, since $G\ge0$ in $D$ by the maximum principle, we have for the normal derivative 
$\nu\cdot A\nabla G\le0$ on $\partial D$ so that integrating \eqref{g.1} yields
$\int_{D}Gdx\le 1$. Hence we may rewrite \eqref{g.1} as $-\nabla\cdot A\nabla(G-c_R)=f:=\delta-G$
with the total variation of the signed measure $f$ bounded by $1+\int Gdx\le 2$. Therefore, testing yields
\begin{equation}\nonumber
\lambda\int_{D}|\nabla\max\{\min\{G-c_R,M\},-M\}|^2dx\le 2M,
\end{equation}
which we reduce to the ball $\{|x|\le R\}$ and split into
\begin{equation}\label{g.17}
\int_{|x|\le R}|\nabla\min\{\max\{\pm(G-c_R),0\},M\}|^2dx\lesssim M.
\end{equation}
By symmetry, it is enough to show that the plus sign in \eqref{g.17} implies
\begin{equation}\label{g.18}
R^{-2}\int_{|x|\le R}u^2dx\lesssim 1,\quad\mbox{where}\quad u:=\max\{G-c_R,0\}.
\end{equation}
Here comes the argument:
By definition of the median $c_R$, $u$ and thus a fortiori $\min\{u,M\}$ vanishes on 
at least half of the ball $\{|x|\le R\}$.
Hence by a Poincar\'e-Sobolev estimate on $\{|x|\le R\}$ we obtain that
\begin{equation}\nonumber
\left(R^{-2}\int_{|x|\le R}\min\{u,M\}^6dx\right)^\frac{1}{6}\lesssim
\left(\int_{|x|\le R}|\nabla\min\{u,M\}|^2dx\right)^\frac{1}{2}
\stackrel{(\ref{g.17})}{\lesssim}M^\frac{1}{2},
\end{equation}
where there is nothing specific to the exponent 6, in fact, any finite exponent larger than 4 
would do.
As in the previous step, this yields the weak-type estimate
\begin{equation}\nonumber
\left(R^{-2}|\{|x|\le R\}\cap\{u>M\}|\right)^\frac{1}{6}\lesssim
\min\{1,M^{-\frac{1}{2}}\},
\end{equation}
which (after taking the sixth power) we integrate against $\int_0^\infty\cdot MdM$ to obtain
the (squared) $L^2$-norm
\begin{equation}\nonumber
R^{-2}\int_{|x|<R}u^2dx\lesssim \int_0^\infty\min\{1,M^{-3}\} MdM\sim 1.
\end{equation}
This establishes \eqref{g.18} and thus \eqref{g.20}.

\medskip

\noindent In order to ``anchor'' the $(d=2)$-estimate \eqref{g.20} on $G$,
we need the following average intermediate-scale estimate on $G$
\begin{equation}\label{g.21}
\left(\int_{|x|\le 1}G^2dx\right)^\frac{1}{2}\lesssim 1.
\end{equation}
As opposed to the previous step, we now use the massive term to our advantage
by testing with $\min\{G,M\}$ as in case of $d>2$:
\begin{equation}\nonumber
\int_{D}\min\{G,M\}^2dx+\lambda\int_{D}|\nabla\min\{G,M\}|^2dx\le M.
\end{equation}
Note that since $\{|x|\le 1\}\subset D$, cf.\ \eqref{gu.30}, we may restrict the estimate to the ball 
$\{|x|\le 1\}$ where we use a Sobolev estimate to obtain
\begin{equation}\nonumber
\left(\int_{|x|\le 1}\min\{G,M\}^6dx\right)^\frac{1}{6}\lesssim M^\frac{1}{2}.
\end{equation}
We then proceed as in the previous step (with $R=1$).

\medskip

\noindent Equipped with \eqref{g.20} and \eqref{g.21}, we now may
complete the average near-field estimate on $G$ in case of $d=2$:
\begin{equation}\label{g.13}
\left(R^{-2}\int_{|x|\le R}G^2dx\right)^\frac{1}{2}
\lesssim \ln(2+\frac{1}{R})\quad\mbox{for all}\;0<R\le 1.
\end{equation}
An elegant way to obtain such a logarithmic estimate, even directly in its pointwise version,
is a dimension reduction from $d=3$ as in Avellaneda \& Lin \cite{Avellaneda-Lin-87};
however, we need the BMO-type bound \eqref{g.20} also for the average near-field estimate
on $\nabla G$ so that we opt for a derivation of \eqref{g.13} from \eqref{g.20}.
We consider dyadic radii $R=2^{-n}$ with $n\in\mathbb{N}_0$.
Let $c_n$ denote the average of $G$ over $\{|x|<2^{-n}\}$. From \eqref{g.20} 
for $R=2^{-n}$ we learn in particular that $|c_{n+1}-c_n|\lesssim 1$, whereas from \eqref{g.21} 
we get in particular $|c_0|\lesssim 1$. Hence we obtain $|c_{n}|\lesssim n+1$
and thus once again from \eqref{g.20}
\begin{equation}\nonumber
\left(R^{-2}\int_{|x|\le 2^{-n}}G^2dx\right)^\frac{1}{2}\lesssim n+1,
\end{equation}
which translates into \eqref{g.13}.

\medskip

\noindent We now obtain the desired {\it pointwise} near-field estimates \eqref{g.3near} on $G$ as follows:
Since $G\ge 0$, $G$ is a subsolution of $-\nabla\cdot A\nabla$ away from the origin and thus \eqref{g.3near} 
follows from (\ref{g.26}) (for $d>2$) and (\ref{g.13}) (for $d=2$) by
applying De Giorgi's result \eqref{g.2} to balls $B$ with center $x$ and radius $R=\frac{|x|}{2}$
(which by \eqref{gu.30} and $|x|\le\frac{2}{3}$ is contained in $D$).

\medskip

\step{2} Far-field estimates on $G$.

\noindent We start with the average version of the far-field estimates
--- all dimensions can be treated simultaneously:
\begin{equation}\label{g.27}
\int_{D\cap\{|x|\ge\frac{1}{3}\}}(\exp(c|x|)G)^2dx\lesssim 1.
\end{equation}
For this purpose, we fix
a smooth cut-off function $\eta$ that vanishes in $\{|x|\le\frac{1}{6}\}$ but is equal to one
on $\{|x|\ge\frac{1}{3}\}$ and will show that
\begin{equation}\label{g.5}
\int_{D}\eta^2\exp(2c|x|)G^2dx\lesssim 1.
\end{equation}
In order to establish
\eqref{g.5}, we follow Caccioppoli's strategy as modified by Agmon \cite{Agmon-82}
and test \eqref{g.1} with $\eta^2\exp(2c|x|)G$ to the effect of
\begin{equation}\nonumber
\int_{D}\eta^2\exp(2c|x|)G^2dx+\int_{D}\nabla(\eta^2\exp(2c|x|)G)\cdot A\nabla Gdx=0.
\end{equation}
Introducing the abbreviation $\tilde\eta:=\eta\exp(c|x|)$ we now use the 
pointwise inequality
\begin{eqnarray*}
\nabla(\tilde\eta^2G)\cdot A\nabla G&=&\tilde\eta^2\nabla G\cdot A\nabla G
-2G\tilde\eta \nabla\tilde\eta\cdot A\nabla G\\
&\ge&\lambda{\tilde\eta}^2|\nabla G|^2-2|G||\tilde\eta||\nabla\tilde\eta||\nabla G|\\
&\ge&-\frac{1}{\lambda}G^2|\nabla\tilde\eta|^2
\end{eqnarray*}
to obtain the integral inequality
\begin{eqnarray*}\nonumber
\lefteqn{\int_{D}\eta^2\exp(2c|x|)G^2dx\le\frac{1}{\lambda}\int G^2|\nabla(\eta\exp(c|x|))|^2dx}\nonumber\\
&\le&
\frac{2c}{\lambda}\int_{D}G^2\eta^2\exp(2c|x|)dx+\frac{2}{\lambda}\int_{D} G^2|\nabla\eta|^2dx.
\end{eqnarray*}
The second RHS term, which by choice of $\eta$ is supported in $\{\frac{1}{6}<|x|<\frac{1}{3}\}$, 
is $\lesssim 1$ by the pointwise
near-field estimates \eqref{g.3near}, the first r.\ h.\ s.\ term can
be absorbed into the LHS provided $c<\frac{1}{2\lambda}$. This establishes \eqref{g.5} and thus \eqref{g.27}.

\medskip

\noindent We now obtain the {\it pointwise} far-field estimates \eqref{g.3far} on $G$
from \eqref{g.27} via De Giorgi's result \eqref{g.2} applied to a ball $B$ with
center $x$ and radius $R=\frac{1}{3}$.

\medskip

\step{3} Average estimates on the {\it gradient}  $\nabla G$.

\noindent The near-field estimates \eqref{gu.2} are easy for $d>2$:
This follows from \eqref{g.3near} via the standard Caccioppoli estimate based on testing \eqref{g.1} with
$\eta^2G$, where $\eta$ is a cut-off function for the annulus $\{R<|x|<2R\}$ in the
annulus $\{\frac{R}{2}<|x|<4R\}$. The massive term produces a good term that we discard.

\medskip

\noindent In case of $d=2$, \eqref{gu.2}
follows from the average near-field estimate \eqref{g.20} on the
oscillation of $G$ via a standard Caccioppoli estimate based on testing \eqref{g.1} with
$\eta^2(G-c)$, where $\eta$ is a cut-off function for the annulus $\{R<|x|<2R\}$ in the
annulus $\{\frac{R}{2}<|x|<4R\}$, and $c$ is the average of $G$ over $\{|x|\le 4R\}$. 
As opposed to the previous step, the massive term gets into our way by generating
the following RHS term, which however is lower order (in $R\ll 1$):
\begin{eqnarray*}
R^{-2}\int_D\eta^2(G-c) Gdx&\lesssim&
\left(R^{-2}\int_{|x|<4R}(G-c)^2dx R^{-2}\int_{|x|<4R}G^2dx\right)^\frac{1}{2}\\
&\stackrel{(\ref{g.20}),(\ref{g.3near})}{\lesssim}&\left(\ln(2+\frac{1}{R})\right)^\frac{1}{2}.
\end{eqnarray*}

\medskip

\noindent The far-field estimates \eqref{gu.3} can again be easily treated for all $d$:
They follows from the average far-field estimates \eqref{g.27} on $G$ 
(employed for $|x|\sim R$, w.\ l.\ o.\ g.\ $R\gg 1$) via a standard Caccioppoli estimate based on testing \eqref{g.1} with
$\eta^2G$, where $\eta$ is a cut-off function for $\{|x|>R\}$ in
$\{|x|>\frac{R}{2}\}$. The massive term produces a good term that we discard.

\medskip


\step{4} Uniqueness argument.

\noindent The uniqueness argument is different from \cite{Grueter-Widman-82} (who do 
not consider the whole-space case with a massive term) in the sense it makes stronger
assumptions, namely \eqref{gu.11}, but uses less machinery, namely no lower pointwise bounds
coming from Harnack's inequality.
By scaling, we may still assume that $T=1$. We fix a uniformly elliptic (but not necessarily smooth)
coefficient field $A$. We consider a Green function $G(x,y)$.

\medskip

\noindent The main technical step of our uniqueness argument is the following:
For any $\e>0$, we consider the mollification of $G(x,y)$ in $y$, say,
\begin{equation}\nonumber
G_\e(x,y)=\e^{-d}\int_{|y'-y|<\e}G(x,y')dy'.
\end{equation}
We claim that
\begin{equation}\label{gu.6}
\int_{\mathbb{R}^d}G_\e^2(x,y)+|\nabla_xG_\e(x,y)|^2dx<\infty
\quad\mbox{for all}\;\e>0,\;y\in\mathbb{R}^d.
\end{equation}
Here comes the argument for \eqref{gu.6}: We note that a dyadic decomposition shows that \eqref{gu.1} and \eqref{gu.2}
(together with \eqref{eq:ptwise-decay-estim}) implies
that for any fixed $\alpha>d-2$, say $\alpha=d-1$, we have
\begin{equation}\label{gu.1}
\int_{\mathbb{R}^d}|x-y|^\alpha(G^2(x,y)+|\nabla_x G(x,y)|^2)dx\lesssim 1\quad\mbox{for all}\;y\in\mathbb{R}^d.
\end{equation}
Since $\alpha<d$, we obtain because of 
$\nabla_xG_\e(x,y)=\e^{-d}\int_{|y'-y|\le\e}\nabla_xG(x,y')dy'$
by Cauchy-Schwarz in $y'$
\begin{eqnarray}\nonumber
\lefteqn{G_\e^2(x,y)+|\nabla_xG_\e(x,y)|^2}\nonumber\\
&\le&\e^{-d}\int_{|y'-y|\le\e}|x-y'|^{-\alpha}dy'
\e^{-d}\int_{|y'-y|\le\e}|x-y'|^\alpha(G^2(x,y')+|\nabla_xG(x,y')|^2)dy'\nonumber\\
&\stackrel{\alpha<d}{\lesssim}&
\e^{-d-\alpha}\int_{|y'-y|\le\e}|x-y'|^\alpha(G^2(x,y')+|\nabla_xG(x,y')|^2)dy'\nonumber
\end{eqnarray}
and thus by \eqref{gu.1}
\begin{eqnarray*}
\lefteqn{\int_{\mathbb{R}^d}G_\e^2(x,y)+|\nabla_x G_\e(x,y)|^2dx}\nonumber\\
&\lesssim&
\e^{-d-\alpha}\int_{|y'-y|\le\e}\int_{\mathbb{R}^d}|x-y'|^\alpha(G^2(x,y')+|\nabla_xG(x,y')|^2)dxdy'
\stackrel{(\ref{gu.1})}{\lesssim}\e^{-\alpha},
\end{eqnarray*}
which is a quantification of \eqref{gu.6}.

\medskip

\noindent We now come to the uniqueness argument proper and consider the difference $u(x,y)$ of
two Green's functions. By assumption, we know that for fixed $y$, $u(\cdot,y)$ and $\nabla_xu(\cdot,y)$
are integrable and satisfy
\begin{equation}\nonumber
u-\nabla_x\cdot A(x)\nabla_x u=0\quad\mbox{distributionally in}\;\mathbb{R}^d.
\end{equation}
This persists for the mollification $u_\e(\cdot,y)$ in the $y$-variable
introduced in the previous step:
\begin{equation}\label{gu.4}
u_\e-\nabla_x\cdot A(x)\nabla_x u_\e=0\quad\mbox{distributionally in}\;\mathbb{R}^d.
\end{equation}
On the other hand, we know from \eqref{gu.6}
that the $u_\e(\cdot,y)$ and $\nabla_x u_\e(\cdot,y)$ are square integrable.
This means that we may test \eqref{gu.4} with $u_\e$ to the effect of
\begin{equation}\nonumber
\int_{\mathbb{R}^d}u_\e^2(x,y)+|\nabla_xu_\e(x,y)|^2dx=0.
\end{equation}
This implies $u_\e(x,y)=0$ for almost every $x$ and all $y$. By the continuity
property \eqref{gu.14}, this yields at first $u_\e(x,y)=0$ for all $x\not=y$
and then in the limit $\e\downarrow 0$ that $u(x,y)=0$ for $x\not=y$,
thus establishing uniqueness.


\subsection{Proof of Lemma~\ref{lem:annealed-estim}}

We follow \cite{Marahrens-Otto-14} and split the proof into four steps.
Let $1\le p\le \bar{p}$ where $\bar p$ is as in Lemma~\ref{lem:int-grad}.

\medskip

\step{1} We claim that by Lemma~\ref{lem:int-grad} we have for any radius $R$
\begin{equation}\label{S.37}
\expec{R^{-d}\int_{R<|y|\le2R} \Big( |\nabla\nabla G_T(0,y)|^{2p}+R^{-2p}|\nabla_xG_T(0,y)|^{2p} \Big) \;dy }^\frac{1}{2p}
\lesssim R^{-d}\exp(-c\frac{R}{\sqrt{T}}).
\end{equation}
Indeed, by stationarity we have 
$$\expec{|\nabla \nabla G_T(x,y)|^{2p}}=\expec{|\nabla \nabla G_T(0,y-x)|^{2p}}\text{ and } \expec{|\nabla_x G_T(0,y)|^{2p}}=\expec{|\nabla_y G_T(-y,0)|^{2p}},$$
 so that
\eqref{S.37} follows by taking the expectation of the $(2p)^{th}$ power of \eqref{eq:int-grad2} and \eqref{eq:int-grad}.

\medskip

\step{2}  Consider the $A$-dependent functions $u=u(x;A)$, $f(x;A)$, $h(x;A)$, and the vector
field $g=g(x;A)$ related by
\begin{equation}\label{S.30}
T^{-1}u-\nabla \cdot A\nabla u=\nabla\cdot g+f+T^{-1}h\quad\mbox{in}\;\R^d.
\end{equation}
Suppose that $f$ and $g$ are supported on an annulus of radius $R$:
\begin{equation}\label{S.31}
f(x)=0,\; g(x)=0\quad\text{unless}\;R<|x|\le2R,
\end{equation}
and that $h$ is bounded by some $\kappa$ and supported on $B_{2R}$.
Then we claim
\begin{align}
\expec{|\nabla u(0)|^{2p}}^\frac{1}{2p}&\lesssim\, \sup_{A\in\Omega}\bigg(R^{-d}\int_{\R^d} \big( |g|^{2p} + R^{2p} |f|^{2p} \big)dy
\bigg)^\frac{1}{2p}
+T^{-1}\min\{R,\sqrt{T}\} \sup_{A\in\Omega} \kappa.
\label{S.32}
\\
\langle|\nabla u(0)|  \rangle            &\lesssim\,\expec{R^{-d}\int_{\R^d} \big( |g|^{2} + R^{2} |f|^{2} \big)dy}^\frac{1}{2}+T^{-1}\expec{\kappa^2}.\label{S.33}
\end{align}
To prove \eqref{S.32}, we start by noting that \eqref{S.30} yields the representation formula
\begin{equation}\nonumber
u(x)=\int_{\R^d} G_T(x,y)(\nabla\cdot g+f+T^{-1}h)(y) \;dy,
\end{equation}
which we use in form of
\begin{equation}\nonumber
\nabla u(0)= -\int_{\R^d} \nabla\nabla G_T(0,y) g(y) \;dy + \int_{\R^d} \nabla_x G_T(0,y) (f(y)+T^{-1}h(y)) \;dy.
\end{equation}
By Cauchy-Schwarz' inequality and the support assumption, this yields
\begin{multline*}
|\nabla u(0)|\le\bigg(\int_{R<|y|\le2R}|\nabla\nabla G_T(0,y)|^2 \;dy \int_{R<|y|\le2R} |g(y)|^2 \;dy \bigg)^{\frac{1}{2}}\\
+ \bigg( \int_{R<|y|\le2R} |\nabla_x G_T(0,y)|^2 \;dy \int_{R<|y|\le2R} |f(y)|^2 \;dy \bigg)^\frac{1}{2}
\\
+T^{-1}\kappa \int_{B_{2R}}|\nabla_x G_T(0,y)|dy.
\end{multline*}
This implies by H\"older's inequality in probability 
\begin{align}
\expec{|\nabla u(0)|^{2p}}^\frac{1}{2p}&\le\Lambda_1\sup_{A\in\Omega}\bigg(R^{-d}\int_{\R^d} \big( |g|^{2p}+R^{2p}|f|^{2p} \big)dy
\bigg)^\frac{1}{2p}+\Lambda_2 T^{-1}\sup_{A\in\Omega} \kappa,
\nonumber\\
\expec{|\nabla u(0)|}     &\le\Lambda_3\expec{R^{-d}\int_{\R^d} \big( |g|^{2}+R^{2} |f|^{2} \big)
dy}^\frac{1}{2}+\Lambda_4T^{-1}\expec{\kappa^2} \,\nonumber
\end{align}
where we have set for abbreviation
\begin{eqnarray*}
\Lambda_1&:=&R^d\expec{R^{-d}\int_{R<|y|\le2R} \big( |\nabla\nabla G_T(0,y)|^{2p} + R^{-2p} |\nabla_x G_T(0,y)|^{2p} \big) \;dy
}^\frac{1}{2p},
\\
\Lambda_2&:=&\expec{\Big(\int_{B_{2R}}|\nabla_x G_T(0,y)|dy\Big)^{2p}}^{\frac{1}{2p}}, 
\\
\Lambda_3&:=&R^d\expec{R^{-d}\int_{R<|y|\le2R} \big( |\nabla\nabla G_T(0,y)|^{2} + R^{-2} |\nabla_x G_T(0,y)|^{2} \big) \;dy
}^\frac{1}{2}, 
\\
\Lambda_4&:=&\expec{\Big(\int_{B_{2R}}|\nabla_x G_T(0,y)|dy\Big)^{2}}^{\frac{1}{2}}, 
\end{eqnarray*}
On the one hand, \eqref{S.37} in Step~1 exactly yields $\Lambda_2\lesssim \Lambda_1\lesssim 1$.
On the other hand, a decomposition of $B_{2R}$ into the dyadic annuli $\{2^{i}<|x|\leq 2^{i+1}\}$
for $i \in (-\infty,I)\cap\Z$ with $I=[\log_2 (2R)]+1$ combined
with the triangle inequality, \eqref{S.37}, and H\"older's inequality yields (using the exponential cut-off for $R\ge \sqrt{T}$):
\begin{eqnarray*}
\Lambda_4\,\leq\, \Lambda_2 &\leq & \sum_{i=-\infty}^I \expec{\Big(\int_{2^{i}<|y|\le 2^{i+1}}|\nabla_x G_T(0,y)|dy\Big)^{2p}}^{\frac{1}{2p}}\\
&\lesssim & \sum_{i=-\infty}^I (2^i)^{d(1-\frac{1}{2p})} \expec{\int_{2^{i}<|y|\le 2^{i+1}}|\nabla_x G_T(0,y)|^{2p}dy}^{\frac{1}{2p}}\\
&\stackrel{\eqref{S.37}}{\lesssim} & \sum_{i=-\infty}^I (2^i)^{d(1-\frac{1}{2p})}  \Big((2^i)^{d+2p(1-d)}\Big)^{\frac{1}{2p}}
\exp(-c\frac{2^i}{\sqrt{T}})\\
&=&\sum_{i=-\infty}^I  (2^i)^{1-d} \exp(-c\frac{2^i}{\sqrt{T}})\,\lesssim \, \min\{2^I,\sqrt{T}\}.
\end{eqnarray*}
The desired estimates \eqref{S.32} and \eqref{S.33} follow.

\medskip

\step{3} Consider an $A$-dependent functions $u=u(x;A)$ satisfying
\begin{equation}\label{S.42}
T^{-1}u-\nabla\cdot A\nabla u=0\quad\mbox{in}\;B_{2R}.
\end{equation}
Then we claim
\begin{align}
\expec{|\nabla u(0)|^{2p}}^\frac{1}{2p}&\lesssim\,(1+\frac{R}{\sqrt{T}} )\sup_{A\in\Omega}\bigg(R^{-d}\int_{B_{2R}}|\nabla u|^{2p}dy\bigg)^\frac{1}{2p}.
\label{S.40}
\\
\expec{|\nabla u(0)|}           &\lesssim\,(1+\frac{R}{\sqrt{T}} )\expec{R^{-d}\int_{B_{2R}}|\nabla u|^{2}dy}^\frac{1}{2}.
\label{S.41}
\end{align}
To see this, consider a cut-off function $\eta$ for $B_R$ in $B_{\frac{3}{2}R}$ such that $|\nabla\eta| \lesssim R^{-1}$ and set
$v:=\eta (u-\bar u)$, where $\bar u$ denotes the average of $u$ on $B_{\frac{3}{2}R}$.
Equation~\eqref{S.42} yields
\begin{equation}\label{u_harmonic_f_g}
T^{-1}v- \nabla\cdot A \nabla v = -\nabla\cdot g + f+T^{-1}h
\end{equation}
with $g:=(u-\ext u)A\nabla\eta$, $f:=-\nabla\eta\cdot A\nabla u$, and $h=\eta \bar u$.
By choice of $\eta$, the functions $g$ and $f$ satisfy the support condition \eqref{S.31} and we have
for all $q\ge 1$
\begin{eqnarray*}
\int_{\R^d} \big( |g|^{2q} + R^{2q} |f|^{2q} \big) dy &\lesssim &\int_{B_{\frac{3}{2}R}} \big( R^{-2q} |u-\bar u|^{2q} + |\nabla u|^{2q}
\big)dy,
\end{eqnarray*}
so that Poincar\'{e}'s inequality on $B_{\frac{3}{2}R}$ applied to the first term of the RHS yields
\begin{equation}\label{S.91}
 \int_{\R^d} \big( |g|^{2q} + R^{2q} |f|^{2q} \big)dy \, \lesssim \,\int_{B_{\frac{3}{2}R}}|\nabla u|^{2q}dy.
\end{equation}
It remains to bound the second RHS term of \eqref{S.32}. To this aim we now take $\eta$ a cut-off function for $B_{\frac{3}{2}R}$ in $B_{2R}$ such that $|\nabla\eta| \lesssim R^{-1}$. Testing \eqref{S.42} with $\eta^2u$ and integrating on $B_{2R}$ then yields
$$
T^{-1}\int_{B_{2R}}\eta^2u^2dy\,=\,-\int_{B_{2R}} \eta^2 \nabla u\cdot A\nabla udy-2\int_{B_{2R}}\eta u \nabla \eta\cdot A \nabla u.
$$
We absorb the second RHS term into the LHS by Young's inequality and get by definition of $\eta$
$$
|\bar{u}|\,=\, \Big|\fint_{B_{\frac{3}{2}R}} udy\Big|\,
\leq \,\Big(\fint_{B_{\frac{3}{2}R}}u^2dy\Big)^{\frac{1}{2}} \,\lesssim \,\sqrt{T} \Big(\fint_{B_{2R}} |\nabla u|^2dy\Big)^{\frac{1}{2}},
$$
so that by Jensen's inequality
\begin{equation}\label{S.92}
T^{-1} \min\{\sqrt{T},R\}  \{\sup_{B_{2R}}|\eta \bar u|\} \,\lesssim \, \frac{R}{\sqrt{T}} 
\bigg(R^{-d}\int_{B_{2R}}|\nabla u|^{2q}dy\bigg)^\frac{1}{2q}.
\end{equation}
By \eqref{S.91} and \eqref{S.92} for $q=p$ and for $q=2$, 
\eqref{S.40} and \eqref{S.41} 
follow from \eqref{S.32}
 and \eqref{S.33}.

\medskip

\step{4} Proof of \eqref{eq:annealed-estim} and \eqref{eq:annealed-estim2}.

\noindent We fix $y\in\mathbb{R}^d\setminus \{0\}$ and apply Step~3 to
$u(x)=G_T(x,y)$ and $R=\frac{1}{6}|y|$.
From \eqref{S.40} we obtain
\begin{equation}\nonumber
\expec{|\nabla_x G_T(0,y)|^{2p}}^\frac{1}{2p}\,\lesssim\, (1+\frac{|y|}{\sqrt{T}})\sup_{A\in\Omega}\bigg(|y|^{-d}\int_{B_{\frac{1}{3}|y|}}
|\nabla_x G_T(x,y)|^{2p} \;dx \bigg)^\frac{1}{2p}.
\end{equation}
Since
\begin{equation}\label{ball_ann}
 B_{\frac{1}{3}|y|}\subset\Big\{x\in\R^d:\frac{2}{3}|y|<|x-y|\le\frac{4}{3}|y|\Big\},
\end{equation}
we obtain by \eqref{eq:int-grad} with $R=\frac{1}{2}|y|$ the desired estimate \eqref{eq:annealed-estim}, i.e.\ we have that
\begin{equation}\nonumber
\expec{|\nabla_x G_T(0,y)|^{2p}}^\frac{1}{2p}\lesssim (1+\frac{|y|}{\sqrt{T}})|y|^{1-d}\exp(-c\frac{|y|}{\sqrt{T}})
\,\lesssim\, |y|^{1-d}\exp(-c\frac{|y|}{\sqrt{T}}),
\end{equation}
for a slightly smaller $c>0$.

\medskip

\noindent Next we turn to the mixed second gradient, and apply Step~3 to the function $u(x)=\nabla_y G_T(x,y)$ with $R=\frac{1}{6}|y|$ and obtain from \eqref{S.41} that
\begin{equation}\nonumber
\langle|\nabla\nabla G_T(0,y)|\rangle\,\lesssim
(1+\frac{|y|}{\sqrt{T}})\bigg\langle\, |y|^{-d}\int_{B_{\frac{1}{3}|y|}}|\nabla\nabla G_T(x,y)|^2 \;dx\bigg\rangle^\frac{1}{2}.
\end{equation}
The inclusion \eqref{ball_ann} yields
\begin{equation}\nonumber
\langle|\nabla\nabla G_T(0,y)|\rangle\,\lesssim\,
(1+\frac{|y|}{\sqrt{T}})\bigg\langle|y|^{-d}\int_{\frac{2}{3}|y|\le|x-y|\le\frac{4}{3}|y|}|\nabla\nabla
G_T(x,y)|^2 \;dx \bigg\rangle^\frac{1}{2}.
\end{equation}
By stationarity in form of 
$\expec{|\nabla \nabla G_T(x,y)|^2}=\expec{|\nabla \nabla G_T(0,y-x)|^2}$ and \eqref{S.37}, this yields the desired estimate \eqref{eq:annealed-estim2}.


\subsection{Proof of Lemma~\ref{lem:modif-ener}}

We split the proof into three steps. We start with the proof of \eqref{eq:modif-ener-1.1},
then show that it implies \eqref{eq:modif-ener-1.2} by a dyadic decomposition of the RHS, and
then turn to \eqref{eq:modif-ener-1.3}, which is a variation of \eqref{eq:modif-ener-1.1}.

\medskip

\step{1} Proof of \eqref{eq:modif-ener-1.1}.

\noindent 
By rescaling length according to $x=\sqrt{T}\hat x$, we see that it is enough
to show that \eqref{eq:modif-ener-3.1} yields \eqref{eq:modif-ener-1.1} for $T=1$. By dyadic iteration,
it is enough to show there exists a constant $\theta(d,\lambda)<1$ such that
\begin{equation}\label{fu.10}
v-\nabla\cdot A\nabla v=0\quad\mbox{in}\;B_{2R}
\end{equation}
implies
\begin{equation}\nonumber
\int_{B_R}(v^2+|\nabla v|^2)dx\le \theta
\int_{B_{2R}}(v^2+|\nabla v|^2)dx,
\end{equation}
which by the Widman hole-filling trick follows from
\begin{equation}\label{fu.9}
\int_{B_R}(v^2+|\nabla v|^2)dx\lesssim
\int_{R< |x|\le 2R}(v^2+|\nabla v|^2)dx.
\end{equation}
In order to obtain \eqref{fu.9}, we test \eqref{fu.10} with $\eta^2(v-\bar v)$,
where $\eta$ is a cut-off function for $B_R$ in $B_{2R}$
and $\bar v$ is the average of $v$ on $\{R< |x|\le 2R\}$, to the effect of
\begin{equation}\nonumber
\int_{B_{2R}}\big(\eta^2(v-\bar v)v+\nabla(\eta^2(v-\bar v))\cdot A\nabla v\big)dx=0.
\end{equation}
For the massive term we use $v(v-\bar v)\ge\frac{1}{2}v^2-\frac{1}{2}\bar v^2$;
for the elliptic term we use $\nabla(\eta^2(v-\bar v))\cdot A\nabla v
=\eta^2\nabla v\cdot A\nabla v+2\eta(v-\bar v)\nabla\eta\cdot A\nabla v
\ge\lambda\eta^2|\nabla v|^2-2\eta|v-\bar v||\nabla\eta||\nabla v|
\ge\frac{1}{2}\lambda\eta|\nabla v|^2-\frac{2}{\lambda}(v-\bar v)^2|\nabla\eta|^2$, so that
we obtain
\begin{equation}\label{cfu.2}
\int_{B_{2R}}\eta^2(v^2+\lambda|\nabla v|^2)dx\le\int_{B_{2R}}\big(\eta^2\bar v^2+\frac{4}{\lambda}|\nabla\eta|^2(v-\bar v)^2\big)dx.
\end{equation}
Using the properties of $\eta$, this yields
\begin{equation}\nonumber
\int_{B_R}(v^2+|\nabla v|^2)dx\lesssim R^d\bar v^2+R^{-2}\int_{R<|x|\le 2R}(v-\bar v)^2dx.
\end{equation}
In order to obtain \eqref{fu.9}, we use Jensen's inequality on the average $\bar v$
yielding
$R^d\bar v^2\lesssim \int_{R<|x|\le 2R}v^2dx$; we use Poincar\'e's estimate on
the annulus yielding $R^{-2}\int_{R<|x|\le 2R}(v-\bar v)^2dx\lesssim\int_{R<|x|\le 2R}|\nabla v|^2dx$.

\medskip

\step{2} Proof of \eqref{eq:modif-ener-1.2}.

\noindent Rescaling lengths according to $x=R\hat x$ (which entails
$v=R\hat v$ and $\sqrt{T}=R\sqrt{\hat T}$) we see that it is enough to establish
\eqref{eq:modif-ener-1.2} for $R=1$ only, that is,
\begin{equation}\label{fu.2}
\left(\int_{B_1}(T^{-1}v^2+|\nabla v|^2)dx\right)^\frac{1}{2}\lesssim 
\left(\int_{\R^d}((|x|+1)^{-\alpha}|g|)^2dx\right)^\frac{1}{2}.
\end{equation}
By the triangle inequality, it is sufficient to establish \eqref{fu.2} 
under the additional condition that
\begin{equation}\label{fu.5}
\mathrm{supp}\,g\subset\{R<|x|\le 2R\},
\end{equation}
in which case \eqref{fu.2} turns into
\begin{equation}\label{fu.4}
\left(\int_{B_1}(T^{-1}v^2+|\nabla v|^2)dx\right)^\frac{1}{2}\lesssim 
(R+1)^{-\alpha}\left(\int_{\R^d}|g|^2dx\right)^\frac{1}{2}.
\end{equation}
By the energy estimate, i.\ e.\ testing \eqref{eq:modif-ener-3.2} by $v$, we have
\begin{equation}\label{fu.6}
\left(\int_{\R^d}(T^{-1}v^2+|\nabla v|^2)dx\right)^\frac{1}{2}\lesssim 
\left(\int_{\R^d}|g|^2dx\right)^\frac{1}{2}.
\end{equation}
This trivially yields \eqref{fu.4} for $R\le 1$, so that we may focus on $R\ge 1$.
For $R\ge 1$, we see that \eqref{fu.6} implies \eqref{fu.4}
using \eqref{eq:modif-ener-1.1} since 
\begin{equation*}\label{fu.8}
T^{-1}v-\nabla\cdot A\nabla v=0\quad\mbox{in}\;B_R.
\end{equation*}

\medskip

\step{3} Proof of \eqref{eq:modif-ener-1.3}.

\noindent 
As in Step~1, by rescaling length according to $x=\sqrt{T}\hat x$,
we may restrict to the case of $T=1$.
By dyadic iteration,
it is enough to show there exists a constant $\theta(d,\lambda)<1$ such that
\begin{equation}\label{v.10}
u-\nabla\cdot A\nabla u=f:=-\xi\cdot x\quad\mbox{in}\;B_{2R}
\end{equation}
implies
\begin{equation}\nonumber
\int_{B_R}(u^2+|\nabla u|^2)dx\,\le\, \theta
\left(\int_{B_{2R}}(u^2+|\nabla u|^2)dx+R^{d+2}\right),
\end{equation}
which follows by the Widman hole-filling trick from
\begin{equation}\label{v.9}
\int_{B_R}(u^2+|\nabla u|^2)dx\lesssim
\int_{R< |x|\le 2R}(u^2+|\nabla u|^2)dx+R^{d+2}.
\end{equation}
In order to obtain \eqref{v.9}, we test (\ref{v.10}) with $\eta^2(u-\bar u)$,
where $\eta$ is a cut-off function for $B_R$ in $B_{2R}$
and $\bar u$ is the average of $u$ on the annulus $\{R< |x|\le 2R\}$, to the effect of
\begin{equation}\nonumber
\int_{B_{2R}}\eta^2(\frac{1}{2}u^2+\frac{1}{2}(u-\bar u)^2+\nabla u\cdot A\nabla u)dx
=\int_{B_{2R}}\big(\eta^2\frac{1}{2}\bar u^2-2\eta(u-\bar u)\nabla\eta\cdot A\nabla u+\eta^2(u-\bar u)f\big)dx.
\end{equation}
By the assumptions on $A$, this turns into the inequality
\begin{equation}\nonumber
\int_{B_{2R}}\eta^2(\frac{1}{2}u^2+\frac{1}{2}(u-\bar u)^2+\lambda|\nabla u|^2)dx\,
\le\,\int_{B_{2R}}\big(\eta^2\frac{1}{2}\bar u^2+2\eta|u-\bar u||\nabla\eta||\nabla u|+\eta^2(u-\bar u)f\big)dx.
\end{equation}
Using Young's inequality, this implies the estimate
\begin{equation}\nonumber
\int_{B_{2R}}\eta^2(u^2+|\nabla u|^2)dx\,
\lesssim\,\int_{B_{2R}}\big(\eta^2\bar u^2+|\nabla\eta|^2|u-\bar u|^2+\eta^2f^2\big)dx.
\end{equation}
Using the properties of $\eta$, this yields
\begin{equation}\nonumber
\int_{B_R}(u^2+|\nabla u|^2)dx\,
\lesssim \, R^d\bar u^2+\int_{B_{2R}}f^2dx+R^{-2}\int_{R<|x|\le 2R}|u-\bar u|^2dx.
\end{equation}
We now use Jensen's inequality on the average $\bar u$
yielding
$R^d\bar u^2\lesssim \int_{R<|x|\le 2R}u^2$; we also appeal to Poincar\'e's estimate on
the annulus yielding $R^{-2}\int_{R<|x|\le 2R}(u-\bar u)^2dx\lesssim\int_{R<|x|\le 2R}|\nabla u|^2dx$.
This entails
\begin{equation}\nonumber
\int_{B_R}(u^2+|\nabla u|^2)dx
\le\int_{B_{2R}}f^2dx+\int_{R<|x|\le 2R}(u^2+|\nabla u|^2)dx.
\end{equation}
Appealing to the special form of $f$ yields \eqref{v.9}.


\subsection{Proof of Lemma~\ref{lem:ptwise-estim}}

On $\{|z'|>\frac{|z|}{2}\}$ the Green function satisfies
\begin{equation}\label{eq:equation-G_T-harmonic}
T^{-1}G_T(\cdot,0)-\nabla \cdot A\nabla G_T(\cdot,0)\,=\,0.
\end{equation}
Let $1\sim R\leq \frac{|z|}{6}$. 
We first prove the result for $d>2$ by combining Caccioppoli's inequality and the pointwise bounds on the Green functions, and then turn to $d=2$ 
using in addition Lemma~\ref{lem:modif-ener}.
Caccioppoli's inequality for \eqref{eq:equation-G_T-harmonic} then yields
\begin{equation}\label{eq:Caccio-Green}
\int_{B_{R}(z)} |\nabla_{z'} G_T(z',0)|^2dz' \,\lesssim \, \int_{B_{\frac{3R}{2}}(z)} G_T^2(z',0)dz',
\end{equation}
and we conclude by the pointwise estimate \eqref{eq:ptwise-decay-estim} for $d>2$.

\medskip

\noindent For $d=2$ we appeal to Lemma~\ref{lem:modif-ener} and use \eqref{eq:modif-ener-1.1}  in the form of:
\begin{equation*}
\left(\int_{B_R(z)}|\nabla_{z'} G_T(z',0)|^2dz' \right)^\frac{1}{2}\lesssim 
|z|^{-\alpha}\left(\int_{B_{\frac{|z|}{3}}(z)}(T^{-1}G_T^2(z',0)+|\nabla G_T(z',0)|^2)dz'\right)^\frac{1}{2}.
\end{equation*}
On the one hand, the pointwise estimate  \eqref{eq:ptwise-decay-estim} for $d=2$ yields for
the first RHS term since $|z|\gtrsim 1$
$$
\int_{B_{\frac{|z|}{3}}(z)}T^{-1}G_T^2(z',0)dz'\,\lesssim \,  \sup_{B_{\frac{|z|}{3}}(z)}\Big(\frac{|z'|}{\sqrt{T}}\Big)^2  \exp(-c\frac{|z'|}{\sqrt{T}})
\ln^2(2+\frac{\sqrt{T}}{|z'|})\,\lesssim \, 1.
$$
On the other hand, for the second RHS term we use Caccioppoli's inequality in the form: For all $c\in \R$,
\begin{equation*}
\int_{B_{\frac{|z|}{3}}(z)}|\nabla G_T(z',0)|^2dz'\,\lesssim \, |z|^{-2} \int_{B_{\frac{|z|}{2}}(z)}(G_T(z',0)-c)^2dz'
+T^{-1}|z|^2|c|.
\end{equation*}
If $|z|\leq \sqrt{T}$, we choose $c=\fint_{B_{\frac{3|z|}{2}}}G_T(z',0)dz'$ and appeal to the oscillation estimate \eqref{g.20} (in its $T$-rescaled version) to the effect of
$$
|z|^{-2} \int_{B_{\frac{|z|}{2}}(z)}(G_T(z',0)-c)^2dz'\,\leq \, |z|^{-2} \inf_{\kappa \in \R} \int_{B_{\frac{3|z|}{2}}(0)}(G_T(z',0)-\kappa)^2dz'\,\lesssim \,1
$$
and to \eqref{eq:ptwise-decay-estim} which implies that
$$
T^{-1}|z|^2|c|\,\lesssim \, \sup_{B_{\frac{3|z|}{2}}}\Big\{\Big(\frac{|z'|}{\sqrt{T}}\Big)^2  \exp(-c\frac{|z'|}{\sqrt{T}})
\ln(2+\frac{\sqrt{T}}{|z'|})\Big\}\,\lesssim \, 1.
$$
If $|z|> \sqrt{T}$, we take $c=0$ and use that $\sup_{B_{\frac{3|z|}{2}}}G_T(z',0)\lesssim 1$ by \eqref{eq:ptwise-decay-estim}.

\section*{Acknowledgements}
The first author acknowledges financial support from the European Research Council under
the European Community's Seventh Framework Programme (FP7/2014-2019 Grant Agreement
QUANTHOM 335410).

\bibliographystyle{plain}


\end{document}